\documentclass[11pt]{amsart}
\usepackage{hyperref}

\usepackage{stmaryrd}

\usepackage{enumerate}

\usepackage[a4paper, twoside=false, vmargin={2cm,3cm}, includehead]{geometry}

\newtheorem{theorem}{Theorem}[section]
\newtheorem*{theorem*}{Theorem}
\newtheorem{lemma}[theorem]{Lemma}
\newtheorem{corollary}[theorem]{Corollary}
\newtheorem{proposition}[theorem]{Proposition}
\newtheorem*{vdcLemma}{van der Corput Lemma}

\theoremstyle{definition}
\newtheorem*{definition}{Definition}
\newtheorem*{notation}{Notation}

\theoremstyle{remark}

\newtheorem*{remark}{Remark}
\newtheorem*{remarks}{Remarks}
\newtheorem{example}{Example}[section]

\newcommand{\baton}[1]{{\mathbb #1}}
\newcommand\C{{\baton C}}
\newcommand\D{{\baton D}}
\newcommand\E{{\baton E}}
\newcommand{\N}{\baton N}
\newcommand{\PP}{{\baton P}}
\newcommand\Q{{\baton Q}}
\newcommand\R{{\baton R}}
\newcommand\T{{\baton T}}
\newcommand\Z{{\baton Z}}

\newcommand{\CA}{{\mathcal A}}

\newcommand{\CF}{{\mathcal F}}
\newcommand{\CH}{{\mathcal H}}
\newcommand{\CI}{{\mathcal I}}
\newcommand{\CK}{{\mathcal K}}

\newcommand{\CN}{{\mathcal N}}

\newcommand{\CS}{{\mathcal S}}
\newcommand{\CX}{{\mathcal X }}
\newcommand{\CY}{{\mathcal Y }}
\newcommand{\CZ}{{\mathcal Z }}

\newcommand{\bc}{\mathbf{c}}
\newcommand{\bh}{\mathbf{h}}
\newcommand{\bI}{{\mathbf{I}}}
\newcommand{\bJ}{{\mathbf{J}}}
\newcommand{\bk}{\mathbf{k}}
\newcommand{\bm}{\mathbf{m}}
\newcommand{\bn}{\mathbf{n}}

\newcommand{\bt}{\mathbf{t}}

\newcommand{\bzero}{{\mathbf{0}}}
\newcommand{\one}{{\mathbf{1}}}

\renewcommand{\vec}{\overrightarrow}

\newcommand{\uepsilon}{{\underline \epsilon}}
\newcommand{\ubh}{\underline{\mathbf{h}}}
\newcommand{\conj}{{\mathcal{C}}}
\newcommand{\ve}{\varepsilon}
\newcommand{\wt}{\widetilde}
\newcommand{\wh}{\widehat}
\newcommand{\inv}{^{-1}}
\DeclareMathOperator{\cor}{Corr}
\newcommand{\e}{\mathrm{e}}
\newcommand{\id}{\mathrm{Id}}
\newcommand{\norm}[1]{\lVert #1\rVert}
\newcommand{\nnorm}[1]{|\!|\!| #1 |\!|\!|}

\newcommand{\Cube}[1]{{\llbracket #1\rrbracket}}
\newcommand{\Cub}[1]{{\llbracket {#1}^*\rrbracket}}
\newcommand{\st}{{\text{\rm st}}}
\newcommand{\er}{{\text{\rm er}}}

\newcommand{\limav}{\mathrm{lim\, 	Av\,}}
\newcommand{\Limav}[1]{\mathrm{lim\, Av}_{#1}\,}
\newcommand{\Av}{\mathrm{Av}}
\newcommand{\Limsup}{\mathrm{limsup}}
\newcommand{\Lim}{\mathrm{lim}}

\newcounter{mycounter}
\numberwithin{equation}{section}

\begin{document}

\title[Weighted multiple ergodic averages]{Weighted multiple ergodic averages and correlation sequences}
\author{Nikos Frantzikinakis}
\address[Nikos Frantzikinakis]{University of Crete, Department of mathematics, Voutes University Campus, Heraklion 71003, Greece} \email{frantzikinakis@gmail.com}
\author{Bernard Host}

\address[Bernard Host]{
Universit\'e Paris-Est Marne-la-Vall\'ee, Laboratoire d'analyse et
de math\'e\-matiques appliqu\'ees, UMR CNRS 8050, 5 Bd Descartes,
77454 Marne la Vall\'ee Cedex, France }
\email{bernard.host@u-pem.fr}

\thanks{The second author was partially supported by Centro de Modelamiento Matem\'atico, Universitad de Chile.}

\begin{abstract}
We study mean convergence results for weighted multiple ergodic averages
defined by commuting transformations with iterates given by integer
polynomials in several variables. Roughly speaking, we prove that   a bounded sequence is a good universal weight for mean convergence of such averages if and only if the averages of this sequence times any nilsequence converge.
Key role in the proof play two decomposition results of independent interest.
 The first states that every  bounded sequence in several variables  satisfying some regularity conditions  is a sum of a nilsequence and a sequence that has small uniformity norm (this generalizes a result of the second author and B.~Kra);
 and the second states that every multiple correlation sequence in several variables is a sum of a nilsequence and a sequence that is small in uniform density (this generalizes a result of the first author). Furthermore, we use the previous results in order to  establish mean convergence and recurrence results for a variety of sequences of dynamical and arithmetic origin  and give some combinatorial implications.
\end{abstract}



\maketitle
\setcounter{tocdepth}{1}
\tableofcontents

\section{Introduction}
Since the early 80's a lot of effort has been put  in the  study of  the limiting behavior  of multiple ergodic averages.
 This study  was partly motivated by combinatorial implications, since positiveness properties of
 such averages
 imply various far reaching extensions of the
celebrated theorem of Szemer\'edi on arithmetic progressions.
After a long series of partial results, most notably those in \cite{Au09, CL84, CL88, FW96, H09, HK05, HK05b, L05c, R95, T08, Z07}, M. Walsh~\cite{W12}, building on previous work of T.~Tao~\cite{T08},
proved the following mean convergence result:
\begin{theorem}[\mbox{\cite{W12}}\footnote{The argument in \cite{W12} is given for Ces\`aro averages and  $d=1$  but the same proof works
in this more general case (see
\cite{Z15a} for details).}]\label{T:Walsh}
Let $d,\ell, s\in \N$, $(X,\CX,\mu)$ be a probability space, and  $T_1,\ldots, T_\ell\colon X\to X$ be invertible commuting measure preserving transformations. Then for every F\o lner sequence $(I_k)_{k\in \N}$ of subsets of $\N^d$, polynomials $p_{i,j}\colon \N^d\to \Z$, $i=1,\ldots, \ell$, $j=1,\ldots, s$,  and functions $f_1,\ldots, f_s\in L^\infty(\mu)$, the averages
\begin{equation}\label{E:averages1}
\frac{1}{|I_k|}\sum_{\bn\in I_k }f_1(\prod_{i=1}^\ell T_i^{p_{i,1}(\bn)}x)
\cdot \ldots\cdot f_s(\prod_{i=1}^\ell T_i^{p_{i,s}(\bn)}x),
\end{equation}
converge in  $L^2(\mu)$ as $k\to+\infty$.
\end{theorem}
\begin{remark} In \cite{W12} the previous result was established under the  weaker hypothesis that the transformations $T_1,\ldots, T_\ell$ generate  a nilpotent group. We believe that our results can be extended to this more general setup but we restrict in this article  to the case of commuting transformations.
\end{remark}

One of the main purposes of this article is to study mean  convergence for weighted versions of the averages \eqref{E:averages1}, that is,  averages of the form
\begin{equation}\label{E:averages2}
\frac{1}{|I_k|}\sum_{\bn\in I_k }w(\bn)\, f_1(\prod_{i=1}^\ell T_i^{p_{i,1}(\bn)}x)
\cdot \ldots\cdot f_s(\prod_{i=1}^\ell T_i^{p_{i,s}(\bn)}x),
\end{equation}
 where
$w\colon \N^d\to \C$ is a bounded sequence.
 We call a sequence $w$ for which the previous averages converge for all choices of systems, functions, and polynomials, a \emph{good universal weight for mean convergence} of the averages \eqref{E:averages2}.

When  $d=1$,
examples of good universal weights
 for some multiple ergodic averages
can be found in \cite{ADM15, AM15, AM15a, AM15b, C09, E13, F15, FH15b,  HK09, Z15b}. Most of  these results
deal with the case where $\ell=1$ and are based on the
theory of characteristic factors that was pioneered by H.~Furstenberg. They  depend crucially on the work of B.~Host and B.~Kra \cite{HK05} and subsequent developments in \cite{HK05b, L05c},  which, in the case where all the transformations are equal, gives a characterization in terms of nilsystems  of the smallest factor of the system that controls the limiting behavior of the averages  \eqref{E:averages2}. Unfortunately, no such characterization is known in the
 case of general commuting transformations (but see \cite{Au11a, Au11b} for related progress), which is the reason why this method is not applicable for our more general setup. Moreover, the method used by M.~Walsh in \cite{W12} does not seem applicable to
  weighted averages
and
  no general  criterion suitable for checking mean convergence of averages of  the form \eqref{E:averages2}    is known.
We fill this gap by showing  in  Theorem~\ref{th:weighted-L2} that a bounded sequence $w\colon \N^d\to \C$
is a good universal weight for mean convergence of the averages \eqref{E:averages2} if and only if  the averages
\begin{equation}\label{E:condition}
\frac{1}{|I_k|}\sum_{\bn\in I_k}w(\bn)\cdot \psi(\bn)
\end{equation}
converge for every  nilsequence $\psi$ in  $d$ variables and F\o lner sequence $(I_k)_{k\in\N}$ in $\N^d$.
Furthermore, when one replaces throughout the averages $\frac{1}{I_k}\sum_{\bn\in I_k }$
  by the Ces\`aro averages $\frac{1}{N^d}\sum_{\bn\in [1,N]^d}$, we prove in Theorem~\ref{th:Cesaro} that a similar criterion holds for weak convergence,
 and a condition somewhat stronger than \eqref{E:condition}
 suffices for mean convergence.   Even for single variable sequences the mean convergence criterion is new and its proof (strangely) depends on decomposition results for sequences in two variables.
 Prior to this work, only the case $d=\ell=1$ was treated  in \cite{C09} for mean convergence, while
   for weak convergence the case where  $d=1$ and $\ell\in \N$ is arbitrary was treated in
 \cite{F15}.

 To prove these results we  use the mean convergence result of
M.~Walsh  as  a black box and
  two decomposition results of independent interest. These are
  Theorems~\ref{th:regular}  and \ref{th:regular+antiunif} which extend
  similar results for single variable sequences from
   \cite[Theorem 2.19]{HK09} and \cite[Theorem~1.2]{F15}.
 Roughly speaking, they state the following:

\begin{enumerate}
\item
If the averages  \eqref{E:condition} converge for every nilsequence $\psi\colon \N^d\to \C$ and  every F\o lner sequence $(I_k)_{k\in \N}$ in $\N^d$,
then the sequence $w\in \ell^\infty(\N^d)$    is
a sum of a nilsequence and a sequence that has small uniformity norm.

\item
Any sequence of the form
$$
\int \prod_{i=1}^\ell f_1(T_i^{p_{i,1}(\bn)}x)
\cdot \ldots\cdot \prod_{i=1}^\ell f_s( T_i^{p_{i,s}(\bn)}x)\, d\mu, \quad \bn \in \N^d,
$$
is the sum of a  nilsequence and a sequence that is small in uniform density.
\end{enumerate}

Regarding the second decomposition,  Theorem~\ref{th:correl-lin} gives more precise information when the iterates are linear, it implies for example that the sequences

\begin{equation}\label{E:three}
\int f\cdot T_1^{m}f\cdot T_2^{n}f\cdot T_3^{r}f\, d\mu, \ \int f\cdot T_1^{m}f\cdot T_2^{n}f\cdot T_3^{m+n}f\, d\mu, \
 \int f\cdot T_1^{n}f\cdot T_2^{n}f\cdot T_3^{n
 }f\, d\mu,
 \end{equation}
are  $1$-step, $2$-step,  and  $3$-step nilsequences respectively modulo small errors in uniform density (simple examples show that the degree of nilpotency is optimal).

Using the previous  criteria we  prove mean convergence results for weighted ergodic averages
with   weights
given by various sequences of dynamical origin,
 bounded multiplicative functions, generalized polynomials, and
 Hardy field sequences (see Sections~\ref{SS:2.2}, \ref{SS:2.3}, \ref{SS:arithmetic}).
 We deduce  some  multiple recurrence results  and combinatorial consequences;
showing for example that every set of integers with positive upper density contains arbitrarily long arithmetic progressions with common
difference of the form $m^2+n^2$ where $m,n $ have   an odd (or an even) number of distinct prime factors
 (see Theorems~\ref{T:convarith}, \ref{T:recarith},  \ref{T:combinatorics}) or $m,n$ are taken from the set  $\{k\in \N\colon \norm{k^a}\in [1/2,3/4]\}$ where $a$ is any positive non-integer (see Theorems~\ref{T:hardymain}, \ref{C:hardyrec}). We also establish  multidimensional variants of these results regarding patterns in positive density subsets of $\Z^\ell$.

In the next section, we give the precise formulation of our main results and define some
 concepts used throughout the article.

\section{Precise statement of main  results}

\subsection{Notation and definitions}
We first introduce some notation that is going to facilitate our presentation.
\subsubsection{Ergodic theory}
Following for example~\cite{G03} we say that a probability  space $(X,\CX,\mu)$ is a \emph{Lebesgue space,} if $X$  can be given the structure of a Polish space  (i.e. metrizable, separable, complete) such that $\CX$  is its Borel $\sigma$-algebra. Throughout the article we make the standard assumption that \emph{all probability  spaces considered are Lebesgue.}

By a system $(X,\CX,\mu,T_1,\dots,T_\ell)$ we mean a Lebesgue probability space $(X,\CX,\mu)$ endowed with $\ell$ invertible commuting measure preserving transformations. For $\vec n=(n_1,\dots,n_\ell)\in\Z^\ell$ we write $T_{\vec n}=T_1^{n_1}\cdot\ldots\cdot T_\ell^{n_\ell}$.
Sometimes we denote by $\vec T$
 the action of $\Z^\ell$ on $X$ and write the system as $(X,\CX,\mu,\vec T)$.
In the sequel, we generally omit the $\sigma$-algebra $\CX$ from our notation.

\subsubsection{Nilmanifolds  and nilsequences}\label{SS:nil}
Let $s\in\N$, $G$ be an $s$-step nilpotent Lie group, and $\Gamma$ be a discrete cocompact subgroup of $G$. Then the quotient space $X=G/\Gamma$ is called an \emph{$s$-step nilmanifold}.
We prefer to denote the  elements of $X$ as points $x,y,\dots$, not  as cosets. The point $e_X$ is the image in $X$ of the unit element of $G$. The natural action of $G$ on $X$ is written $(g,x)\mapsto g\cdot x$ and the unique measure on $X$ invariant under this action is called the \emph{Haar measure} of $X$ and is denoted by  $m_X$.

 Let $\tau_1,\dots,\tau_d$ be  commuting elements of $G$. For $i=1,\dots,d$ let $T_i$ be the translation $x\mapsto\tau_i\cdot x$ by $\tau_i$ on $X$. Then the system  $(X,m_X,T_1,\dots,T_d)$ is called an \emph{$s$-step nilsystem}.
Nilsystems have been extensively studied and basic properties were established by  Auslander~\cite{AGH63}, Parry~\cite{P69,P70}, Lesigne~\cite{Le91},  Leibman~\cite{L05,L05b}, and others.

\begin{definition}[\cite{BHK05}]
If $X=G/\Gamma$ is an $s$-step nilmanifold, $\Psi\in C(X)$,  and $\tau_1,\ldots, \tau_d\in G$ are commuting elements, then the  sequence
$\big(\Psi(\tau_1^{n_1}\cdot\ldots\cdot\tau_d^{n_d}\cdot e_X)\big)_{n_1,\ldots, n_d \in \N}$ is called an \emph{$s$-step nilsequence in $d$ variables}. Also, for notational convenience, we   define a \emph{$0$-step nilsequence}  to be a constant sequence.
\end{definition}
\begin{remarks}
$\bullet$ In  \cite{BHK05} the notion ``\emph{basic} $s$-step nilsequence'' is used for what we call here an ``$s$-step nilsequence''.

$\bullet$   By \cite[Paragraph 1.11]{L05}  the nilmanifold $X$ is isomorphic to a sub-nilmanifold of a nilmanifold $\tilde{X}=\tilde{G}/\tilde{\Gamma}$ where $\tilde{G}$ is a connected and simply connected $s$-step nilpotent Lie group and all elements of $G$ are represented in $\tilde{G}$.
Hence, whenever needed, we can   assume that the group $G$ is connected and simply connected.
\end{remarks}
In recent years, nilsequences have played a key role in
 ergodic theory and additive combinatorics. They form the right substitute for
 linear exponential sequences needed to formalize certain inverse theorems which are used in the course of studying  various multilinear expressions in analysis and number theory.

\subsubsection{F\o lner sequences and related averages}
We recall first some notions and introduce some notation.
 \begin{notation}
We write $[N]$ for the interval $\{1,2,\dots,N\}$ in $\N$.
\end{notation}
\begin{definition}
A \emph{F\o lner sequence in} $\N^d$
 is a sequence $\bI=(I_j)_{j\in \N}$ of finite subsets of $\N^d$
 that satisfies
$$
\lim_{j\to \infty} \frac{|(I_j+\bk)\triangle I_j|}{|I_j|}= 0\quad \text{for every }\bk\in\Z^d,
$$
where $\triangle$ denotes the symmetric difference and
$I_j+\bk:=\{\bn+\bk\colon\bn\in I_j\}$.
\end{definition}
An example of a F\o lner sequence in $\N$ is a sequence of intervals whose lengths tend to infinity. If $N_j\to+\infty$ and $(\bk_j)_{j\in \N}$ is a sequence in $\N^d$, then $I_j:=\bk_j+[N_j]^d$, $j\in \N$, is a F\o lner sequence in $\N^d$.
In all subsequent results and proofs  we can replace the  general F\o lner sequences by these particular examples.

 If  $a\colon \N^d\to \C$ is a sequence  and ${\bf I}=(I_j)_{j\in \N}$ is a F\o lner sequence in $\N^d$,  we let $$
\Limav{\bI} a(\bn):=\lim_{j\to+\infty}\frac 1{|I_j|}\sum_{\bn\in I_j}a(\bn),
$$
assuming of course that the previous limit exists. If the previous limit exists for every F\o lner sequence ${\bf I}$, then it is independent of ${\bf I}$;
we denote its common value with
$$
\limav a(\bn)
$$
and say that  \emph{the averages of $a$ converge}.
When it is unclear with respect to which variable we take the averages we use the notation
$$
\Limav{\bn,\bI} a(\bn), \quad \Limav{\bn} a(\bn).
$$
Furthermore, we use the notation
$$
\Limsup\bigl|\Av_\bI\, a(\bn)\bigr|:=\limsup_{j\to+\infty}\Bigl|\frac 1{|I_j|}\sum_{\bn\in I_j} a(\bn)\Bigr|
$$
and
$$
\Limsup\bigl| \Av\, a(\bn)\bigr|:=\sup_\bI\bigl(\Limsup\bigl|\Av_\bI \,a(\bn)\bigr|\bigr)
$$
where the sup is taken over all  F\o lner sequences $\bI=(I_j)_{j\in \N}$ of subsets of  $\N^d$.

We use similar notation for  limits in $L^2(\mu)$ involving averages of functions $(f_\bn)_{\bn\in\N^d}$ in $L^2(\mu)$  and write
$$
\Limav{\bI} f_{\bn}, \quad \limav f_{\bn}, \quad
\Limsup\bigl\Vert\Av_\bI\, f_\bn\bigr\Vert_{L^2(\mu)},\quad \Limsup\bigl\Vert \Av\, f_{\bn}\bigr\Vert_{L^2(\mu)}
$$
for the corresponding limits, where in the first two cases convergence takes place in $L^2(\mu)$
and in the last two cases we use $L^2(\mu)$ norms in place of the absolute values.

\subsection{Convergence results for uniform averages}\label{SS:2.2}
First we give  convergence criteria for weighted ergodic averages which are  defined using uniform averages, that is, averages
over arbitrary F\o lner sequences in $\N^d$.
\begin{definition}
If $p_i\colon \N^d\to\Z$, $i=1,\ldots, \ell,$ are polynomials, we call the map $\vec p\colon \N^d\to\Z^\ell$ defined by $\vec p :=(p_1,\dots,p_\ell)$ a \emph{polynomial mapping} from $\N^d$ to $\Z^\ell$. The \emph{degree} $\mathrm{deg}(\vec p)$ of $\vec p$ is $\max_{i=1,\ldots, \ell}(\mathrm{deg}(p_i))$.
\end{definition}
We first state a strengthening of Theorem~\ref{T:Walsh} that will be used frequently in this article.
It is proved in Section~\ref{subsec:Proof-Prop1}.
\begin{proposition}
\label{prop:weight-nil}
For every $d,\ell,s\in\N$, polynomial mappings  $\vec{p_i}\colon \N^d\to\Z^\ell$, $i=1,\ldots, s,$ nilsequence
  $\psi\colon \N^d\to \C$, system $(X,\mu,T_1,\dots,T_\ell)$, and functions $f_1,\dots,f_s\in L^\infty(\mu)$,   the limit
$$
\limav \psi(\bn)\cdot T_{\vec{p_1}(\bn)}f_1\cdot\ldots\cdot T_{\vec {p_s}(\bn)}f_s
$$
exists in $L^2(\mu)$.
\end{proposition}
\begin{remark}
For $d=1$ this was proved in \cite[Section~2.4]{F15}.
\end{remark}

Our main convergence criterion for uniform averages is the next result which  is proved in Section~\ref{subsec:Proof-Th2}.
\begin{theorem}
\label{th:weighted-L2}
Let $d,\ell,s, t\in\N$. Then there exists a positive integer $k=k(d,\ell,s,t)$, such that  the following holds:
If  $w\in \ell^\infty(\N^d)$ is a  sequence and  for every $k$-step nilsequence $\psi\colon \N^d\to \C$   the limit
\begin{equation}
\label{eq:av-w-psi}
\limav w(\bn)\, \psi(\bn)  \ \text{ exists},
\end{equation}
  then
for every  system $(X,\mu,T_1,\dots,T_\ell)$,
 functions $f_1,\dots,f_s\in L^\infty(\mu)$, and polynomial mappings
 $\vec{p_i}\colon\N^d\to\Z^\ell$, $i=1,\dots,s$, of degree at most $t$, the limit \begin{equation}
\label{eq:weighted-L2}
\limav w(\bn)\cdot T_{\vec{p_1}(\bn)}f_1\cdot\ldots\cdot T_{\vec {p_s}(\bn)}f_s
\end{equation}
exists in  $L^2(\mu)$.
Furthermore, if the limit in~\eqref{eq:av-w-psi} is zero for every $k$-step nilsequence $\psi$ in $d$ variables, then the limit in~\eqref{eq:weighted-L2} is always zero.
Lastly, if the polynomial mappings $\vec{p_i}\colon \N^d\to\Z^\ell,$ $i=1,\ldots, s$, are linear, then we can take $k=s$.
\end{theorem}
\begin{remarks}
$\bullet$ For $d=\ell=t=1$ this result was proved in \cite{HK09} and for $d=\ell=1$ and $t\in \N$  arbitrary in \cite{C09}.

$\bullet$ In Theorem~\ref{th:regular} we give a characterization  using ``uniformity seminorms'' of sequences satisfying the hypothesis of Theorem~\ref{th:weighted-L2}.
\end{remarks}
The next result shows that
the hypothesis of Theorem~\ref{th:weighted-L2} is  necessary
in order to have weak convergence of the averages in \eqref{eq:weighted-L2} for all linear polynomial mappings.
\begin{proposition}
\label{prop:converse}
Let $d, s \in\N$ and $w\in \ell^\infty(\N^d)$ be a sequence. Suppose that for every system $(X,\mu,T_1,\dots,T_d)$, functions $f_0,\dots,f_s\in L^\infty(\mu),$ and linear forms  $\vec{L_i}\colon \N^d\to\Z^d$, $i=1,\dots,s$, the limit
$$
\limav w(\bn)\int f_0\cdot T_{\vec{L_1}(\bn)}f_1\cdot\ldots\cdot T_{\vec {L_s}(\bn)}f_s\, d\mu
$$
 exists. Then the limit
$$
\limav w(\bn)\, \psi(\bn)
$$
exists for every $s$-step nilsequence $\psi$ in $d$ variables.
\end{proposition}
\begin{remark}
For $d=1$ this was proved in \cite{F15}.
\end{remark}
Next for $d=1$ we give some examples of sequences of weights in  $\ell^\infty(\N)$ for which  Theorem~\ref{th:weighted-L2}  is applicable.
\subsection*{Examples}
Let $(Y,S)$ be a minimal uniquely ergodic system with invariant measure $\nu$ and  for every $s\in \N$, let $(Z_s,\nu_s,S)$ be the ``factor of order $s$''  defined in~\cite{HK05}. Suppose that for every $s\in \N$  the factor map $\pi_s\colon Y\to Z_s$ is continuous. Then for every  $\Psi \in C(Y)$ and every $y_0\in Y$, the sequence $w\colon N\to \C$ defined by $w(n):=\Psi(S^ny_0),$ $n\in\N$, satisfies the hypothesis of
 Theorem~\ref{th:weighted-L2}~\cite[Proposition~7.1]{HK09}. Examples of this type include:
\begin{enumerate}
\item
The Thue-Morse sequence which is the indicator function of those integers that have an odd sum of digits when expanded in base $2$ (see \cite[Proposition~2.21]{HK09}).
\item Bounded generalized polynomials (see~\cite[Corollary~2.23]{HK09}).\footnote{A \emph{generalized polynomial} is a real valued function that is obtained
from the identity function and real constants by using the
operations of addition, multiplication, and taking the integer part.} These  include sequences of the form $(\{p(n)\})_{n\in\N}$ or $(e(p(n)))_{n\in\N}$, where $p\colon \N\to \Z$
    is an arbitrary generalized polynomial,
 $\{x\}$ denotes the fractional part of $x$, and $e(t):=e^{2\pi i t}$.
\smallskip
\end{enumerate}



\subsection{Convergence results for Ces\`aro averages}\label{SS:2.3}
The  assumptions of Theorem~\ref{th:weighted-L2} are in
many cases too strong to be of use (this is the case for the examples $\text{(i)}$-$\text{(v)}$ below)
and we would like to have a criterion that uses convergence assumptions  of certain Ces\`aro averages instead of uniform averages.
We obtain such a result by   utilizing  tools different from those used in the proof of
 Theorem~\ref{th:weighted-L2} (Theorem~\ref{th:weighted-L2} relies on  Theorem~\ref{th:regular} which  necessitates the hypothesis \eqref{eq:av-w-psi}).
 The key new ingredients are decomposition results for multiple correlation sequences that are stated in Section~\ref{SS:mcs} below.
\begin{theorem}
\label{th:Cesaro}
Let $d,\ell,s,t\in\N$. Then there exists a positive integer $k=k(d,\ell,s,t)$, such that
the following holds:
\begin{enumerate}
\item
\label{it:Cesaro1}
If $w\in \ell^\infty(\N^d)$ and
for every $k$-step nilsequence $\psi\colon \N^d\to \C$
the limit
\begin{equation}
\label{eq:Cesaro-av1}
\lim_{N\to+\infty}\frac 1{N^d}\sum_{\bn \in [N]^d}w(\bn)\, \psi(\bn) \ \text{ exists,}
\end{equation}
 then  for every system $(X,\mu,T_1,\dots,T_\ell)$, functions $f_0,\dots, f_s\in L^\infty(\mu)$, and polynomial mappings
 $\vec{p_i}\colon\N^d\to\Z^\ell$, $i=1,\dots,s$, of degree at most $t$, the limit
\begin{equation}
\label{eq:Cesaro-av2}
\lim_{N\to+\infty}\frac 1{N^d}\sum_{\bn \in [N]^d}w(\bn)\cdot \int f_0\cdot  T_{\vec{p_1}(\bn)}f_1\cdot\ldots\cdot T_{\vec{p_s}(\bn)}f_s\,d\mu
\end{equation}
exists.
\item
\label{it:Cesaro2}
If $w\in \ell^\infty(\N^d)$  and
 for every $k$-step nilsequence $\psi\colon \N^{2d}\to \C$ 
the limit
\begin{equation}
\label{eq:Cesaro-av3}
\lim_{N,N'\to+\infty}\frac 1{(NN')^d}
\sum_{\bn\in[N]^d,\ \bn'\in[N']^d}
w(\bn)\, \overline{w(\bn')}\,
\psi(\bn,\bn')
\end{equation}
exists,  then  for every system $(X,\mu,T_1,\dots,T_\ell)$, functions  $f_1,\ldots, f_s\in L^\infty(\mu)$, and polynomial mappings
 $\vec{p_i}\colon\N^d\to\Z^\ell$, $i=1,\dots,s$, of degree at most $t$, the limit
\begin{equation}
\label{eq:Cesaro-av4}
\lim_{N\to+\infty}\frac 1{N^d}\sum_{\bn \in [N]^d}w(\bn)  \cdot T_{\vec{p_1}(\bn)}f_1\cdot\ldots\cdot T_{\vec{p_s}(\bn)}f_s
\end{equation}
 exists in  $L^2(\mu)$.

\end{enumerate}

 Furthermore, if the limit in~\eqref{eq:Cesaro-av1} (resp. \eqref{eq:Cesaro-av3}) is zero for every $k$-step  nilsequence $\psi \colon \N^d\to \C$  (resp. $\psi\colon \N^{2d}\to \C$), then the limit in~\eqref{eq:Cesaro-av2} (resp. \eqref{eq:Cesaro-av4}) is always zero.

  Lastly, in \eqref{it:Cesaro1} (resp. \eqref{it:Cesaro2}), if  the polynomial mappings $\vec{p_i}$ are linear, then  we can take $k=s$ (resp. $k=2s-1$), and if in addition
$\ell=s$ and $T_{\vec{p_i}(\bn)}=T_i^{L_i(\bn)}$ for  $i=1,\ldots, s$,  where $L_1,\ldots, L_s$ are linear forms spanning a subspace of dimension $r$, then we can take $k=s-r+1$
(resp. $k=2s-2r+1$).

\end{theorem}
\begin{remarks}
$\bullet$ For  $d=1$ the first part of this result was proved in \cite{F15}.

$\bullet$ For single variable polynomials,  in order to prove Part $(ii)$ of this result, we rely on decomposition results of correlations sequences involving polynomials in two variables.
\end{remarks}
An analogue of Proposition~\ref{prop:converse}, with Ces\`{a}ro averages in place of uniform averages, holds  with the same proof. This implies that   the condition~\eqref{eq:Cesaro-av1}
is also necessary in order for the limit \eqref{eq:Cesaro-av2} to exist for all linear polynomial mappings.

Sequences  $w\in \ell^\infty(\N^d)$ that satisfy the hypothesis~\eqref{eq:Cesaro-av1} and~\eqref{eq:Cesaro-av3}
 of Theorem~\ref{th:Cesaro} (but do not satisfy the condition~\eqref{eq:av-w-psi} of Theorem~\ref{th:weighted-L2} even for $\psi=1$) include the following:
\begin{enumerate}

\item
Any  sequence of the form
$(g(S_{\vec n} y))_{\vec n\in \N^d}$
 where $(Y,\nu, S_1,\ldots, S_d)$ is a system and $g\in L^\infty(\nu)$, for $y\in Y$ belonging to a set
of full measure that depends only on the system and the function $g$. For $d=1$ this was proved in \cite[Theorem~2.22]{HK09} for hypothesis~\eqref{eq:Cesaro-av1} but a similar argument also gives hypothesis~\eqref{eq:Cesaro-av3}
and works (using Theorem~\ref{th:structure}) for general $d\in \N$.

\item
 Any ``good''  multiplicative function  $\phi\colon\N^d\to\C$ (see Section~\ref{SS:arithmetic} and Theorem~\ref{th:double-av1}).
  For $d=1$ an alternate proof which depends on~\cite{FH15a}, and thus on  deep results from~\cite{GT12c}
  and~\cite{GTZ12},  was given in~\cite{FH15b}.

  \item The indicator function of all vectors of $\N^d$ whose coordinates have an even (or an odd) number of distinct prime
  factors; or more generally, the indicator function of any set $S$ defined as in Theorem~\ref{T:convarith} below
(this follows from Theorem~\ref{th:double-av1} and the argument in Section~\ref{SS:pf}).

\item  Any sequence of the form
 $(e(\sum_{i=1}^d n_i^{a_i}))_{n_1,\ldots,n_d\in\N}$, where $a_1,\ldots, a_d$ are  positive non-integers. In this case the limits in ~\eqref{eq:Cesaro-av1} and~\eqref{eq:Cesaro-av3}
  are always  $0$ (see Theorem~\ref{T:hardymain} and Proposition~\ref{prop:hardyequi}). Moreover, for $d=1$,  we give necessary and sufficient conditions for  a sequence of the form
$(e(f(n)))_{n\in\N}$, where $f$ is a  Hardy field function of at most polynomial growth, to be a good universal
weight for mean convergence of the averages~\eqref{eq:Cesaro-av4}   (see Corollary~\ref{C:hardymain}).

\item  The indicator function of any set of the form
$$
S:=\{n_1,\ldots, n_d\in \N\colon \norm{f_1(n_1)}\in [a_1,b_1], \ldots, \norm{f_d(n_d)}\in [a_d,b_d]\},
 $$
 where $0\leq a_i<b_i\leq 1/2$, $\norm{x}:=d(x,\Z)$,  and     $f_i$ are   Hardy field  functions
    of at most polynomial growth that stay away from polynomials (this follows from
    Proposition~\ref{prop:hardyequi} and an approximation argument that uses the estimate   \eqref{E:approximate}).
Furthermore, the set $S$ is good for multiple recurrence and mean convergence  and  the $L^2(\mu)$ limit
    $$
\lim_{N\to+\infty} \frac{1}{|S\cap [N]^d|}\sum_{\bn \in S\cap [N]^d}  T_{\vec{p_1}(\bn)}f_1\cdot\ldots\cdot T_{\vec {p_s}(\bn)}f_s
$$
    is equal to the limit obtained when $S$ is replaced by $\N^d$ (see Theorem~\ref{C:hardyrec}).
    \end{enumerate}

\subsection{Multiple correlations in ergodic theory}\label{SS:mcs}
Multiple correlation sequences are a well studied object in ergodic theory
 and form  an indispensable tool in the study of various multiple ergodic averages.
 For single variable sequences, structural results have been obtained in
\cite{BHK05, F15,  L10, L15, Mei90}; we extend some of these results to sequences in several variables. These extensions
turn out to be  key for the proof of
the convergence criterion given in  Theorem~\ref{th:Cesaro}.
Our argument follows closely the method used
in~\cite{F15} to obtain similar results for single variable sequences;
 but also  some refinements obtained (for example Theorem~\ref{th:correl-lin}) require new methodology.

\begin{notation}
For a bounded sequence $a\colon \N^d\to \C$ we let
\begin{equation}\label{E:seminorm2}
\norm a_2:=
\big(\Limsup\;\Av\,|a(\bn)|^2\big)^{\frac{1}{2}}.
\end{equation}
\end{notation}

\begin{definition}
A bounded sequence $a\colon \N^d\to \C$ is an \emph{approximate $s$-step nilsequence in $d$ variables} if for every $\ve>0$ it admits a decomposition as $a=a_\st+a_\er$, where
\begin{enumerate}
\item
$a_\st\colon \N^d\to \C$ is an $s$-step nilsequence in $d$ variables with $\norm{a_\st}_\infty\leq \norm a_\infty$;
\item
$\norm{a_\er}_2\leq \ve$.
\end{enumerate}
\end{definition}
  The subscripts ``st'' and ``err'' are used to indicate  ``structured'' and ``error'' respectively.
In  Section~\ref{subsec:Proof-Th5} we show:
\begin{theorem}
\label{th:correl-poly}
Let $d,\ell,s,t\in\N$. Then there exists a positive integer $k=k(d,\ell,s,t)$, such that
for every system $(X,\mu,T_1,\dots,T_\ell)$, functions $f_0,\dots, f_s\in L^\infty(\mu)$,
and polynomial mappings $\vec{p_i}\colon\N^d\to\Z^\ell$, $i=1,\dots,s,$ of degree at most $t$,
 the sequence $a\colon \N^d\to \C$ given by
\begin{equation}
\label{eq:correl-poly}
a(\bn):=\int f_0\cdot T_{\vec{p_1}(\bn)}f_1\cdot\ldots\cdot
T_{\vec{p_s}(\bn)}f_s\,d\mu, \quad \bn \in \N^d,
\end{equation}
is an approximate
 $k$-step nilsequence in $d$ variables.

Furthermore, if the polynomial mappings $\vec{p_1}, \ldots, \vec{p_s}$ are linear, then we can take $k=s$.
\end{theorem}
\begin{remark}
 For  $d=1$ this result was proved in \cite{F15}.
\end{remark}
The degree of nilpotency provided in  the last part of Theorem~\ref{th:correl-poly} is not always optimal. In  Section~\ref{subsec:Proof-Th6b}
we establish the following improvement for particular correlation sequences:
\begin{theorem}
\label{th:correl-lin}
For $d,\ell\in\N$ let $(X,\mu,T_1,\dots,T_\ell)$ be  a system, $f_0,\dots,f_\ell\in L^\infty(\mu)$ be functions, and $L_1,\dots,L_\ell\colon\N^d\to\Z$ be linear forms spanning a space of dimension $r$.  Then  the sequence $a\colon \N^d\to \C$   given by
$$
a(\bn):=\int f_0\cdot T_1^{L_1(\bn)}f_1\cdot\ldots\cdot T_\ell^{L_\ell(\bn)}f_\ell\,d\mu, \quad \bn \in \N^d,
$$
is an approximate $(\ell-r+1)$-step nilsequence in $d$ variables.
\end{theorem}
\begin{remark}
Examples of sequences for which this theorem applies and gives the  optimal degree of nilpotency are the three sequences in \eqref{E:three}.
\end{remark}

A crucial ingredient in the proof of the previous two decomposition results is
Theorem~\ref{th:regular+antiunif} which gives a characterization involving uniformity seminorms of
approximate nilsequences.
Furthermore, the proof of Theorem~\ref{th:correl-lin} uses a structural result for the generalized Kronecker factor of a not necessarily ergodic system   that is of independent interest (see  Theorem~\ref{th:basis-eigen}).

Lastly, we give an interesting corollary of Theorem~\ref{th:correl-poly}.
  For $d \in \N$ we consider various subsets of $\ell^\infty(\N^d)$. The first is the set
$$
\mathcal{N}_d:=\big\{(\psi(\bn))_{\bn\in  \N^d} \colon \psi \text{ is a } \text{nilsequence in } d \text{ variables}\big\}.
$$

With
$\mathcal{MC}_\text{d,pol}$ we denote the set that contains all  sequences of the form
$$
\Big(\int  f_0\cdot  T_{\vec{p_1}(\bn)}
f_1\cdot \ldots \cdot T_{\vec{p_s}(\bn)}
f_s\, d\mu\Big)_{\bn \in \N^d}
$$
for arbitrary systems
 $(X, \mu, T_1,\ldots, T_\ell)$, functions
$f_0,  \ldots, f_s \in L^\infty(\mu)$,  polynomial mappings
 $\vec{p_1}, \ldots, \vec{p_s} \colon \Z^d \to \Z^\ell$, and $\ell, s\in \N$.

We also define as $\mathcal{MC}_\text{d,lin}$ the set of multiple correlation sequences defined as above using linear polynomial mappings only.
\begin{theorem}
\label{T:4}
For every $d\in \N$ the sets $\mathcal{N}_d$,  $\mathcal{MC}_\text{d,lin}$,  $\mathcal{MC}_\text{d,pol}$
are subspaces of $\ell^\infty(\N^d)$ and
we have
$$
\overline{\mathcal{N}_d}=\overline{\mathcal{MC}_\text{d,lin}}=\overline{\mathcal{MC}_\text{d,pol}}
$$
where the closure is taken with respect to the seminorm $\norm{\cdot}_2$ defined in \eqref{E:seminorm2}.
\end{theorem}


\subsection{Multiple correlations for sequences in $\N^d$ and $\Z_N^d$}\label{SS:2.5} 
We use the  decomposition results of the previous subsection in order to deduce similar results for multiple correlations of bounded sequences in $\N^d$.
\begin{definition}
Let $\CA$ be a finite collection
of bounded complex valued  sequences  in $\ell$ variables and $\bI=(I_j)_{j\in \N}$ be a F\o lner sequence in $\N^\ell$. We say that the collection $\CA$ \emph{admits correlations along $\bI$} if for every $s\in\N$ and all $\bh_1,\dots,\bh_s\in\N^\ell$, the limit
$$
\Limav{\bk,\bI} \prod_{j=1}^sb_j(\bk+\bh_j)
$$
exists, where  for $j=1,\dots,s$ the sequence $b_j$ or the sequence $\overline{b_j}$ belongs to $\CA$.\end{definition}

Combining Theorem~\ref{th:correl-poly} with the correspondence principle stated in
Proposition~\ref{P:correspondence} below, we deduce the following statement:
\begin{theorem}
\label{th:cor-infty}
 Let $d,\ell,s, t\in\N$. Then there exists a positive integer $k=k(d,\ell,s, t)$ such that the following holds: If
 $a_1,\dots,a_s\colon \Z^\ell\to \C$ are  bounded sequences admitting correlations along  a F\o lner sequence  $\bI$ in $\N^\ell$ and $\vec{p_i}\colon \N^d\to\Z^\ell$, $i=1,\ldots,s$, are  polynomial mappings of degree at most $t$,  then  the sequence $b\colon\N^d\to\C$ defined by
 $$
b(\bn):=
\Limav{\bk,\bI}\prod_{i=1}^s
a_i\bigl(\bk+ \vec{p_i}(\bn)\bigr), \quad \bn\in \N^d,
$$
 is an approximate $k$-step nilsequence in $d$ variables.

Moreover, if the polynomial mappings are linear, then we can take $k=s-1$.
 \end{theorem}

 If $A$ is a finite set we let  $\E_{n\in A}:=\frac{1}{|A|}\sum_{n
 \in A}$. Decomposition results of similar nature also hold in the finite world, for example, the following is true:
\begin{theorem}\label{th:cor-finite}
Let $d,\ell,  s,t \in\N$. Then there exists a positive integer $k=k(d, \ell, s, t)$ such that the following holds:
For every $\varepsilon>0$ there exists a $k$-step  nilmanifold $X=X(d,\ell, s,t,\varepsilon)$ such that
for every $N\in \N$,
finite sequences  $a_1,\dots,a_s\colon \Z_N^\ell\to \C$ of modulus   at most  $1$,
and polynomial mappings
  $\vec{p_i}\colon \N^d\to\Z^\ell$, $i=1,\ldots, s,$ having integer coefficients and   degree at most $t$,
the sequence $b\colon\N^d\to\C$ defined by
 \begin{equation}\label{E:corfin}
b(\bn):=
\E_{\bk\in \Z_N^\ell}
\prod_{i=1}^s
a_i\bigl(\bk+ \vec{p_i}(\bn)\bigr), \quad \bn\in \Z_N^d,
\end{equation}
admits a decomposition of the form $b=b_\st+b_\er$ where
\begin{enumerate}
\item
$b_\st \colon \N^d\to \C$ is a convex combination of
$k$-step nilsequences  defined by functions on $X$ with  Lipschitz norm at most $1$;
\item
$
\E_{\bn\in \Z_N^d}|b_\er(\bn)|\leq\ve$.
\end{enumerate}

Furthermore, if the polynomial mappings  are linear, then we can take $k=s-1$.
\end{theorem}
\begin{remark}
It is important that the nilmanifold $X$ and the  Lipschitz norm of the function defining the nilsequence are independent of $N\in \N$.
\end{remark}
\subsection{Applications to arithmetic}\label{SS:arithmetic} Next we give some applications  of number theoretic
and combinatorial flavor.
 \begin{definition}
A function $\phi\colon \N\to \C$ is called \emph{multiplicative} if
 $$
 \phi(mn)=\phi(m)\, \phi(n) \ \text{ whenever } \  (m,n)=1.
 $$
 It is called \emph{completely multiplicative} if this relation holds for all $m,n\in\N$.

We say that a  multiplicative function  $\phi\colon \N\to \C$  that is bounded by $1$ is \emph{good} if  the limit
\begin{equation}
\label{eq:av-prog}
\lim_{N\to+\infty} \frac 1N\sum_{n=1}^N \phi(an+b)
\end{equation}
exists for all $a\in \N,b\in \Z_+$. It is called \emph{aperiodic} if all these limits are equal to $0$.
\end{definition}

For $d\in \N$  a function $\phi\colon \N^d\to \C$ is called \emph{multiplicative} if it is of the form
 $$
 \phi(n_1,\ldots,n_d)=\phi_1(n_1)\cdots \phi_d(n_d), \quad n_1,\ldots, n_d\in \N,
 $$ for some multiplicative functions $\phi_i\colon \N\to
 \C, i=1,\ldots, d, $ which we call  the \emph{components of} $\phi$.
 We call a  multiplicative function $\phi\colon \N^d\to \C$  \emph{good} if all its component functions are good  and
\emph{aperiodic}
if at least one of its component functions is aperiodic.

By a  classical result of Wirsing~\cite{Wi67}, every real valued multiplicative function
that is bounded by $1$ is good. A result of
Hal\'asz~\cite{Hal68}  allows to characterize   good and aperiodic multiplicative functions. Let $\PP$ be the set of primes.
A Dirichlet character is a periodic completely multiplicative function which takes the value $1$ at $1$.
\begin{notation}[\cite{GS15}]
If $\phi_1,\phi_2\colon \N\to \C$ are  multiplicative functions, bounded by $1$, we define $\D(\phi_1,\phi_2)\in[0,+\infty]$ by
$$
\D(\phi_1,\phi_2)^2:=\sum_{p\in\PP} \frac{1}{p}\,\bigl(1-\Re\bigl(\phi_1(p) \overline{\phi_2(p)}\bigr)\bigr).
$$
\end{notation}
\begin{remark}
Note that if $|\phi_1|=|\phi_2|=1$, then
$
\D(\phi_1,\phi_2)^2=\sum_{p\in \PP} \frac{1}{2p}\, |\phi_1(p)- \phi_2(p)|^2.
$
\end{remark}
The next result can be deduced from~\cite[Theorem~6.3]{E79}.
\begin{theorem}\label{T:Halasz}
Let $\phi\colon \N\to \C$ be a  multiplicative function that is bounded by $1$. Then $\phi$ is good   if and only if for every Dirichlet character $\chi$  we either  have
\begin{enumerate}
\item
\label{it:good1}
$\D(\phi\, \chi,n^{it})=+\infty$ for every $t\in\R$, or
\item
\label{it:good3}
for some $t\in \R$ we  have
$\D(\phi\chi,n^{it})<\infty$  and
   $\chi(2)^k\phi(2^k)=-2^{ikt}$  for all $k\in \N$, or
   \item
$\sum_{p\in\PP}\frac{1}{p}(1-\phi(p)\chi(p))$ converges.
\end{enumerate}

Moreover, $\phi$ is aperiodic if and only if either the condition~$(i)$
or~$(ii)$
 is satisfied for every Dirichlet character $\chi$.
\end{theorem}
For a more complete  discussion of these notions, see~\cite[Section~ 2.5]{FH15b}.
The next result is proved in Section~\ref{subsec:mult}.
\begin{theorem}
\label{th:cesaro-mult}
Let  $d\in \N$ and $\phi\colon \N^d\to \C$ be a good    multiplicative function. Then for every $\ell,s\in \N$, system $(X,\mu,T_1,\dots,T_\ell)$,
functions $f_1,\ldots, f_s\in L^\infty(\mu)$,  and polynomial mappings $\vec{p_i}\colon\N^d\to\Z^\ell$, $i=1,\dots,s$,   the limit
\begin{equation}\label{E:phiav}
\lim_{N\to+\infty}\frac 1{N^d}\sum_{\bn \in [N]^d} \phi(\bn) \cdot
   T_{\vec{p_1}(\bn)}f_1\cdot\ldots\cdot T_{\vec{p_s}(\bn)}f_s \quad \text{exists in } L^2(\mu).
\end{equation}

 Furthermore, if the multiplicative function $\phi$  is aperiodic, then the limit is equal to $0$.
\end{theorem}
\begin{remarks}
$\bullet$  This result and its consequences below were proved in \cite{FH15b} for $d=1$ using a deep structural result for multiplicative functions from \cite{FH15a}. The current argument relies on the convergent criterion of Theorem~\ref{th:Cesaro} and  uses much softer number theoretic input (we only use  Theorem~\ref{th:double-av1}).

$\bullet$ Conversely, if for $d=\ell=s=1$ the averages in \eqref{E:phiav} converge weakly, then examples of periodic systems show that $\phi$ has to be good, and if the averages in \eqref{E:phiav} converge weakly to $0$, then $\phi$ has to be aperiodic.

$\bullet$ Similar statements, with similar proofs,  hold if in \eqref{E:phiav} we use averages of the form
$\frac 1{N_1\cdots N_d}
\sum_{\bn \in [N_1]\times \cdots \times [N_d]}$ and take the limit as $N_1,\ldots, N_d\to+\infty$. A similar comment applies for the next two results.
\end{remarks}
 \begin{definition} For $a\in \Z_+$ and $b\in \N$   we let
$S_{a,b}$ consist of  those $n\in \N$ whose number of  distinct prime factors is congruent to $a$ $\! \! \! \mod{b}$.
\end{definition}
We can also define $S_{a,b}$ by counting prime factors with multiplicity; then   all results stated below continue to hold with similar proofs.
The next  result is proved in Section~\ref{SS:pf}.
\begin{theorem}\label{T:convarith}
Let $d\in \N$, $a_i,c_i\in \Z_+$, $b_i\in \N$,  $i=1,\ldots, d$,      and let
\begin{equation}\label{E:S}
S:=(S_{a_1,b_1}+c_1)\times \cdots \times (S_{a_d,b_d}+c_d).
\end{equation}
Then for all  $\ell,s\in\N$, polynomial mappings  $\vec{p_1},\ldots, \vec{p_s}\colon \N^d\to\Z^\ell$,
 system $(X,\mu,T_1,\dots,T_\ell)$, and functions $f_1,\dots,f_s\in L^\infty(\mu)$, the limit
$$
\frac{1}{|S\cap [N]^d|}\sum_{\bn \in S\cap [N]^d}  T_{\vec{p_1}(\bn)}f_1\cdot\ldots\cdot T_{\vec {p_s}(\bn)}f_s
$$
exists in $L^2(\mu)$ and is equal to the limit  obtained when one replaces $S$ with $\N^d$.
\end{theorem}
\begin{remark}
 It follows from our argument  that  $\lim_{N\to+\infty}|S\cap [N]^d|/N^d=(\prod_{i=1}^db_i)^{-1}$.
\end{remark}
In Section~\ref{SS:pf} we deduce from this result the following multiple recurrence statement:
\begin{theorem}\label{T:recarith}
 We use the notation of Theorem~\ref{T:convarith} and assume in addition that $\vec{p_i}({\bf 0})=\vec{0}$ for $i=1,\ldots, s$. Then for $S$ as in \eqref{E:S}, for  every $A\in \CX$ with $\mu(A)>0$ we have
$$
\lim_{N\to +\infty} \frac{1}{|S\cap [N]^d|}\sum_{\bn \in S\cap [N]^d} \mu(A\cap T_{-\vec{p_1}(\bn)}A\cap \cdots\cap T_{-\vec {p_s}(\bn)}A)>0.
$$
\end{theorem}

Lastly,
we give some combinatorial implications of the previous
multiple recurrence result. We define the {\em
upper Banach density} $d^*(E)$ of a set $E\subset\Z^\ell$
by $d^*(E):= \limsup_{|I|\to +\infty}\frac{|E\cap I|}{|I|}$,
where the $\limsup$ is taken over all parallelepipeds $I\subset\Z^\ell$
whose side lengths tend to infinity and the {\em lower natural density} as
$ \liminf_{N\to\infty}\frac{|E\cap [-N,N]^\ell|}{(2N+1)^\ell}$.
Using a
 modification of the  correspondence principle of H.~Furstenberg  (proved as  in~\cite{BL96})
  we deduce from Theorem~\ref{T:recarith} the following result:
\begin{theorem}
\label{T:combinatorics}
 We use the notation of Theorem~\ref{T:convarith} and assume in addition that $\vec{p_i}({\bf 0})=\vec{0}$ for $i=1,\ldots, s$.
Then for  $S$ as in \eqref{E:S}, for every set
$E\subset\Z^\ell$ with
$d^*(E) > 0$,  the set
$$
\big\{\bn\in S\colon
 d^*\bigl(E\cap (E-\vec{p_1}(\bn))\cap\ldots\cap (E-\vec{p_s}(\bn))
\bigr)>0
\big\}
$$
has positive lower natural density.
\end{theorem}
 Applying this for $d=2$,  $a_i=2$, $b_i=0$ or $1$,  $c_i=0$,  and $\vec{p_i}(n_1,n_2)=i(n_1^2+n_2^2)$ for $i=1,\dots,s,$ we obtain the refinement of Szemer\'edi's theorem mentioned towards the end of the introduction.
\subsection{Open problems}
When  all the maps $T_1,\ldots, T_\ell$ are powers of the same transformation and $d=1$,
 a strengthening of Theorem~\ref{th:correl-poly} holds
which shows that the error term can be taken to converge to $0$ in uniform density
  (see \cite{BHK05, L10, L15, L15b}).   It is not clear whether a similar
 result holds for arbitrary commuting transformations, even when $d=1$, $\ell=2$,  and the polynomial mappings are linear.

\medskip
\noindent {\bf Problem 1.} \
{\em Let $(X,\mu,T_1,T_2)$ be a system and $f_0,f_1,f_2\in L^\infty(\mu)$. Is it true that the sequence
$a\colon \N\to \C$ defined by
$$
a(n):=\int f_0\cdot T_1^nf_1\cdot T_2^nf_2\, d\mu, \quad  n\in \N,
$$
can be decomposed as $a=a_{\st}+a_{\er}$ where $a_\st$ is
 a  uniform limit of $2$-step nilsequences and $\norm{a_\er}_2=0$?}
\medskip

When $T_2=T_1^2$ this is shown to be the case in  \cite{BHK05, L15, L15b}.

Theorem~\ref{th:Cesaro} shows that condition \eqref{eq:Cesaro-av1} is sufficient for
weak convergence of the averages in \eqref{eq:Cesaro-av4}, but we needed the
stronger hypothesis \eqref{eq:Cesaro-av3} in order to guarantee mean convergence. This is probably an artifact of
our proof.

\medskip
\noindent {\bf Problem 2.} \
{\em Show that condition \eqref{eq:Cesaro-av1} is sufficient for
mean convergence of the averages in \eqref{eq:Cesaro-av4}.}

\medskip

When $d=1$ and all the maps $T_1,\ldots, T_\ell$ are powers of the same transformation this is shown to be the case in
\cite[Theorem~1.3]{C09}.

In Theorem~\ref{th:correl-poly}, even in seemingly simple cases,  it is not clear what the  optimal
 dependence of $k$ on $d,\ell,s,t$ is, even when the polynomial mappings are linear.
It is expected (but we are unable to verify this) that this optimal dependence can already be inferred
from the case where all the $T_1,\ldots, T_\ell$ are powers of the same transformation (a case which is much more tractable using the theory of characteristic factors).
We record here a relevant open problem.

\medskip
\noindent {\bf Problem 3.} \
{\em Let $(X,\mu,T_1,T_2)$ be a system and $f,g,h\in L^\infty(\mu)$. Show that the sequence
$a\colon \N^2\to \C$ defined by
$$
a(m,n):=\int f\cdot T_1^m T_2^ng\cdot T_1^n T_2^mh\, d\mu, \quad m,n\in \N,
$$
is an approximate $1$-step nilsequence in $2$ variables.}
\medskip

 When $T_1,T_2$ are powers of the same transformation  this can be verified by an argument similar to the one used in the proof of Theorem~\ref{th:correl-lin}. Note also that  Theorem~\ref{th:correl-poly} gives that the sequence $a$ is an approximate  $2$-step nilsequence in $2$ variables.
\subsection{Notation and conventions} For readers convenience, we gather here   some notation used throughout the article.

  \noindent We denote by $\N$ the set of positive integers, by  $\Z_+$ the set of non-negative integers, and by $\R_+$ the set of non-negative real numbers.

 \smallskip

\noindent For $N\in \N$ we denote by $[N]$ the set $\{1,\ldots,
N\}$.

\smallskip

\noindent With $\ell^\infty(\N^d)$ we denote the space of all bounded sequences
$a\colon \N^d\to \C$.

\smallskip

\noindent
If $A$ is a finite set we let  $\E_{n\in A}:=\frac{1}{|A|}\sum_{n
 \in A}$.


\smallskip

\noindent We write $\conj\colon \C\to\C$ for the complex conjugation.

\smallskip
\noindent If $x$ is a real, $\e(x)$ denotes the number $e^{2\pi i
x}$ and  $\norm x$ denotes the distance between $x$ and the nearest
integer.


\smallskip

\noindent Given $d\in \N$ we write $\bn=(n_1,\dots,n_d)$ for a point
of $\Z^d$.

\smallskip
\noindent  We typically use the letter $\psi$ to denote nilsequences.

\section{Uniformity seminorms and  decomposition results}
\label{sec:uniformity-seminorms}
In this section we extend to sequences in $\ell^\infty(\N^d)$ some results established in~\cite{HK09} and~\cite{F15} for sequences in $\ell^\infty(\N)$. The statements and the proofs are completely similar and we only  give the necessary definitions and sketch the main steps of the proofs.

\subsection*{Some definitions and notation}

We write $\conj\colon \C\to\C$ for the complex conjugation; then $\conj^k z=z$  if $k$ is even and $\conj^k z=\overline{z}$ if $k$ is odd.
We let $\Cube k:=\{0,1\}^k$ and $\Cub k:=\Cube k\setminus\{\bzero\}$.
Elements of $\Cube k$ are written as  $\uepsilon=(\epsilon_1,\dots,\epsilon_k)$. We let $|\uepsilon|:=\epsilon_1+\cdots+\epsilon_k$.
Elements of $(\N^d)^k$ are written as $\ubh=(\bh_1,\dots,\bh_k)$ where $\bh_i\in\N^d$ for  $i=1,\ldots, k$.  For $\ubh\in(\N^d)^k$ and $\uepsilon\in\Cube k$ we let $\uepsilon\cdot\ubh:=\epsilon_1\bh_1+\dots+\epsilon_k\bh_k\in\N^d$.

\subsection{The definition of the seminorms}
\label{subsec:def-seminorms}
We follow Section~2 of~\cite{HK09}.

 We say that a finite or countable family $\CF$ of bounded sequences 
 in $\ell^\infty(\N^d)$ \emph{admits correlations} along a  F\o lner sequence $\bI=(I_j)_{j\in \N}$  in $\N^d$ if   the limit
$$
\Limav{\bn,\bI}\bigl(\prod_{i=1}^m
b_i(\bn+\bh_i)\bigr)
$$
 exists for every $m\in \N$, all $\bh_1,\dots,\bh_m\in\N^d,$ and all sequences $b_1,\dots,b_m\in\ell^\infty(\N^d)$ such that either  $b_i$ or $\overline{b_i}$ belongs to $\CF$ for
 $i=1,\dots, m$. We remark that from every F\o lner sequence $\bI$ we can extract a subsequence $\bI'$ so that a given family of sequences admits correlations along $\bI'$.

Suppose that the sequence $a\in \ell^\infty(\N^d)$ admits correlations along $\bI$. Then
for $k\in \N$ and $\ubh=(\bh_1,\dots,\bh_k)\in(\N^d)^k$ we write
$$
\cor_\bI(a; \ubh):=
\Limav{\bn,\bI}\bigl(
\prod_{\uepsilon\in\Cube k}\conj^{|\uepsilon|}a(\bn+\uepsilon\cdot\ubh)\bigr)
$$
and
$$
\norm{a}_{\bI,k}:=\Big( \lim_{H\to+\infty}\frac 1{H^{dk}}\sum_{\bh_1,\dots,\bh_k\in[H]^d}\cor_\bI(a; \ubh)\Big)^{1/2^k}.
$$
In \cite[Proposition 2.4]{HK09} it is shown that for $d=1$ the previous limit exists and is non-negative;  the proof is similar for general $d\in \N$.
Furthermore, the map $a\mapsto\norm a_{\bI,k}$ is subadditive (\cite[Proposition 2.5]{HK09} for $d=1$), that is,  if  the sequences $a, b,$ and $a+b$ admit correlations along $\bI$, then $\norm{a+b}_{\bI,k}\leq \norm{a}_{\bI,k}+\norm{b}_{\bI,k}$.

For $a\in \ell^\infty(\N^d)$ we  define
$$
\norm a_{U^k(\N^d)}:=\sup_{\bI}\norm a_{\bI,k}
$$
where the supremum is taken over  all F\o lner sequences $\bI$  in $\N^d$ for which the sequence $a$ admits correlations. Then the map
$a\mapsto \norm a_{U^k(\N^d)}$ is a seminorm on $\ell^\infty(\N^d)$, we call it the  \emph{uniformity seminorm of order $k$} of $a$.

\subsection*{Interpretation} Next, we interpret the previous definitions and results in dynamical terms. We use a variant of Furstenberg’s correspondence principle that enables to transfer results from ergodic theory to results about bounded sequences of complex numbers. We  follow
 the method used in~\cite[Section~6.1]{HK09} when $d=1$;  similar arguments work for general $d\in \N$ and we summarize them here.   For notational convenience we restrict to the case where the family $\CF$ contains only a single sequence $a\colon \N^d\to \C$, the general case being completely similar.  Let $D$ be the closed disk in $\C$ of radius $\norm{a}_\infty$
 and let $D^{\Z^d}$ be endowed with the product topology and with the natural shifts $T_1,\dots,T_d$ given by
$$
(T_ix)(n_1,\dots,n_d)=x(n_1,\dots,n_{i-1}, n_i+1,n_{i+1},\dots,n_d)
$$
for $i=1,\dots,d$, where we use the notation $\bn=(n_1,\dots,n_d)\in\Z^d$ and $x=(x(\bn))_{\bn\in\Z^d}\in D^{\Z^d}$.
We define  the continuous function $f\colon D^{\Z^d}\to \C$ by $f(x):=x(\bzero)$. Furthermore, we define
  the point
 $\omega$ in $D^{\Z^d}$ by
 $\omega(\bn):=a(\bn)$ for  $\bn\in\N^d$ and $\omega(\bn)=0$ otherwise.
 Then  we have
 $$
f(T_\bn \omega)=a(\bn), \quad \bn\in\N^d.
$$
Let $X$ be the closed orbit of $\omega$ under $T_1,\dots,T_d$. Then $(X,T_1,\dots,T_d)$ is a topological dynamical system and $\omega$ is a transitive point of this system, meaning
it has a dense orbit in $X$.

Let $\bI=(I_j)_{j\in \N}$ be a F\o lner sequence  in $\N^d$
 for which the sequence $a$ admits correlations. Let $\mu$ be a $w^*$-limit point for   the sequence of measures
 $$
\mu_j:=\frac 1{|I_j|}\sum_{\bn\in I_j}\delta_{T_\bn \omega}, \quad j\in \N.
$$
 Then $\mu$ is a probability measure on $X$, invariant under $T_1,\dots,T_d$, and by construction, for every $m\in \N$, all $\eta_1,\dots,\eta_m\in\{0,1\},$ and all $\bh_1,\dots,\bh_m\in\Z^d$, we have
 $$
\Limav{\bn,\bI}\big(\prod_{i=1}^m
\conj^{\eta_i}a(\bn+\bh_i)\big)=
\Limav{\bn,\bI}\big(\prod_{i=1}^m\conj^{\eta_i}
T_{\bh_i}f(T_\bn\omega)\big)
=\int \prod_{i=1}^m \conj^{\eta_i}T_{\bh_i}f\,d\mu.
$$
In particular, for $\ubh=(\bh_1,\dots,\bh_k)\in(\N^d)^k$, we have
$$
\cor_\bI(a; \ubh)=\int \prod_{\uepsilon\in\Cube k}\conj^{|\uepsilon|}T_{\uepsilon\cdot \ubh} f\,d\mu
$$
and thus
$$\norm a_{\bI,k}=\nnorm f_{\mu,k}$$
 where $\nnorm\cdot_{\mu,k}$ is the seminorm on $L^\infty(\mu)$ defined in~\cite{HK05} in the ergodic case and in~\cite{CFH11} in the general case.\footnote{The seminorms were defined for a single transformation but the properties we use extend immediately to the case of several commuting transformations.} We recall the definition and some properties of these seminorms in Appendix~\ref{sec:seminorms}. Note also that if $\psi\colon \N^d\to \C$ is a nilsequence of the form $(\Phi(T_{\bn}x))_{\bn \in \N^d}$, then $\psi$ admits correlations along every F\o lner sequence $\bI$  and
 \begin{equation}\label{E:psi}
 \norm{\psi}_{\bI,k}=\norm {\psi}_{U^k(\N^d)}=\norm{\Phi}_{\mu, k}
 \end{equation}
 where the last seminorm is defined with respect to the action of $T_\bn, \bn \in \N^d,$ on $X$.
From the properties~ \eqref{eq:seminorm1} and \eqref{eq:seminorm3}  we deduce
\begin{gather}
\label{eq:unif-seminorm1}
\Limsup\bigl|\Av_\bI\, a(\bn)\bigr\vert
\leq \norm a_{\bI,1};\\
\label{eq:unif-seminorms2}
\norm a_{\bI,k+1}^{2^k+1}=\lim_{H\to+\infty}
\frac 1{H^d}\sum_{\bh\in [H]^d}
\norm{\sigma_\bh a\cdot\overline a}_{\bI,k}^{2^k}, \quad k\in \N,
\end{gather}
where $\sigma_\bh a(\bn):=a(\bn+\bh)$ for $\bh,\bn\in \N^d$.

Then we have
(see \cite[Proposition~4.5]{HK09} for $d=1$)
\begin{align}
\notag
\norm a_{U^k(\N^d)}
&=
\sup_{\mu\text{ invariant probability measure on }X}
\nnorm f_{\mu,k}\\
\label{eq:norm-a-norm-f}
& =
\sup_{\mu\text{ invariant ergodic probability measure on }X}
\nnorm f_{\mu,k}
\end{align}
where the last equality follows by using the ergodic decomposition of the measure $\mu$.

\subsection{Tools}
 For an ergodic system $(X,\mu,T)$ the structure theorem of~\cite{HK05} links the seminorms $\nnorm\cdot_{\mu,k}$ with the factors of the system $(X,\mu,T)$ that are $(k-1)$-step nilsystems. This result was generalized to $\Z^d$-actions by Griesmer~\cite[Lemma 4.4.3 and Theorem 4.10.1]{Gr09} and an alternate proof that is based on finitestic inverse theorems was recently given by Tao   \cite[Remark~4]{T15}. We record an immediate corollary of this result that is more convenient for our purposes.
\begin{theorem}[\cite{HK05} for $d=1$, \cite{Gr09,T15} for general $d$]
\label{th:structure}
Let $d,k\in \N$,  $(X,\mu,T_1,\dots,T_d)$ be an ergodic system, $f\in L^\infty(\mu)$, and $\ve>0$.
Then there exists a $(k-1)$-step nilsystem
$(Y,\nu,T_1,\dots,T_d)$, a factor map $\pi\colon X\to Y$, and a continuous function $\Phi$ on $Y$ such that
$$
\norm{\Phi}_\infty\leq\norm f_\infty\ \text{ and }\ \nnorm{f-\Phi\circ\pi}_{\mu,k}\leq\ve.
$$
\end{theorem}
The next result can be considered as a strengthening of the correspondence principle of   Furstenberg.  We recall that a topological dynamical system $(Y,T_1,\dots,T_d)$ is \emph{distal} if for all $y\neq y'\in Y$ we have $\inf_{\bn\in \Z^d}d_Y(T_\bn,y,T_\bn y')>0$, where $d_Y$ is the distance on $Y$ defining its topology. It is known that every nilsystem is distal~\cite{AGH63}.

\begin{proposition}[\mbox{\cite[Proposition 6.1]{HK09}} for $d=1$]
\label{prop:corresp-distal}
Let $d\in \N$,  $(X,T_1,\dots,T_d)$ be a topological dynamical system,  $\omega\in X$ be a   transitive point, and $\mu$ be an invariant ergodic measure on $X$. Moreover, let $(Y,T_1,\dots,T_d)$ be a distal topological dynamical system,  $\nu$ be an invariant measure on $Y$, and $\pi\colon X\to Y$  be a measure theoretic factor map. Then there exists a point $y_0\in Y$  and a F\o lner sequence
$\bI=(I_j)_{j\in \N}$  on $\N^d$  such that
$$
\lim_{j\to+\infty} \frac{1}{|I_j|}\sum_{\bn\in I_j} f(T_\bn \omega) \, g(T_\bn y_0)=
\int_Xf \cdot g\circ\pi\,d\mu
$$
for every $f\in C(X)$ and every $g\in C(Y)$.
\end{proposition}
The important point in this statement is that we do not assume that  the map $\pi$ is continuous. The proof is exactly the same as in the case of a single transformation.  The  classical property of distal systems that we use is:

\begin{theorem}[\mbox{\cite[Chapter 5]{Aus88}}]
\label{th:Ellis}
Let $d\in \N$ and $(Y,T_1,\dots,T_d)$ be a distal system. Then
 for every $y_1\in Y$ and every sequence $(\bm_i)_{i\in \N}$ with values in  $\N^d$,  there exists $y_0\in Y$ and a subsequence
 $(\bm'_i)_{i\in \N}$ of   $(\bm_i)_{i\in \N}$  such that $T_{\bm'_i}y_0$ converges to $y_1$.
\end{theorem}
From Theorem~\ref{th:structure}
 and the discussion of Section~\ref{subsec:def-seminorms} we deduce:
\begin{proposition}\cite[Proposition~6.2 for $d=1$]{HK09}
\label{prop:approx}
Let $d,k\in \N$, $a\in\ell^\infty(\N^d)$ be a sequence, and $\ve >0$. Then there exists a F\o lner sequence $\bI=(I_j)_{j\in \N}$ and a $(k-1)$-step nilsequence $\psi_1$ in $d$ variables such that the sequences $a$ and $a-\psi_1$ admit correlations along $\bI$ and we have
$$
\norm a_{\bI,k}\geq\norm a_{U^k(\N^d)}-\ve, \quad
\norm{\psi_1}_\infty\leq \norm a_\infty, \quad  \text{and} \quad
\norm{a-\psi_1}_{\bI,k}\leq \ve.
$$
\end{proposition}
\begin{proof}
As explained above, there exists a system
$(X,T_1,\dots,T_d)$, a transitive point $\omega\in X,$ and a continuous function $f$ on $X$ such that $a(\bn)=f(T_\bn \omega)$ for every $\bn\in\N^d$. Note that then $\norm{a}_{\infty}=\norm{f}_\infty$. Moreover, by \eqref{eq:norm-a-norm-f} there exists an invariant ergodic probability measure $\mu$ on $X$ with  $\nnorm f_{\mu,k}\geq \norm a_{U^k(\N^d)} -\ve$.

Let the nilsystem $(Y,\nu,T_1,\dots,T_d)$, the factor map $\pi$,  and the function $\Phi$ be defined as in Theorem~\ref{th:structure}. Let $\sigma$ be the measure on $X\times Y$ which is the image of $\mu$ under the map $\text{id}\times\pi$. This measure is ergodic under the product action and thus admits a generic point $(x_1,y_1)$. Since $\omega$ is a transitive point of $X$, there exists a sequence $(\bm_j)_{j\in \N}$  in $\N^d$ such that $T_{\bm_j}\omega\to x_1$.  By Theorem~\ref{th:Ellis}, substituting a subsequence for the sequence $(\bm_j)_{j\in \N}$, we can assume that there exists  a point $y_0\in Y$ such that $T_{\bm_j}y_0\to y_1$ and thus $(T\times T)_{\bm_j}(\omega,x_0)\to (x_1,y_1)$. Since the point $(x_1,y_1)$ is generic, substituting again a subsequence for $(\bm_j)_{j\in \N}$ and defining the F\o lner sequence $\bI$ by $I_j=\bm_j+[j]^d$, $j\in \N$,  we obtain that the sequence of probability measures
$$
\frac 1{|I_j|}\sum_{\bn\in I_j}\delta_{T_\bn \omega}\times\delta_{T_\bn y_0}, \quad j\in \N,
$$
on $X\times Y$ converges weak* to  a probability measure $\sigma$.
Let the nilsequence $\psi_1$ be defined by $\psi_1(\bn):=\Phi(T_\bn y_0)$, $\bn\in \N^d$. We have $\norm{\psi_1}_\infty\leq \norm{\Phi}_\infty \leq \norm{f}_\infty=\norm a_\infty$.
By applying the preceding discussion to the product system on $X\times Y$, for  the point $(\omega,y_0)$,  the F\o lner sequence $\bI$,  and  the function given by $F(x,y)=f(x)$, we obtain $\norm a_{\bI,k}=\nnorm F_{\sigma,k}$.
Letting $G(x,y)=f(x)-\Phi(y)$ we have $\norm{a-\psi_1}_{\bI,k}=\nnorm G_{\sigma,k}$. By the definition of $\sigma$ and of the seminorms $\nnorm\cdot_k$ we have $\nnorm F_{\sigma,k}=\nnorm f_{\mu,k}$ and
  by the definition of $\Phi$ we have $\nnorm G_{\sigma,k}=\nnorm{f-\Phi\circ\pi}_{\mu,k}\leq\ve$. This completes the proof.
\end{proof}

\subsection{Anti-uniformity}
We introduce certain classes of sequences that are  asymptotically approximately orthogonal to
$k$-uniform sequences; in Theorem~\ref{th:structure}  we give a characterization of such sequences in terms of $(k-1)$-step nilsequences.

\begin{definition}
Let $a\in\ell^\infty(\N^d)$. We say that
\begin{itemize}
\item
 The sequence $a$ is \emph{strongly $k$-anti-uniform} if there exists a constant  $C\geq 0$ such that for every $b\in \ell^\infty(\N^d)$  we have
\begin{equation}
\label{eq:antinuiform}
\Limsup\bigl\vert\Av\,a(\bn)\, b(\bn)\bigr\vert\leq C\norm b_{U^k(\N^d)}.
\end{equation}
In this case, we write $\norm a_{U^k(\N^d)}^*$ for the smallest constant $C$ such that~\eqref{eq:antinuiform} holds.
\item
The sequence  $a$ is \emph{$k$-anti-uniform} if for every $\ve >0$ there exists $C=C(\ve)\geq 0$ such that
for every $b\in \ell^\infty(\N^d)$ we have
 $$
\Limsup\bigl\vert\Av\, a(\bn)\, b(\bn)\bigr\vert\leq C\norm b_{U^k(\N^d)}+\ve\norm b_\infty.
$$
\end{itemize}
\end{definition}

\begin{proposition}[\mbox{\cite[Section~5]{HK09}} for $d=1$]
\label{prop:dual-antiunif}
 Let $d, k\in \N$,  $(X,T_1,\dots,T_d)$ be an ergodic $(k-1)$-step nilsystem, and $f_{\uepsilon}\in C(X)$ for  $ \uepsilon\in\Cub k$.
Then  the limit
\begin{equation}
\label{eq:phi}
\Phi(x):=\lim_{H\to+\infty}\frac 1{H^{dk}}\sum_{ \bh_1,\dots,\bh_k \in[0,H)^d}\;
\prod_{\uepsilon\in\Cub k}f_{\uepsilon}(T_{\uepsilon\cdot\ubh}x)
\end{equation}
exists for every $x\in X$ and  the convergence is  uniform in $x\in X$ (hence, $\Phi\in C(X)$). If $x_0\in X$  and $a\in\ell^\infty(\N^d)$ is the sequence defined by
$$
a(\bn):=\Phi(T_\bn x_0), \quad  \bn\in\N^d,
$$
then $a$ is strongly $k$-anti-uniform and
$$
\norm a_{U^k(\N^d)}^*\leq
\prod_{\epsilon\in\Cub k} \nnorm {f_\epsilon} _k\leq \prod_{\epsilon\in\Cub k}\norm{f_\epsilon}_\infty.
$$
\end{proposition}
\begin{proof}[Sketch of the proof]
The first part of the result is proved in \cite[Corollary~5.2]{HK09}.

To prove the second part,   let $b\in \ell^\infty(\N^d)$ and $\bI=(I_j)_{j\in \N}$ be a F\o lner sequence such that  $b$ admits correlations along $\bI$.
 We use~\eqref{eq:phi} for  $x:=T_\bn x_0$, take the averages for $\bn\in I_j$ and exchange the limits in $j$ and in $H$; this can be achieved because of the uniform convergence in \eqref{eq:phi} (see
 \cite[Theorem~5.4]{HK09}). By an iterated use of the Cauchy-Schwarz inequality  (this is estimate (12) in \cite{HK09}) we obtain that
\begin{equation}
\label{E:direct}
\Limsup\bigl\vert \Av_\bI\, a(\bn)\, b(\bn)\bigr\vert
\leq \norm b_{\bI,k}\cdot
\prod_{\epsilon\in\Cub k}\nnorm{f_{\epsilon}}_k.
\end{equation}
Taking the supremum over all F\o lner sequences $\bI$ in the left hand side, we obtain the announced bound.
For $d=1$ the details can be found in \cite[Section~5.4]{HK09}, the proof for general $d\in \N$ is similar.
\end{proof}
Let  $(X,T_1,\dots,T_d)$ be a $(k-1)$-step nilsystem defining the nilsequence $\psi\colon \N^d\to \C$.
By~\cite[Proposition~5.6]{HK09} (see also~\cite[Proposition~3.2]{CFH11}), the linear span of
the functions defined as in~\eqref{eq:phi} is dense in $C(X)$ with the uniform norm. By Proposition~\ref{prop:dual-antiunif}, the sequence $\psi$ is a uniform limit of strongly $k$-anti-uniform sequences. We deduce:

\begin{corollary}
\label{cor:direct}
 Every $(k-1)$-step nilsequence is $k$-anti-uniform.
\end{corollary}
\begin{remark}
 Alternatively,  this follows by combining Proposition~\ref{prop:represent} and Proposition~\ref{prop:unif-in-norm}.
\end{remark}

For $d=1$, the first statement of the next Proposition is~\cite[Theorem 2.16]{HK09} and the second statement is~\cite[Theorem~2.1]{F15}.

\begin{proposition}
\label{prop:unif-approx2}
Let $d,k\in \N$, $a\in\ell^\infty(\N^d)$ be a sequence,  and $\delta>0$. Then there exists a F\o lner sequence ${\bf I}=(I_j)_{j\in \N}$ such that the following holds:
\begin{enumerate}
\item
\label{it:unif-approx2a}
There exists a $(k-1)$-step nilsequence $\psi_2\colon \N^d\to\C$  such that
$$
\norm{\psi_2}_{U^k(\N^d)}^*\leq 1\ \text{ and }\
\Limsup\bigl\vert\Av_\bI\,a(\bn) \, \psi_2(\bn)\bigr\vert\geq \norm a_{U^k(\N^d)}-\delta.
$$
 Moreover, if  $b\in \ell^\infty(\N^d)$ is a sequence that  admits correlations along a F\o lner sequence $\bJ$,
then we have 
\begin{equation}\label{E:unifI}
\Limsup\bigl\vert\Av_\bJ\, \psi_2(\bn)\, b(\bn) \bigr\vert \leq  \norm{b}_{\bJ,k}.
\end{equation}
\item
\label{it:unif-approx2b}
If $\norm a_\infty\leq 1$, then there exists a $(k-1)$-step nilsequence $\psi_3\colon \N^d\to \C$ such that
$$
\norm{\psi_3}_\infty\leq 1\ \text{ and } \
\Limsup\bigl\vert\Av_\bI\, a(\bn)\, \psi_3(\bn)\bigr\vert\geq \norm a_{U^k(\N^d)}^{2^k}-\delta.
$$
\end{enumerate}
\end{proposition}
\begin{proof}
Let the $(k-1)$-step nilsequence $\psi_1\colon \N^d\to \C$ and the F\o lner sequence
$\bI$ be given by Proposition~\ref{prop:approx} for an $\ve>0$ that will be specified later.  Then
 $a-\psi_1$ admits correlations along  $\bI$  and
\begin{equation}\label{E:39}
\norm{a-\psi_1}_{\bI,k}\leq \ve, \quad  \norm{\psi_1}_{\bI,k}\geq  \norm{a}_{U^k(\N^d)}-2\ve, \quad \text{and} \quad  \norm{\psi_1}_\infty\leq \norm{a}_\infty.
\end{equation}
 The sequence $\psi_1$ has the form  $\psi_1(\bn)=\Phi_1(T_\bn x_0), \bn\in\N^d,$ for some $(k-1)$-step ergodic nilsystem $(X,\mu,T_1,\ldots, T_d)$, point $x_0\in X$,  and function $\Phi_1\in C(X)$ with  $\norm{\Phi_1}_\infty\leq \norm{a}_\infty$.
 For $x\in X$, we  define the function $\Phi$ on $X$ as in~Proposition~\ref{prop:dual-antiunif} taking
 $f_\uepsilon := \conj^{|\uepsilon|}\Phi_1$ for  $\uepsilon\in\Cub k$. Then $\Phi\in C(X)$ and  $\int \Phi\,\Phi_1\,d\mu=\nnorm{\Phi_1}_{\mu,k}^{2^k}$. Let $\psi_3\colon \N^d\to \C$ be the nilsequence defined by
 $$
 \psi_3(\bn):=\Phi(T_\bn x_0),\quad  \bn\in\N^d.
 $$
If we assume that
 $\norm a_\infty\leq 1$, then we have $\norm{\Phi_1}_\infty\leq 1$, thus $\norm\Phi_\infty\leq 1$ and $\norm{\psi_3}_\infty\leq 1$.
 By unique ergodicity of $(X,\mu,T_1,\dots,T_d)$ we have
\begin{multline*}
\Limav{\bI} \psi_1(\bn) \, \psi_3(\bn)=\int\Phi_1\, \Phi \,d\mu=\nnorm{\Phi_1}_{\mu,k}
^{2^k}
=\norm{\psi_1}_{\bI,k}^{2^k} \quad \text{by } \eqref{E:psi}\\
\geq(\norm a_{U^k(\N^d)}-2\ve)^{2^k}
\geq \norm a_{U^k(\N^d)}^{2^k}-\delta/2
\end{multline*}
by \eqref{E:39} if $\ve$ is chosen sufficiently small. Moreover, since $a-\psi_1$ admits correlations along $\bI$,  by \eqref{E:direct} we have
$$
 \Limsup\bigl\vert\Av_\bI\,  \psi_3(\bn)\, (a-\psi_1)(\bn) \bigr\vert \leq  \norm{a-\psi_1}_{\bI,k} \,
 \norm{\Phi_1}_{\infty}^{2^k-1} \leq \ve \leq \delta/2
 $$
 by \eqref{E:39} if $\ve \leq \delta/2$.
Combining the last two estimates we get   Part~\eqref{it:unif-approx2b}.

We move now to the proof of  Part~\eqref{it:unif-approx2a}.
 Notice first that using
 Proposition~\ref{prop:dual-antiunif} and \eqref{E:psi} we get that
 $$
 \norm{\psi_3}_{U^k(\N^d)}^*\leq\nnorm{\Phi_1}_{\mu,k}^{2^k-1}= \norm{\psi_1}_{U^k(\N^d)}^{2^k-1}.
  $$
  Moreover, if  the sequence  $b\in \ell^\infty(\N^d)$   admits correlations along a F\o lner sequence $\bJ$, then   \eqref{E:direct} and \eqref{E:psi}   give
 $$
 \Limsup\bigl\vert\Av_\bJ\,  \psi_3(\bn)\, b(\bn) \bigr\vert \leq
 \norm{b}_{\bJ,k} \,
 \nnorm{\Phi_1}_{\mu,k}^{2^k-1}=
  \norm{b}_{\bJ,k} \,
 \norm{\psi_1}_{U^k(\N^d)}^{2^k-1}.
 $$
  Defining
  $$
  \psi_2:=(\norm{\psi_1}_{U^k(\N^d)}^{2^k-1})^{-1}\cdot \psi_3
  $$
  we deduce
 that $\norm{\psi_2}_{U^k(\N^d)}^*\leq 1$ and   estimate  \eqref{E:unifI} is satisfied. Furthermore,  as before we get
 $$
\Limav{\bI} \psi_1(\bn) \, \psi_2(\bn)=\norm{\psi_1}_{U^k(\N^d)}^{2^k}/ \norm{\psi_1}_{U^k(\N^d)}^{2^k-1}
=\norm{\psi_1}_{U^k(\N^d)}\geq\norm a_{U^k(\N^d)}-\delta/2
$$
by \eqref{E:39}  if $\ve\leq \delta/4$.
 Using \eqref{E:direct}  and \eqref{E:psi} we get
$$
 \Limsup\bigl\vert\Av_\bI\,  \psi_2(\bn)\, (a-\psi_1)(\bn) \bigr\vert \leq  \norm{a-\psi_1}_{\bI,k} \,
 \frac{\nnorm{\Phi_1}_{\mu,k}^{2^k-1}}{ \norm{\psi_1}_{U^k(\N^d)}^{2^k-1}}=\norm{a-\psi_1}_{\bI,k}\leq \ve \leq \delta/2
$$
 by \eqref{E:39}  if $\ve \leq \delta/2$.  Combining the last two estimates finishes the proof of    Part~\eqref{it:unif-approx2a}.
\end{proof}

\begin{corollary}
\label{cor:unif-approx2}
Let $d,k\in \N$ and  $a\in\ell^\infty(\N^d)$ be such that the averages of $a(\bn)\, \psi(\bn)$ converge to $0$ for every $(k-1)$-step nilsequence $\psi\colon \N^d\to \C$. Then
$\norm a_{U^k(\N^d)}=0$.
\end{corollary}

\subsection{Regular sequences and their structure}
 Next we introduce certain  classes of sequences for which we  are able to prove the two decomposition results of this section.
\begin{definition}
Let $d\in \N$, $k\in \Z_+$. We say that a sequence $a\in\ell^\infty(\N^d)$ is \emph{$k$-regular} if the limit
$$
\limav a(\bn)\, \psi(\bn)
$$
exists for every $k$-step nilsequence $\psi\colon \N^d\to \C$.
\end{definition}
\begin{remarks}  $\bullet$ Every nilsequence is $k$-regular for every $k\in \N$.

$\bullet$ The product of two $k$-regular sequences may not be $k$-regular.\footnote{Let $a(n):=e(n^{k+1}\alpha),$  $n\in \N$, where  $\alpha\in \R\setminus \Q$ and $b:=\overline{a}\cdot \sum_{k=0}^\infty {\bf 1}_{[2^{2k}, 2^{2k+1})}$.
The sequences $a,b$ are
 $k$-regular for every $k\in \N$, but the sequence  $a\cdot b$ is not even $0$-regular since the averages
 $\frac{1}{N}\sum_{n=1}^N a(n)\, b(n)$ do not converge as $N\to +\infty$.}
\end{remarks}

\begin{theorem}[\mbox{\cite[Theorem 2.19]{HK09}} for $d=1$]
\label{th:regular}
Let $d,k\in \N$ and $a\in \ell^\infty(\N^d)$ be a sequence. Then
the following are equivalent:
\begin{enumerate}
\item\label{it:converges}
$a$ is $(k-1)$-regular.
\item\label{it:unif}
For every $\delta >0$, the sequence $a$ can be written as
$a=\psi+u$, where $\psi$ is a $(k-1)$-step nilsequence with $\norm{\psi}_\infty\leq \norm{a}_\infty$,
 and
$\norm{u}_{U^k(\N^d)}\leq\delta$.
\end{enumerate}
\end{theorem}
\begin{proof}
 The implication $\eqref{it:unif} \Longrightarrow  \eqref{it:converges}$ is a simple consequence of Corollary~\ref{cor:direct} and the fact that the product of two nilsequences is a nilsequence and thus has convergent averages.

We establish now the
  converse implication.
 For $\delta>0$ let $\psi:=\psi_1$ and $\bI$ be as in Proposition~\ref{prop:approx}  with $\ve:=\delta/2$. Then we have
 $\norm{a-\psi}_{\bI,k}\leq \delta/2$.
   Since $a=\psi+(a-\psi)$,  it remains to show that  $\norm{a-\psi}_{U^k(\N^d)}\leq\delta$. Suppose that this is not the case.
 Then
 Part~\eqref{it:unif-approx2a} of Proposition~\ref{prop:unif-approx2} provides a F\o lner sequence $\bJ$ and a $(k-1)$-step nilsequence $\psi_2$ with
 $$
\Limsup\bigl\vert\Av_\bJ (a(\bn)-\psi(\bn))\, \psi_2(\bn)\bigr\vert
\geq \norm{a-\psi}_{U^k(\N^d)}-\delta/3\geq  2\delta/3.
$$
Since $a$ is a $(k-1)$-regular sequence and  $\psi \psi_2$ is a $(k-1)$-step nilsequence,
 the sequence
$(a-\psi)\psi_2$  has convergent averages. We deduce that
$$
|\Limav{\bI} (a(\bn)-\psi(\bn))\, \psi_2(\bn)|=|\Limav{\bJ} (a(\bn)-\psi(\bn))\, \psi_2(\bn)|\geq 2\delta/3.
$$
On the other hand, we have
$$
|\Limav{\bI} (a(\bn)-\psi(\bn))\, \psi_2(\bn)|\leq \norm{a-\psi}_{\bI,k}\leq \delta/2
$$
where   the second estimate follows from our data and the first estimate follows from the estimate \eqref{E:unifI}
since the sequence $a-\psi$ admits correlations along $\bI$.
Combining the above we get a contradiction.
\end{proof}

 \subsection{The structure of regular anti-uniform sequences}

The next result gives a characterization of regular anti-uniform sequences in $\ell^\infty(\N^d)$.
Assuming the multiparameter  inverse theorem of Proposition~\ref{prop:unif-approx2},
its proof is identical to the proof of Theorem~1.2 in  \cite{F15} where the case $d=1$ was treated; we only sketch the main idea of the proof below.
\begin{theorem}[\mbox{\cite[Theorem~1.2]{F15}} for $d=1$]
\label{th:regular+antiunif}
Let $d, k \in \N$ and   $a\in\ell^\infty(\N^d)$ be a sequence. Then  the following properties are equivalent:
 \begin{enumerate}
 \item
\label{it:regul+anti}
$a$ is $(k-1)$-regular and $k$-anti-uniform.
\item $a$ is an approximate  $(k-1)$-step nilsequence.
\end{enumerate}
\end{theorem}
\begin{proof}[Idea of the proof]
The implication $\eqref{it:unif} \Longrightarrow  \eqref{it:converges}$  follows from Corollary~\ref{cor:direct} and the fact that a product of two nilsequences is a nilsequence, hence a  regular sequence.

We explain now the main idea of  the  converse implication.
 Let $a\in\ell^\infty(\N^d)$ be $(k-1)$-regular and $k$-anti-uniform. Let $\CN$  be the linear space of $(k-1)$-step nilsequences in $d$ variables and let $\CH$ be the linear span of $\CN$ and $a$. Then it follows from our $(k-1)$-regularity assumption and Theorem~\ref{th:regular} that  for all $c,c'\in\CH$ the averages of $c(\bn)\, \overline{c'(\bn)}$ converge and we write $\langle c, c'\rangle$ for this limit. We remark that $\norm c_2^2=\langle c, c\rangle$.

If the space $\CH$ endowed with the ``scalar product'' $\langle\cdot, \cdot\rangle$ was a Hilbert space and  $\CN$ was a closed subspace, then we could define $\psi$ to be the orthogonal projection of $a$ on $\CN$. Then $\langle a-\psi, \psi'\rangle=0$ for all $\psi'\in \mathcal{N}$ and Corollary~\ref{cor:unif-approx2}  would imply that  $\norm{a-\psi}_{U^k(\N^d)}=0$. Since, by assumption, $a$ is
 $k$-anti-uniform, we deduce that  $\norm{a-\psi}_2=0$, that is, $a$ is a $(k-1)$-step nilsequence plus a sequence that converges to $0$ in uniform density. In our present setup there is lack of completeness, so we choose $\psi$ to be an ``approximate orthogonal projection'', in the sense  that  $\norm{a-\psi}_2$ is sufficiently close to the distance of $a$ to $\CN$, and we obtain using Part~$\text{(i)}$ of Proposition~\ref{prop:unif-approx2} the announced decomposition. For the details see the proof of  Theorem~\cite[Theorem~1.2]{F15}  which contains a proof for $d=1$ under the assumption of strong $k$-anti-uniformity (it is called $k$-anti-uniformity there); the same argument works without change for general $d\in \N$ under the weaker assumption of $k$-anti-uniformity.
\end{proof}


\section{Correlations are regular sequences}
The goal of this section is to show that modulo small $\ell^\infty$-errors, nilsequences can be represented as multiple correlation sequences, and then use the  mean convergence result of Walsh (Theorem~\ref{T:Walsh}) in order to show that multiple correlation sequences are regular sequences.
\subsection{Producing nilsequences as correlations}
The argument we use below is analogous to the one used in \cite{F15} to handle single variable nilsequences.
\begin{lemma}[\mbox{\cite[Lemma~14.2]{GT08}}]
\label{lem:prog}
Let $d,k\in \N$ and $X=G/\Gamma$ be a $(k-1)$-step nilmanifold. Then there exists a continuous map $P\colon X^k\to X$ such that
\begin{equation}
\label{eq:prog}
P(hg\cdot e_X, h^2g\cdot e_X,\dots,h^kg\cdot e_X)=g\cdot e_X
\ \text{ for all }g,h\in G.
\end{equation}
\end{lemma}
\begin{remark} The result in \cite[Lemma~14.2]{GT08} gives $P(g\Gamma,hg\Gamma, h^2g\Gamma, \ldots, h^{\ell-1} g\Gamma)=h^{\ell} g\Gamma$.
Inserting $h^{-\ell}g $ in place of $g$, then $h^{-1}$ in place of $h$, and rearranging coordinates,  we get  \eqref{eq:prog}.
\end{remark}
\begin{proposition}
\label{prop:represent}
Let $d,k\in \N$ and $\psi \colon \N^d\to \C$ be a $(k-1)$-step nilsequence. Then for every $\ve>0$ there exists a system $(Y,\nu,S_1,\dots,S_d)$ and functions $F_1,\dots,F_k\in L^\infty(\nu)$, such that the sequence $b\colon \N^d\to \C$ defined by
\begin{equation}
\label{eq:def-b}
b(\bn):=\int \prod_{j=1}^k
\bigl(\prod_{i=1}^d S_i^{\ell_jn_i}\bigr)F_j\,d\nu, \quad \bn=(n_1,\ldots,n_d) \in \N^d,
\end{equation}
where $\ell_j=k!/j$ for $j=1,\dots,k$, satisfies
$$
\norm{\psi-b}_{\ell^\infty(\N^d)}\leq\ve.
$$
\end{proposition}
\begin{proof}
The sequence $\psi\colon \N^d\to \C$ has the form
$$
\psi(\bn)=\Psi\big(\prod_{i=1}^d\tau_i^{n_i}\cdot e_X\big), \quad \bn=(n_1,\ldots, n_d)\in \N^d,
$$
for some  $(k-1)$-step nilmanifold $X=G/\Gamma$, commuting elements
$\tau_1,\dots,\tau_d$  of $G$, and function
$\Psi\in C(X)$.
Let $\ve>0$.
As remarked  after the definition of the nilsequence in Section~\ref{SS:nil}, we can assume that the group $G$
is connected.

For $i=1,\dots,d$, let   $g_{i}\in G$ be  such that $g_{i}^{k!}=\tau_i$ (such elements exist since $G$ is connected, hence  divisible).  Let $\bm=(m_1,\ldots,m_d)\in\N^d$.
Using Lemma~\ref{lem:prog}
with  $g:=\prod_{i=1}^d g_i^{k!n_i}$ and $h:=\prod_{i=1}^d g_i^{m_i}$ and writing $H=\Psi\circ P$, we have  $H\in C(X^k)$ and  we obtain
$$
\psi(\bn)=H\Bigl(\Big(\big(\prod_{i=1}^d g_i^{m_i+k!n_i}\big), \big(\prod_{i=1}^d g_i^{2m_i+k!n_i}\big),\ldots, \big(\prod_{i=1}^dg_i^{k m_i+k!n_i}\big)\Big)\cdot e_{X^k}\Bigr)
$$
for all $\bm$ and $\bn\in\N^d$. Letting
$$
\wt{\alpha}_i:=(g_i,g_i^2,\dots,g_i^k)\in  G^k, \quad i=1,\ldots, d,
$$
 averaging with respect to $\bm\in \N^d$, and using the equidistribution results for sequences on nilmanifolds (see  \cite{Le91} or  \cite{L05}),
 we get that for every $\bn\in\N^d$ we have
\begin{align*}
\psi(\bn)
&=\lim_{M\to+\infty}\frac 1{M^d}\sum_{\bm\in  [M]^d}
H\Bigl(\Big(\big(\prod_{i=1}^dg_i^{k!n_i}\big), \big(\prod_{i=1}^d g_i^{k!n_i}\big),\ldots, \big(\prod_{i=1}^d g_i^{k!n_i}\big)\Big)\cdot \wt{\alpha}_1^{m_1}\cdot\ldots\cdot\wt{\alpha}_d^{m_d}\cdot e_{X^k}\Bigr)
\\
&
=
\int_{\wt Y}
H\Bigl(
\bigl(\prod_{i=1}^d g_i^{k!n_i}\bigr)\cdot x_1\,,\,
\bigl(\prod_{i=1}^dg_i^{k!n_i}\bigr)\cdot x_2\,,\,\dots\,,\,
\bigl(\prod_{i=1}^dg_i^{k!n_i}\bigr)\cdot x_k\Bigr)\, dm_{\wt Y}(x_1,\dots,x_k)
\end{align*}
where  $\wt Y$ is the closure of the sequence $\{\wt\alpha_1^{m_1}\cdot\ldots\cdot \wt\alpha_d^{m_d}\cdot e_{X^k} \colon m_1,\ldots, m_d \in\N\}$ in $X^k$ and $m_{\wt Y}$ is  the  Haar measure of this sub-nilmanifold of $X^k$.
For $i=1,\dots,d,$ let $\wt{S_i}\colon \wt Y\to \wt Y$ be the translation by $\wt\alpha_i$.
Note that  $(\wt Y,m_{\wt Y},\wt{S}_1,\dots,\wt{S}_d)$ is a nilsystem.

The continuous function $H$ on $X^k$ can be approximated uniformly by linear combinations of functions of the form $f_1\otimes\dots\otimes f_k$ where $f_j\in C(X)$ for $j=1,\dots,k$. Since finite linear combinations of sequences of the form \eqref{eq:def-b} have the same form (see the proof of Theorem~\ref{T:4} in Section~\ref{subsec:proofT4} below),  it remains to show that any sequence $\psi'\colon \N^d\to \C$, given by
$$
\psi'(\bn):=\int_{\wt Y}\prod_{j=1}^k f_j\Bigl(\bigl(
\prod_{i=1}^dg_i^{k!n_i}\bigr)\cdot x_j\Bigr)\, dm_{\wt Y}(x_1,\dots,x_k), \quad \bn=(n_1,\ldots, n_d)\in \N^d,
$$
has the form \eqref{eq:def-b}. To this end,
for $j=1,\dots,k,$ let $ F_j\in C(X^k)$ be given by  $ F_j(\wt x):=f_j(x_j)$ for $\wt x=(x_1,\dots,x_k)\in X^k$.
Recall that $\ell_j=k!/j$ for $j=1,\dots,k$. Since
$\wt{S}_i\wt{x}=\wt{\alpha}_i\wt{x}=(g_i,g_i^2,\dots,g_i^k)\wt{x}$, $i=1,\ldots, d$,
the $j^{\text{th}}$ coordinate of the element $(\prod_{i=1}^d \wt S_i^{\ell_jn_i})\wt x$ is $\prod_{i=1}^dg_i^{j\ell_jn_i}\cdot x_j=\prod_{i=1}^d g_i^{k!n_i}\cdot x_j$
and thus for $j=1,\ldots, k$
we have
$$F_j\bigl((\prod_{i=1}^d \wt S_i^{\ell_jn_i})\wt x\bigr)=
f_j\big(\prod_{i=1}^dg_i^{k!n_i}\cdot x_j\big), \quad n_1,\ldots, n_d\in \N.
$$ Therefore, we have
$$
\psi'(\bn)= \int_{\wt Y}\prod_{j=1}^k
 F_j\bigl((\prod_{i=1}^d \wt S_i^{\ell_jn_i})\wt x\bigr)
\,dm_{\wt Y}(\wt x), \quad \bn=(n_1,\ldots, n_d)\in \N^d.
$$
This completes the proof.
\end{proof}

\subsection{Proof of Proposition~\ref{prop:converse}}
\label{subsec:Proof-Prop3}
Let $d,s\in\N$ and $w\in \ell^\infty(\N^d)$ satisfy the hypothesis of Proposition~\ref{prop:converse}. Let $\psi$ be an $s$-step nilsequence in $d$ variables.

We set $k:=s+1$. Let $\ve >0$ and let the system
$(Y,\nu,S_1,\dots,S_d)$ and the functions $F_1,\dots,F_k\in L^\infty(\nu)$ be as in
Proposition~\ref{prop:represent}. Letting $h_j:=F_{j+1}$, for $j=0,\ldots, s$, the sequence $b\colon\N^d\to \C$ defined by
\eqref{eq:def-b}  can be rewritten as
$$
b(\bn)=\int h_0\cdot\prod_{j=0}^{s}S_{(\ell_{j+1}-\ell_1)\cdot\bn} h_j\,d\nu, \quad \bn \in \N^d.
$$
 By hypothesis the averages of $w(\bn)\, b(\bn)$ converge. Since $|b(\bn)-\psi(\bn)|$ is uniformly bounded by $\ve$, the oscillations of the averages
of $w(\bn)\, \psi(\bn)$ are uniformly bounded by $2\ve$.
Since this holds for every $\ve>0$, the averages of $w(\bn)\, \psi(\bn)$ converge, completing the proof.
\qed

\subsection{Proof of Proposition~\ref{prop:weight-nil}}
\label{subsec:Proof-Prop1}
By Proposition~\ref{prop:represent} it suffices to prove that the limit
$$
\limav b(\bn)\, T_{\vec{p_1}(\bn)}f_1\cdot\ldots\cdot T_{\vec{p_s}(\bn)}f_s
$$
exists in $L^2(\mu)$ for every sequence $(b(\bn))_{\bn\in\N^d}$ defined as in~\eqref{eq:def-b}.

By Theorem~\ref{T:Walsh} the limit
$$
\Limav{\bn}\Big(\prod_{j=1}^kF_j\big(\prod_{i=1}^d S_i^{\ell_jn_i}y\big)\cdot \prod_{m=1}^sf_m\bigl(T_{\vec{p_m}(\bn)}x\bigr)\Big)
$$
exists in $L^2(\nu\times\mu)$. Taking the integral over $Y$ with respect to $\nu$ we obtain the announced result.
\qed

\section{The structure of systems of order $1$}
\label{sec:orth-basis}
In this section we prove a structural result which is going to be used in the proof of Theorem~\ref{th:correl-lin} in the next section.
Here we work only with systems $(X,\mu, T)$ with a single transformation.
We denote by $\CI(T)$ the $\sigma$-algebra of $T$-invariant subsets of $X$ and  write the ergodic decomposition of $\mu$ under $T$ as
$$
\mu=\int \mu_x\,d\mu(x),
$$
where for $\mu$-a.e. $x\in X$ the measure $\mu_x$ is invariant and ergodic under $T$ and  the map $x\mapsto\mu_x$ is invariant under $T$.
Let $f\in L^1(\mu)$. Then for $\mu$-a.e. $x\in X$ the function $f$ is defined $\mu_x$-a.e. and belongs to $L^1(\mu_x)$; the map $x\mapsto\int f\,d\mu_x$ is measurable with respect to $\CI(T)$, and we have
$$
\E_\mu(f\mid\CI(T))(x)=\int f\,d\mu_x \ \text{ for }\ \mu\text{-a.e. } x\in X,
$$
meaning that for every set $A\in\CI(T)$ we have$$
\int_A f\,d\mu=\int_A\Bigl(\int f\,d\mu_x\Bigr)\,d\mu(x).
$$

The factor $\CZ_1$ is defined in Appendix~\ref{subse:Zk} and reduces to the Kronecker factor for ergodic systems.

\begin{definition}
We say that  a system $(X,\mu, T)$ has \emph{order} $1$ if $\CZ_1=\CX$.
\end{definition}
Throughout this section we work only with systems of order $1$.
Generalizing the construction of a Fourier basis of a compact Abelian group, our goal is to  construct a
 ``relative orthonormal'' basis for systems of order $1$ consisting of ``relative'' eigenfunctions. This is the context of  Theorem~\ref{th:basis-eigen} below.

\subsection{Relative orthonormal basis}
\label{SS:ROB}

\begin{definition}
A \emph{relative orthonormal system} (with respect to the $T$-invariant $\sigma$-algebra $\CI(T)$) is a countable family $(\phi_j)_{j\in \N}$ of functions belonging to $L^2(\mu)$ such that
\begin{enumerate}
\item
\label{it:basis1}
$\E_\mu(|\phi_j|^2\mid\CI(T))$  has value $0$ or $1$ $\mu$-a.e. for every $j\in \N$;
\item
\label{it:basis2}
 $\E_\mu(\phi_j\overline{\phi_k}\mid\CI(T))=0$ \ $\mu$-a.e. for all  $j,k\in \N$ with $j\neq k$.
\setcounter{mycounter}{\value{enumi}}
\end{enumerate}
The family  $(\phi_j)_{j\in \N}$ is a \emph{relative orthonormal basis} if it also satisfies
\begin{enumerate}
\setcounter{enumi}{\value{mycounter}}
\item
\label{it:basis3}
the linear space spanned by all functions of the form $\phi_j\psi,$ where $j\in \N$ and $\psi\in L^\infty(\mu)$  varies over all  $T$-invariant functions, is dense in $L^2(\mu)$.
\end{enumerate}
\end{definition}
We do not make the apparently natural assumption that $\E_\mu(|\phi_j|^2\mid \CI(T))=1$ $\mu$-a.e., as  there does not exist in general a relative orthonormal basis satisfying this additional condition (consider a system with ergodic components given by rotations on cyclic groups of different order).  This creates a few minor complications in the statements and the  proofs below.
We remark that the definition allows some of the elements of the base to be identically $0$. This explains why we can assume without loss of generality that all relative orthonormal systems are countably infinite, and thus indexed by $\N$.

\begin{definition} Given a relative orthonormal system  $(\phi_j)_{j\in \N}$ and $f\in L^2(\mu)$ we let
$$
f_j:=\E_\mu(f\overline{\phi_j}\mid\CI(T)), \quad j\in \N.
$$
 If $(\phi_j)_{j\in \N}$ is a relative orthonormal basis, we say that the $T$-invariant functions $f_j, j\in \N,$ are the \emph{coordinates} of $f$ in this basis.
\end{definition}
\begin{example}
\label{Ex1}
 On $\T^2$ with the Haar measure $m_{\T^2}$ let $T\colon \T^2\to \T^2$ be given by $T(x,y)=(x,y+x)$.  Then $(e(jy))_{j\in \Z}$ is a relative orthonormal basis for $L^2(m_{\T^2})$.  The coordinates $(f_j)_{j\in \Z}$ of a   function $f\in L^2(m_{\T^2})$ are given by $f_j(x):=\int f(x,y)\,  e(-jy)\, dy,$ $j\in\Z$.
\end{example}
We remark that if $\norm{f}_{L^\infty(\mu)}\leq 1$, then $\norm{f_j}_{L^\infty(\mu)}\leq 1$ for every $j\in \N$.
We will use repeatedly
that condition~\eqref{it:basis1} implies the  identity
\begin{equation}\label{E:trivial}
 f_j=f_j \cdot \E_\mu(|\phi_j|^2\mid\CI(T)) \ \mu\text{-a.e.,} \quad j \in  \N;
\end{equation}
in particular, $f_j(x)=0$ $\mu$-a.e. on the set where $\E(|\phi_j|^2\mid\CI(T))(x)=0$. We remark also that for $\mu$-a.e. $x\in X$ for every $j\in \N$  we have the identities
$$
f_j(x)=\int f\, \overline{\phi_j}\, d\mu_x \quad \text{and} \quad  \E(|f|^2\mid \CI(T))(x)=\norm{f}_{L^2(\mu_x)}^2.
$$

Next we establish a relative version of the Parseval and Fourier identity and give  necessary and sufficient conditions for a relative orthonormal system to be a relative orthonormal basis.
\begin{proposition}
\label{prop:relative}
Let $(X,\mu,T)$ be a system with ergodic decomposition $\mu=\int \mu_x\, d\mu(x)$ and  $(\phi_j)_{j\in  \N}$ be a relative orthonormal system for $L^2(\mu)$. Then
\begin{enumerate}
\item
\label{it:relative1}
For every $f\in L^2(\mu)$ and  for $\mu$-a.e.  $x\in X$ the series
\begin{equation}
\label{eq:relative1}
Pf:= \sum_{j\in\N} f_j\phi_j
\end{equation}
converges in $L^2(\mu_x)$ and we have
\begin{equation}
\label{eq:relativeb}
\norm{Pf}_{L^2(\mu_x)}^2=\sum_{j\in \N}|f_j(x)|^2\ \text{ and }\ \norm f_{L^2(\mu_x)}^2=\norm{Pf}_{L^2(\mu_x)}^2+\norm{f-Pf}_{L^2(\mu_x)}^2.
\end{equation}
\item
\label{it:relative2}
$(\phi_j)_{j\in  \N}$ is a relative orthonormal basis if and only if for every $f\in L^2(\mu)$ and for $\mu$-a.e. $x\in X$ the series~\eqref{eq:relative1} converges to $f$ in $L^2(\mu_x)$. In this case, we have
\begin{equation}
\label{eq:relative2}
\norm f_{L^2(\mu_x)}^2=\sum_{j\in\N}|f_j(x)|^2\quad \mu\text{-a.e..}\end{equation}

\item
\label{it:relative3}
$(\phi_j)_{j\in  \N}$ is a relative orthonormal basis if and only if for every $f\in L^2(\mu)$  we have
\begin{equation}
\label{eq:relative3}
\norm f_{L^2(\mu)}^2=\sum_{j\in\N} \norm{f_j}_{L^2(\mu)}^2.
\end{equation}
\end{enumerate}
\end{proposition}

\begin{proof}
We prove \eqref{it:relative1}.
Let $f\in L^2(\mu)$;  then $f\in L^2(\mu_x)$ for $\mu$-a.e. $x\in X$.
 Recall that $\int|\phi_j|^2\,d\mu_x=\E_\mu(|\phi_j|^2\mid\CI(T))(x)\in\{0,1\}$ for $\mu$-a.e. $x\in X$.
Let $E_x=\{j\in\N\colon \int|\phi_j|^2\,d\mu_x=1\}$. Then for $\mu$-a.e.  $x\in X$ the functions $(\phi_j)_{j\in E_x}$ form an orthonormal system in the Hilbert space $L^2(\mu_x)$ and viewing $f$ as an element of $L^2(\mu_x)$, the functions  $(f_j)_{ j\in E_x}$ are the coordinates of $f$ in this system. Note also that   by~\eqref{E:trivial}, we have $f_j=0$ if $j\notin E_x$. This establishes~\eqref{it:relative1}.

We prove \eqref{it:relative2}. Suppose first that for every $f\in L^2(\mu)$  for $\mu$-a.e. $x\in X$ we have $Pf=f$  $\mu_x$-a.e.. It follows that the series in~\eqref{eq:relative1} converges  to $f$ in $L^2(\mu)$. Therefore, $f$ belongs to the closed linear subspace of $L^2(\mu)$ spanned by all functions of the form $\phi_j\psi$ where $j\in \N$ and $\psi\in L^\infty(\mu)$ varies over all $T$-invariant functions. By definition, $(\phi_j)_{j\in  \N}$ is a relative orthonormal system.

We establish now the converse implication.  Suppose that  $(\phi_j)_{j\in  \N}$ is a relative orthonormal basis and let $f\in L^2(\mu)$. For $\mu$-a.e. $x\in X$, the function $Pf$ defined by~\eqref{eq:relative1} satisfies $\int Pf\, \overline{\phi_j}\,d\mu_x=f_j$ for every $j\in \N$, hence $\E_\mu((f-Pf)\, \overline{ \phi_j}\mid\CI(T))=0$ and $\int (f-Pf)\, \overline{ \phi_j}\, \overline\psi\,d\mu=0$ for every $T$-invariant function $\psi\in L^\infty(\mu)$. Therefore, the  function $f-Pf$ is orthogonal in $L^2(\mu)$ to the linear space spanned by all functions of the form $\phi_j\psi$ where $j\in \N$ and $\psi\in L^\infty(\mu)$ varies over all   $T$-invariant functions. Thus, $f-Pf=0$ by hypothesis.
Inserting $f$ in place of $Pf$ in \eqref{eq:relativeb}  gives identity~\eqref{eq:relative2}.

We prove \eqref{it:relative3}. If $(\phi_j)_{j\in  \N}$ is a relative orthonormal basis, then
\eqref{eq:relative3} follows by integrating the  identity in \eqref{eq:relative2} over $X$ with respect to $\mu$.
To prove the converse implication, let $f\in L^2(\mu)$. Using \eqref{eq:relative3}, inserting the first identity of \eqref{eq:relativeb} in the second,  and
integrating over $X$ with respect to $\mu$, we  deduce that $f=Pf$ for $\mu$-a.e. $x\in X$. Hence, Property~\eqref{it:basis3} of the definition of a relative orthonormal basis
is satisfied.
\end{proof}
\subsection{Relative orthonormal basis of eigenfunctions}
\begin{definition}
Let $\lambda\in L^\infty(\mu)$ be a $T$-invariant function  and $\phi\in L^\infty(\mu)$. We say that $\phi$ is an \emph{eigenfunction} with eigenvalue $\lambda$ if
\begin{enumerate}
\item
\label{it:eigen1}
 $|\phi(x)|$ has value $0$ or $1$ for $\mu$-a.e. $x\in X$;
\item
\label{it:eigen2}
$\lambda(x)=0$ for $\mu$-a.e. $x\in X$ such that $\phi(x)=0$;
\item
\label{it:eigen3}
$\phi\circ T=\lambda\cdot\phi$ $\mu$-a.e..
\end{enumerate}
\end{definition}
The role of Property~\eqref{it:eigen2} is to avoid ambiguities. Note also that  Property~\eqref{it:eigen3} does not imply anything about the value of $\lambda(x)$ at the points $x\in X$ where $\phi(x)=0$.

 Property~\eqref{it:eigen3} gives that   for $\mu$-a.e. $x\in X$ such that $\phi(x)\neq 0$ we have $\phi(T\inv x)\neq 0$ and thus $|\phi(x)|=|\phi(T^{-1}x)|=1$ by Property $\text{(i)}$ above and $|\lambda(x)|=1$
by Property $\text{(iii)}$.
On the other hand,  for $\mu$-a.e. $x\in X$ such that $\phi(x)=0$ we have $\phi(T\inv x)= 0$ by Property $\text{(iii)}$
and $\lambda(x)=0$ by Property $\text{(ii)}$.
Therefore, the function $|\phi|$ is $T$-invariant and
\begin{equation}
\label{eq:eigen4}
|\phi|=|\lambda|.
\end{equation}

Next, we state  the main result of this section and  we   prove it  in Section~\ref{sec:proof}.

\begin{theorem}
\label{th:basis-eigen}
Let $(X,\mu,T)$  be a system of order $1$. Then  $L^2(\mu)$ admits a relative orthonormal basis of eigenfunctions.
\end{theorem}
\begin{remarks}
$\bullet$
It is  true and not difficult to prove that if a system has a relative orthonormal basis of eigenfunctions, then it has   order $1$.

$\bullet$ For ergodic systems, Theorem~\ref{th:basis-eigen} is well known, but it is not easy to deduce the general case from the ergodic one. The reason is that although
the ergodic components of a system of order $1$ are ergodic rotations (Proposition~\ref{prop:ae-rotation}), we cannot  simply ``glue'' together their eigenfunctions, because of
 measurability issues.
\end{remarks}

\begin{example}
For the system described in the Example~\ref{Ex1} we have that  $(e(jy))_{j\in \Z}$ is a relative orthonormal basis of eigenfunctions with eigenvalues $(e(jx))_{j\in \Z}$.
\end{example}

The next proposition will be used in the proof of Theorem~\ref{th:correl-lin}.
\begin{proposition}
\label{prop:basis-eigen2}
Let $(X,\mu,T)$ be a system of order $1$ with ergodic decomposition $\mu=\int \mu_x\, d\mu(x)$. Suppose that  $(\phi_j)_{j\in \N}$ is an orthonormal basis of eigenfunctions,  $f\in L^\infty(\mu)$,  and let $(f_j)_{j\in \N}$  be the coordinates of $f$ in this base. Then  for $\mu$-a.e. $x\in X$ we have
\begin{equation}
\label{eq:basis-eigen2}
\nnorm f_{T,\mu_x, 2}^4=\sum_{j\in \N}|f_j(x)|^4 \quad \mu\text{-a.e.}.
\end{equation}
\end{proposition}
\begin{remark}
After integrating \eqref {eq:basis-eigen2} over $X$ with respect to $\mu$ it follows from \eqref{E:seminonerg} that
$\nnorm f_{T,\mu, 2}^4=\sum_{j\in \N}\norm{f_j}_{L^4(\mu)}^4$.
\end{remark}

In the proof of Proposition~\ref{prop:basis-eigen2} we will use the following basic fact:
\begin{lemma}
\label{E:different}
Let $(X,\mu,T)$ be a system, $(\phi_j)_{j\in \N}$ be a relative orthonormal system of eigenfunctions,
and $(\lambda_j)_{j\in \N}$ be the corresponding eigenvalues. Then for $\mu$-a.e. $x\in X$  we have $\lambda_j(x)\overline{\lambda_k(x)}\neq 1$ for  all $j,k\in \N$ with $j\neq k$.
\end{lemma}
\begin{proof}[Proof of Lemma~\ref{E:different}]
Note that  the set $A=\{x\colon \lambda_j(x)\overline{\lambda_k(x)}= 1\}$   is $T$-invariant  and the function $\one_A\phi_j\overline{\phi_k}$ is $T$-invariant by Part~\eqref{it:eigen3} of the definition of an eigenfunction, and thus equal to $0$ by Part~\eqref{it:basis2} of the definition of a relative orthonormal system. On the other hand, for $\mu$-a.e.  $x\in A$ we have
$|\phi_j(x)\overline{\phi_k(x)}|=1$ by~\eqref{eq:eigen4}. Thus, $\mu(A)=0$
and the claim follows.
\end{proof}

\begin{proof}[Proof of Proposition~\ref{prop:basis-eigen2}]
By Part~\eqref{it:relative2} of Proposition~\ref{prop:relative}  
we have for $\mu$-a.e.
 $x\in X$ that
$$
f=\sum_{j\in \N}f_j \cdot\phi_j
$$
where $f_j =\int f\overline{\phi_j}\, d\mu_x$ and the convergence takes place in $L^2(\mu_x)$. 
 It follows that for
$\mu$-a.e. $x\in X$ for every $n\in\N$ we have
 $$
 T^nf\cdot\overline f=\sum_{j,k\in \N}\lambda_j^n\, f_j\, \overline{ f_k}\,  \phi_j\, \overline{\phi_k}
$$
where  convergence takes place in $L^1(\mu_x)$.
Using this, identity \eqref{E:trivial}, and that    $\int \phi_j\overline{\phi_k}d\mu_x=0$ for   $j,k\in \N$ with $j\neq k$, we deduce that for  $\mu$-a.e. $x\in X$ we have
$$
 \int T^nf\cdot\overline f\, d\mu_x=\sum_{j\in \N} \lambda_j^n|f_j|^2.
$$
Note that by \eqref{eq:relative2} the  above series converges absolutely for $\mu$-a.e. $x\in X$.
Hence, for $\mu$-a.e. $x\in X$ we have
$$
   \Big|\int T^nf\cdot\overline f\, d\mu_x\Big|^2= \sum_{j,k\in \N} (\lambda_j\overline{\lambda_k})^n  |f_j|^2|f_k|^2.
$$
Averaging in $n\in\N$ and using  that by \eqref{eq:seminorm2} we have
$$
\nnorm f_{T,\mu_x,2}^4=\lim_{N\to +\infty}\frac{1}{N}\sum_{n=1}^N \Big|\int T^nf\cdot\overline f\, d\mu_x\Big|^2,
$$
we obtain for  $\mu$-a.e. $x\in X$ that
$$
\nnorm f_{T,\mu_x,2}^4=
\sum_{j,k\in \N} \Big(\lim_{N\to+\infty}\frac 1N\sum_{n=1}^N(\lambda_j\overline{\lambda_k})^n\Bigr)\cdot |f_j|^2|f_k|^2
$$
where the interchange of limits is justified  because $\sum_{j,k\in \N}|f_j|^2|f_k|^2$ converges $\mu$-a.e. by \eqref{eq:relative2}.

If $j\neq k$, since $|\lambda_j\overline{\lambda_k}|\leq 1$ and
 $\lambda_j\overline{\lambda_k}\neq 1$ by Lemma~\ref{E:different},   we have for $\mu$-a.e. $x\in X$ that
$$
\lim_{N\to+\infty}\frac 1N\sum_{n=1}^N(\lambda_j\overline{\lambda_k})^n=0.
$$
Suppose now that $j=k$. For $\mu$-a.e. $x\in X$ the following holds:
If $\lambda_j(x)=0$, then $\phi_j(x)=0$ by~\eqref{eq:eigen4} and $f_j(x)=0$ by~\eqref{E:trivial}. If $\lambda_j(x)\neq 0$, then
$|\lambda_j(x)|=1$ by~\eqref{eq:eigen4}  and Part~\eqref{it:eigen1}
of the definition of an eigenfunction. In both cases we have $|\lambda_j(x)|^{2n} |f_j(x)|^4=|f_j(x)|^4$. Combining the above we get
\eqref{eq:basis-eigen2}.
\end{proof}

\subsection{A Borel selection result}
The proof of Theorem~\ref{th:basis-eigen} needs some technical preliminaries.
We will need  the following selection theorem of Lusin-Novikov.
\begin{theorem}[see for example \mbox{\cite[Theorem 18.10]{Ke95}}]
\label{T:Lusin}
Let $X,Y$ be Polish (i.e. complete separable  metric) spaces and $P\subset X \times Y $ be a
Borel set such that every vertical section $P_x= \{y \colon \, (x,y) \in
P\}$ is a countable set. Then the vertical projection $A$ of $P$ on $X$ is Borel and
 there exists a Borel function
$f\colon A\rightarrow Y$ such that $f(x)\in P_x$ for every $x\in A$.
\end{theorem}
\begin{proposition}
\label{prop:select}
Let  $X,K$  be Polish spaces  and $F\colon X\times K\to \R_+$ be a bounded Borel function.  Suppose that for every $x\in X$ the
set
$$
P_x:=\{y\in K\colon F(x,y)>0\}
$$
is countable and
\begin{equation}
\label{eq:Sigma}
\Sigma(x):=\sum_{y\in P_x}F(x,y)<+\infty.
\end{equation}
Let $N(x):=|P_x|\in[0,+\infty]$.
 Then there exists a sequence $(t_j)_{j\in \N}$ of Borel maps  $X\to K$ such that for  every $x\in X$ the values $t_j(x)$, $1\leq j<1+N(x)$, are pairwise distinct and
$P_x=\{t_j(x)\colon 1\leq j\leq N(x)\}$.
\end{proposition}
\begin{remarks}
$\bullet$ The  condition $j<1+N(x)$ means $j\leq N(x)$ if $N(x)<+\infty$ and $j$ is arbitrary if $N(x)=+\infty$.

$\bullet$ Note that if $P_x=\emptyset$, then the values of $t_j(x)$ are not determined by the statement.
\end{remarks}

Proposition~\ref{prop:select} is a consequence of  the next lemma.
\begin{lemma}
\label{lem:select}
Let $F, K, X$ be as in Proposition~\ref{prop:select} and for $x\in X$ let
$$
S(x):=\sup_{y\in K}F(x,y).
$$  Then there exists a Borel map $t\colon X\to K$ such that
$F(x,t(x))=S(x)$  for every $x\in X$.
\end{lemma}
\begin{proof}[Proof of Lemma~\ref{lem:select}]
Note that \eqref{eq:Sigma} implies that this supremum $S(x)$ is attained.
Let $\CX$ and $\CK$ be the Borel $\sigma$-algebras of the spaces $X$ and $K$ respectively.

We first claim that the  function $S$ is Borel.
Indeed, let $s\geq 0$. The set
$\{(x,y)\in X\times K\colon F(x,y)> s\}$ belongs to $\CX\otimes\CK$ and has countable fibers. By  Theorem~\ref{T:Lusin} its projection on $X$ belongs to $\CX$. Since this projection is the set $\{x\in X\colon S(x)>s\}$, the map $S$ is Borel.

For $m\in \N$, let
$$
A_m:=\{(x,y)\in X\times K\colon F(x,y)\geq (1-2^{-m})S(x)\}.
$$
Then $A_m$ belongs to $\CX\otimes\CK$, the projection of $A_m$ on $X$ is onto, and this projection is countable to one.  By Theorem~\ref{T:Lusin}, for every $m\in \N$ there exists a Borel map $t_m\colon X\to K$ such that $(x,t_m(x))\in A_m$ for every $x\in X$, that is,
$$
F(x,t_m(x))\geq (1-2^{-m})S(x) \ \text{ for every } x\in X.
$$
If $x$ is such that $S(x)=0$, then $t_m(x)=0$ for every $m\in \N$. If not, then  $F(x,t_m(x))\geq S(x)/2$ for every $m\in \N$ and thus the set
$\{t_m(x)\colon m\geq 1\} $ contains at most $2\Sigma(x)S(x)\inv$ distinct elements.
It follows that for every $x\in X$ the sequence $(t_m(x))_{m\in \N}$ is eventually constant. The limit value $t(x)$ of this sequence is therefore well defined, it is a Borel map from $X$ to $K$, and satisfies $F(x,t(x))=S(x)$. This completes the proof.
\end{proof}

\begin{proof}[Proof of Proposition~\ref{prop:select}]
We build by induction the family of Borel maps $t_j\colon X\to K, j\in \N$.
Let $t_1\colon X\to K$ be given by  Lemma~\ref{lem:select}. We have that $F(x,t_1(x))=0$ if and only if $P_x$ is empty, that is, if  $N(x)=0$. We define
\begin{gather*}
F'(x,t)
:=\begin{cases}
0&\text{if }t=t_1(x);\\
F(x,t)&\text{otherwise,}
\end{cases}\\
 P'_x:=\{t\colon F'(x,t)>0\} \ \text{\ for }x\in X.
\end{gather*}
The function $F'$ is Borel and for every $x\in X$ for which $P_x$ is non-empty, this set is the disjoint union of $P'_x$ and $\{t_1(x)\}$. We replace the function $F$ with $F'$, and Lemma~\ref{lem:select} provides a map $t_2\colon X\to K$.  Iterating, we obtain a sequence $(t_j)_{j\in \N}$ of Borel maps   $X\to K$ that satisfy
\begin{gather*}
F(x,t_1(x))\geq F(x,t_2(x))\geq\dots\geq F(x,t_j(x));\\
F(x,t_j(x))=0\text{ if and only if }j>N(x);\\
\text{if }N(x)\geq j\text{ then } t_{j+1}(x)\notin\{t_1(x),\dots,t_j(x)\};\\
F(x,t_{j+1}(x))=\sup\bigl\{F(x,t)\colon t\notin\{t_1(x),\dots,t_j(x)\}\bigr\}.
\end{gather*}

We have  $\{t_j(x)\colon 1\leq j\leq N(x)\}\subset P_x$ and we claim that equality holds. Suppose that this is not the case and let $t\in P_x\setminus \{t_j(x)\colon 1\leq j\leq N(x)\}$. Then, by construction, $t_j(x)\geq t$ for $j\leq N(x)$ and thus $\sigma(x)\geq N(x)t$. It follows that $N(x)$ is finite. By construction,
$|\{t_j(x)\colon 1\leq j\leq N(x)\}|=N(x)=|P_x|$ which is a contradiction. This completes the proof.
\end{proof}

\subsection{Proof of Theorem~\ref{th:basis-eigen}}
\label{sec:proof}
In this subsection we use a different presentation of the ergodic decomposition  of a system that is more convenient for our purposes.\footnote{The measures $\mu_x$ of the top of Section~\ref{sec:orth-basis} are written as $\mu_{\pi(x)}$ here.} Recall that we  assume that $(X,\CX,\mu)$ is a Lebesgue space. It is well known (see for example~\cite[Theorems~8.7 and~A.7]{G03}) that there exists a Lebesgue space $(Y,\CY,\nu)$,
 a measure preserving map $\pi\colon X\to Y$ that satisfies
$\pi\circ T=\pi$,
 $\CI(T)=\pi\inv(\CY)$ up to $\mu$-null sets,  and for $y\in Y$ a probability measure $\mu_y$ on $X$ such that
\begin{enumerate}
\item
 The map $y\mapsto\mu_y$ is Borel, meaning that  for every bounded Borel function on $X$, the function $y\mapsto\int f\,d\mu_y$ is Borel.
 \item \label{it:ergo-decomp2}
 For every bounded Borel function on $X$,
$$
\E_\mu(f\mid(\CI(T))(x)=\int f\,d\mu_{\pi(x)}\ \text{ for }\mu\text{-a.e. }x\in X.
$$

\item
For $\nu$-a.e.  $y\in Y$ the measure $\mu_y$ on $X$ is concentrated on $\pi\inv(\{y\})$ and  is  invariant and ergodic under $T$.
\setcounter{mycounter}{\value{enumi}}
\end{enumerate}

Taking the integral in~\eqref{it:ergo-decomp2}, we obtain
$$
\mu=\int_Y \mu_y\, d\nu(y).
$$  By density,
\begin{enumerate}
\setcounter{enumi}{\value{mycounter}}
\item
For $f\in L^1(\mu)$, we have $f\in L^1(\mu_y)$ for $\nu$-a.e. $y\in Y$, $\norm f_{L^1(\mu)}=\int\norm f_{L^1(\mu_y)}\,d\nu(y)$,  and the equality in~\eqref{it:ergo-decomp2} remains valid.
 \setcounter{mycounter}{\value{enumi}}
 \end{enumerate}

 Lastly, since our standing assumption in this section is that the system $(X,\mu,T)$ has order $1$, by Proposition~\ref{prop:ae-rotation}
   we have:
 \begin{enumerate}
 \setcounter{enumi}{\value{mycounter}}
\item
For $\nu$-a.e.  $y\in Y$ the system $(X,\mu_y,T)$ is an ergodic rotation.
 \end{enumerate}

We first prove the following intermediate result:
\begin{lemma}
\label{lem:eigen-ortho}
Let $(X,\mu,T)$ be a system of order $1$ and $f\in L^2(\mu)$. Then there exists a relative orthonormal system of eigenfunctions  $(\phi_j)_{j\in \N}$ such that
\begin{equation}\label{E:f}
f=\sum_{j\in \N} \E_\mu(f\overline{\phi_j}\mid\CI(T))\cdot\phi_j
\end{equation}
where convergence takes place in $L^2(\mu)$.

Moreover, for every $j\in \N$ the eigenfunction $\phi_j$ belongs to the smallest closed $T$-invariant subspace of $L^2(\mu)$ containing
the set $\{ f\phi\colon \phi\in L^\infty(\mu) \text{ is } T\text{-invariant}\}$.
\end{lemma}
\begin{remark}
For notational convenience, the orthonormal system we build below is indexed by $\Z_+$ instead of $\N$. 
\end{remark}
\begin{proof}
We can assume that $f$ is a Borel function defined everywhere.
The set $Y_0$ of $y\in Y$ such that $\mu_y$ is invariant under $T$ is Borel and has full measure. Since the map $y\mapsto\norm f_{L^2(\mu_y)}$ is Borel, and since $\int\norm f_{L^2(\mu_y)}^2\,d\nu(y)=\norm f_{L^2(\mu)}<+\infty$, the set $Y_1$ of points $y\in Y_0$ such that $f\in L^2(\mu_y)$ is Borel and has full measure. Substituting $Y_1$ for $Y$, we are reduced to the case where the measure $\mu_y$ is $T$-invariant and $f\in L^2(\mu_y)$ for every $y\in Y$. Since $(Y,\CY,\nu)$ is a Lebesgue space, we can assume that $Y$ is a Polish space and $\CY$ is its Borel $\sigma$-algebra.

For  $y\in Y$, we write  $\sigma_y$ for the spectral measure of $f$  with respect to  the system $(X,\mu_y,T)$; it is the finite positive measure on $\T$ defined by
$$
\wh{\sigma_y}(n):=\int T^nf\cdot\overline f\,d\mu_y, \quad n\in\Z.
$$
For $\nu$-a.e.  $y\in Y$ this measure is atomic because $(X,\mu_y,T)$ is a rotation. For every  $t\in\T$ and every $y\in Y$  the limit
$$
F(y,t):=\lim_{N\to+\infty}\frac 1N\sum_{n=1}^N\wh{\sigma_y}(n)\, \e(-nt)
$$
exists and we have
$$
F(y,t)=\sigma_y(\{t\}).
$$

Since for every $n\in \N$ the function $T^nf\cdot\overline f$ is Borel, we get that  the map $y\mapsto \wh{\sigma_y}(n)$ is Borel on $Y$.
Thus, the function $F$ satisfies the hypothesis of Proposition~\ref{prop:select} for the Polish space $Y\times \T$. Henceforth,  we use the notation of this proposition and
let $t_j\colon Y\to \T$, $j\in \N$, be the Borel maps obtained. For $j\in\Z_+$, let
\begin{gather*}
A_j:=\{y\in Y\colon N(y)>j\}=\{y\in Y\colon \sigma_y(\{t_j(y)\})>0\};\\
 \lambda_j(y):=\one_{A_j}(y)\, \e(t_j(y))\ \text{ for }y\in Y.
\end{gather*}
For $y\in A_j,$ $j\in \Z_+$,  we have  $\lambda_j(y)\in P_y$, that is,
$\sigma_y(\{t_j(y)\})>0$.
 By the Wiener-Wintner Theorem, the limit
$$
\psi_j(x):=\lim_{N\to+\infty}\frac 1N\sum_{n=1}^N f(T^nx)\, \overline{\lambda_j(\pi(x))}^n
$$
 exists in $L^2(\mu)$ (and  for $\mu$-a.e. $x\in X$).
  We remark that if $\nu(A_j)=0$, then the function $\psi_j$ is equal to $0$ $\mu$-a.e.,  and the same holds for  the functions $\theta_j$ and $\phi_j$ defined below. For $\nu$-a.e. $y\in A_j$  we have
\begin{multline}
\label{eq:psij1}
\int f\, \overline{\psi_j}\,d\mu_y=
\lim_{N\to+\infty}\frac 1N\sum_{n=1}^N
\int
f\cdot T^n\overline f\,d\mu_y \cdot
\lambda_j(y)^{n}\\
=\lim_{N\to+\infty}\frac 1N\sum_{n=1}^N
\wh{\sigma_y}(-n)\, \lambda_j(y)^n=\sigma_{y}(\{t_j(y)\}).
\end{multline}

%
 Furthermore, for $j\in \Z_+$ the function $\psi_j$ satisfies
 \begin{equation}\label{E:eigen}
 \psi_j(Tx)=\psi_j(x)\, \lambda_j(\pi(x)).
 \end{equation}
  It follows that $|\psi_j|$ is $T$-invariant. Expressing $|\psi_j|^2$ as
  \begin{equation}\label{E:psitheta}
 |\psi_j|^2=\theta_j\circ\pi, \quad \text{for some } \ \theta_j\in L^\infty(Y,\nu),
 \end{equation}
  a similar computation gives
\begin{multline}
\label{eq:psij2}
\theta_j(y)= \int|\psi_j|^2\,d\mu_y
=\lim_{N\to+\infty}\frac 1{N^2}\sum_{m,n\in[N]}
\int T^nf\cdot T^m\overline f\,d\mu_y \cdot \lambda_j(y)^{m}\, \overline{\lambda_j(y)}^n
\\
=\sigma_{y}(\{t_j(y)\})>0
\text{ for }\nu\text{-a.e. }y\in A_j.
 \end{multline}
 We let
$$
\phi_j(x):=\begin{cases} |\psi_j(x)|\inv\psi_j(x) & \text{for }x\in \pi\inv(A_j);\\
0&\text{for }x\notin \pi\inv(A_j).
\end{cases}
$$
It is immediate to check that $\phi_j$ is an eigenfunction for the eigenvalue $\lambda_j$. By construction, for $i\neq j$ we have $\lambda_i(x)\neq\lambda_j(x)$ for $\mu$-a.e. $x\in X$ except when $\lambda_i(x)=\lambda_j(x)=0$;  thus $\lambda_i(x)\overline{\lambda_j(x)}\neq 1$ $\mu$-a.e..  On the other hand,
 $\E_\mu(\phi_i\overline{\phi_j}\mid\CI(T))= \lambda_i\overline{\lambda_j}\E_\mu(\phi_i\overline{\phi_j}\mid\CI(T))$. Therefore,
 $$
\E_\mu(\phi_i\cdot \overline{\phi_j}\mid\CI(T))=0\ \text{ for all }\, i\neq j, \ i,j\in \Z_+.
$$
Furthermore, for every $j\in \Z_+$ by the definition of the $\phi_j$ we have $\E_\mu(|\phi_j|^2|\CI(T))$ takes the values $0$ or $1$. Hence,  $(\phi_j)_{j\in \Z_+}$ is a relative orthonormal system.

Next we establish identity \eqref{E:f}.
  Arguing as in Part~\eqref{it:relative3} of Proposition~\ref{prop:relative} it suffices to show that
\begin{equation}\label{E:fL2}
\norm f_{L^2(\mu)}^2=\sum_{j\in \Z_+} \bigl\Vert \E_\mu(f\cdot\overline{\phi_j}\mid\CI(T))\bigr\Vert_{L^2(\mu)}^2.
\end{equation}
Using  the definition of the function $\phi_j$ and then \eqref{E:psitheta} combined with  \eqref{eq:psij2},  we get
$$
\sum_{j\in \Z_+} \Bigl|\int f\cdot\overline{\phi_j}\,d\mu_y\Bigr|^2=
\sum_{j\in \Z_+}\one_{A_j}(y)\frac{1}{|\psi_j|^2}\Bigl|\int f\cdot\overline{\psi_j}\,d\mu_y\Bigr|^2
=
\sum_{j\in\Z_+}\one_{A_j}(y)\frac{1}{\sigma_y(\{t_j(y\})}\Bigl|\int f\cdot\overline{\psi_j}\,d\mu_y\Bigr|^2.
$$
By \eqref{eq:psij1}, for $\nu$-a.e. $y\in Y$ the last sum is equal to
$$
\sum_{j\in \Z_+}\one_{A_j}(y)\, \sigma_y(\{t_j(y)\})=\sum_{j\in \Z_+} \sigma_y(\{t_j(y)\})=\sigma_y(\T)=\wh{\sigma_y}(0)= \int |f|^2\,d\mu_y,
$$
where we used the definition of the set $A_j$ to get the first identity, and that for $\nu$-a.e. $y\in Y$  the measure  $\sigma_y$ is atomic and the defining property of the maps $t_j, j\in \Z_+,$ to get the second identity.
 Integrating the established identity over $Y$ with respect to $\nu$ we obtain \eqref{E:fL2}.

  Lastly, the last claim of the lemma follows by the construction of the functions $\psi_j$ and $\phi_j$ for $j\in \Z_+$.
This completes the proof.
\end{proof}
We are now ready to complete the proof of  Theorem~\ref{th:basis-eigen}.
\begin{proof}[End of proof of Theorem~\ref{th:basis-eigen}]
Let $(f_k)_{k\in \N}$ be a dense sequence in $L^2(\mu)$. For every $k\in \N$, we build by induction a countable family $\CF_k$ of functions
in $L^2(\mu)$  such that for every $k\in \N$ we have
 \begin{enumerate}
\item
\label{i:induc1}
 $\CF_1\cup\dots\cup\CF_k$ is a relative orthonormal system of eigenfunctions;
\item
\label{i:induc2}
the function $f_k$ belongs to the closed subspace $\CH_k$ of $L^2(\mu)$ spanned by all functions of the form
$ \phi\, w$ where $\phi\in
\CF_1\cup\dots\cup\CF_k$ and $w\in L^\infty(\mu)$  varies over all $T$-invariant functions.
\end{enumerate}

For $k=1$, the result is given by Lemma~\ref{lem:eigen-ortho} with $f=f_1$. Suppose that the result holds for $k\in\N$, we shall show that it holds for $k+1$. We can decompose $f_{k+1}$ as
 $$
 f_{k+1}=g+f\ \text{ where }g\in\CH_k\text{ and }f\perp \CH_k.
$$
 The space $\CH_k$ is invariant under multiplication by bounded $T$-invariant functions and thus $\E_\mu(f\cdot\overline h\mid\CI(T))=0$ for every $h\in\CH_k$.
 On the other hand,
every $\phi\in\CF_1\cup\dots\cup \CF_k$ is an eigenfunction, and thus, by definition, $T\phi$ belongs to the space $\CH_k$. It follows that $\CH_k$ is invariant under $T$ and thus
\begin{equation}
\label{eq:Fk+1}
\E_\mu(T^nf\cdot\overline h\mid\CI(T))=0\ \text{ for every $n\in \N$ and every } h\in\CH_k.
\end{equation}
Applying Lemma~\ref{lem:eigen-ortho} to the function $f$, we obtain a relative orthonormal system $\CF_{k+1}=(\phi_j)_{ j\in \N}$ of eigenfunctions such that $f$ belongs to the closed linear span of all  functions of the form $\phi_j w$ where $j\in \N$ is arbitrary and $w\in L^\infty(\mu)$ varies over all $T$-invariant functions. Moreover, for every $j\in \N,$ $\phi_j$ belongs to the smallest $T$-invariant subspace of $L^2(\mu)$ containing all functions of the form $fw$ where $w$ varies over all $T$-invariant functions  in $L^\infty(\mu)$. Hence, by~\eqref{eq:Fk+1}, for every $j\in \N$ and every $h\in\CH_k$ we have $\E_\mu(\phi_j \cdot \overline h\mid\CI(T))=0$; this holds, in particular, for all functions $h$ belonging to $\CF_1\cup\dots\cup\CF_k$. Therefore, $\CF_1\cup\dots\cup \CF_k\cup\CF_{k+1}$ is a relative orthonormal system.

The closed subspace of $L^2(\mu)$ spanned by functions of the form $\phi\, w$, where $\phi\in \CF_1\cup\dots\cup \CF_k\cup\CF_{k+1}$ and $w\in L^\infty(\mu)$ varies over all  $T$-invariant functions,  contains $f$ and $g$, and thus $f_{k+1}$. This completes the induction.

 We choose an enumeration of the countable set  $\CF:=\bigcup_{k=1}^\infty\CF_k$  and write it as $(\phi_j')_{j\in\N}$; then $(\phi_j')_{j\in\N}$  is a relative orthonormal system. The closed subspace of $L^2(\mu)$ spanned by
all functions of the form $\phi_j'\, w$, where $j\in\N$ and $w\in L^\infty(\mu)$ varies over all  $T$-invariant functions, contains  the functions $f_k$  for every $k\in \N$, and thus is equal to $L^2(\mu)$. Hence, $(\phi_j')_{j\in\N}$  is a relative orthonormal basis of eigenfunctions, completing the proof.
\end{proof}

\section{Decomposition of correlation sequences}
\label{sec:correl-seq}
In this section we prove the decomposition results stated in Sections~\ref{SS:mcs} and \ref{SS:2.5}.
\subsection{Anti-uniformity in norm} We start with some preparatory results.
The main tool in  verifying anti-uniformity properties is  the following inner product  space variant of  a well known estimate of van der Corput
(for a proof see \cite{B87}):
\begin{vdcLemma}
Let $d\in\N$, $\CH$ be an inner product space,  $\xi\colon \N^d\to \CH$ be a bounded sequence, and
let  $\bI$
 be a F\o lner sequence in $\N^d$. Then
$$
\Limsup\bigl\Vert\Av_\bI\, \xi_\bn\bigr\Vert^2
\leq
4\limsup_{H\to+\infty}\frac 1{H^d}
\sum_{\bh\in[H]^d}
\Limsup\bigl\vert\Av_{\bn,\bI}\,\langle \xi_{\bn+\bh}, \xi_\bn\rangle\bigr\vert.
$$
\end{vdcLemma}

\begin{proposition}
\label{prop:unif-in-norm}
Let $d,\ell,s, t\in\N$.  Then there exists a positive integer  $k=k(d,\ell,s,t)$ such that for every system $(X,\mu,T_1,\dots,T_\ell)$, functions $f_1,\dots, f_s\in L^\infty(\mu)$ bounded by $1$,
polynomial mappings $\vec{p_1},\ldots, \vec{p_s}\colon\N^d\to\Z^\ell$ of degree at most $t$,  and sequence  $w\in \ell^\infty(\N^d)$, we have
$$
\Limsup\bigl\Vert\Av\,
w(\bn) \, T_{\vec{p_1}(\bn)}f_1\cdot\ldots\cdot
T_{\vec{p_s}(\bn)}f_s\bigr\Vert_{L^2(\mu)}
\leq  4\norm w_{U^{k+1}(\N^d)}.
$$

Furthermore, if the polynomial mappings  are linear, then we can take $k=s$.
\end{proposition}
\begin{proof}[Sketch of the proof] Let $\bI=(I_i)_{j\in\N}$ be a F\o lner in $\N^d$.
After passing to a subsequence,  we can assume that the sequence $w$ admits correlations along $\bI$.
 It suffices  to show that
$$
\Limsup\bigl\Vert\Av_\bI\,w(\bn) \cdot T_{\vec{p_1}(\bn)}f_1\cdot\ldots\cdot
T_{\vec{p_s}(\bn)}f_s\bigr\Vert_{L^2(\mu)}
\leq  4\norm w_{\bI,k+1}
$$
for some $k=k(d,\ell,s,t)$.
To verify this, one  applies an inductive argument, often called PET induction, introduced by V. Bergelson in~\cite{B87}.
Each step uses the van der Corput Lemma in $L^2(\mu)$, invariance of the measure under some of the transformations, and the Cauchy-Schwarz inequality.
The details are similar to several other arguments in the literature
(see for example the proof of~\cite[Lemma~3.5]{FHK13})
and so we do not give the proof.

In the case of linear polynomials, it can be shown
by induction on $s\in \Z_+$ that one can take $k=s$; for $s=0$ the statement is trivial
 and
the inductive step can be carried out  as the inductive step  in the proof of Theorem~\ref{th:correl-lin} below.
\end{proof}

\subsection{Proof of Theorem~\ref{th:correl-poly}}
\label{subsec:Proof-Th5}
By Theorem~\ref{th:regular+antiunif} it suffices to show that the sequence $a\colon \N^d\to \C$ given by
$$
a(\bn):= \int f_0\cdot T_{\vec{p_1}(\bn)}f_1\cdot\ldots\cdot
T_{\vec{p_s}(\bn)}f_s\,d\mu, \quad \bn\in \N^d,
$$
is $k$-regular and $(k+1)$-anti-uniform for some $k\in \N$ that depends only on the integers $d,\ell,s,t$. The regularity is given by Proposition~\ref{prop:weight-nil} and the  anti-uniformity by Proposition~\ref{prop:unif-in-norm}.
\qed

\subsection{Proof of Theorem~\ref{th:correl-lin}}
\label{subsec:Proof-Th6b}
We start with some preparatory results.
\begin{lemma}
\label{lem:estimate}
 Let $d, \ell \in \N$ and $L_1,\ldots, L_\ell\colon \N^d\to \Z$ be linearly independent linear forms. Then there exists
a constant  $C:=C(d, L_1,\ldots, L_\ell)$ such that the following holds:
If   $(X, \mu, T_1,\ldots,  T_\ell)$ is a system
and $f_0,\ldots, f_\ell\in L^\infty(\mu)$ are functions bounded by $1$,
 then
\begin{equation}
\label{eq:estimate}
\limav \Big|\int f_0\cdot T_1^{L_1(\bn)}f_1\cdot\ldots \cdot T_\ell^{L_\ell(\bn)} f_\ell\, d\mu\Big|^2
\leq C\,\min_{1\leq i\leq\ell} \nnorm{f_i}^2_{T_i,\mu,2}.\end{equation}
\end{lemma}
\begin{proof}
We first note that the limit on the left hand side of~\eqref{eq:estimate} can be rewritten as
$$
\limav\Bigl(\int (f_0\otimes\overline{f_0})\cdot \prod_{i=1}^\ell
 (T_i\times T_i)^{L_i(\bn)}(f_i\otimes\overline{f_i})\, d(\mu\times\mu)\Bigr)
$$
and thus exists  by Theorem~\ref{T:Walsh}. Therefore, in~\eqref{eq:estimate} we can restrict to  averages   taken on  the cubes $[N]^d$, $N\in \N$, that is, it suffices to obtain bounds for the following limit
\begin{equation}\label{E:csr}
\lim_{N\to+\infty}\frac{1}{N^d} \sum_{\bn\in [N]^d}\Big|\int f_0\cdot T_1^{L_1(\bn)}f_1\cdot\ldots \cdot T_\ell^{L_\ell(\bn)} f_\ell\, d\mu\Big|^2.
\end{equation}

 Next, we  claim  that it suffices to consider the case where $d=\ell$. Indeed, since the linear forms are linearly independent, we have $d\geq\ell$, and
if $d>\ell$, then there exist linear forms $L_{\ell+1},\ldots, L_d \colon \N^d\to \Z$ such that the  linear forms $L_1,\ldots, L_d$ are linearly independent.
Then applying the $d=\ell$ case of the result for this set of linear forms  and the functions $f'_0,\ldots, f'_d$ defined by $f'_i:=f_i$ for $i=0,\ldots, \ell$ and $f'_i:=1$ for $i=\ell+1,\ldots, d,$  we get the asserted estimate.
Henceforth, we assume that $d=\ell$.

Let $\vec L\colon\N^d\to\N^d$ be defined  by  $\vec L(\bn):=(L_1(\bn),\ldots, L_d(\bn)),$ $\bn\in \N^d$.
Since the linear forms $L_1,\ldots, L_d$ are linearly independent,  the linear map $\vec L$ is injective. Furthermore, there exists  a positive integer $M=M(L_1,\ldots, L_d)$ such that $\vec L([N]^d)\subset [-MN,MN]^d$ for every $N\in \N$.
This easily implies that the  limit in \eqref{E:csr}  is bounded by
\begin{multline*}
(3M)^d\lim_{N\to +\infty} \frac{1}{(2N+1)^d}\sum_{\bn \in [-N,N]^d}  \Big|\int f_0\cdot T_1^{n_1}f_1\cdot\ldots \cdot T_d^{n_d} f_d\, d\mu\Big|^2\\
= (3M)^d\lim_{N\to +\infty} \frac{1}{(2N+1)^d}\sum_{\bn \in [-N,N]^d}
\int (f_0\otimes\overline{f_0})\cdot \prod_{i=1}^d
 (T_i\times T_i)^{n_i}(f_i\otimes\overline{f_i})\, d(\mu\times\mu).
\end{multline*}
By the ergodic theorem, the last limit is equal to
$$
\int (f_0\otimes\overline{f_0})\cdot \prod_{i=1}^d
\E_\mu\bigl(f_i\otimes\overline{f_i}\mid(\CI(T_i\times T_i)\bigr)\,d(\mu\times\mu).
$$
For $i=1,\dots,d$ this quantity is bounded by
 $$
 \bigl\Vert \E_\mu\bigl(f_i\otimes\overline{f_i}\mid(\CI(T_i\times T_i)\bigr)\bigr\Vert_{L^2(\mu)}=\nnorm{f_i\otimes\overline{f_i}}_{T_i\times T_i,\mu\times\mu,1}\leq \nnorm{f_i}_{T_i,\mu,2}^2
$$
 where we used that all the functions are bounded by $1$ and the estimate \eqref{eq:seminorm4}. This completes the proof.
\end{proof}

\begin{proposition}
\label{prop:independent}
Let   $d, \ell\in \N$ and $L_1,\ldots, L_\ell\colon \N^d\to \Z$ be linearly independent linear forms. Then for every $\ve>0$ there exists a constant
 $C:=C(d, L_1,\ldots, L_\ell, \varepsilon)>0$ such that the following holds:
If
 $(X,\mu, T_1, \ldots, T_\ell)$ is a system and $f_0,\ldots, f_\ell\in L^\infty(\mu)$ are functions bounded by $1$,  then for every  sequence $w\in \ell^\infty(\N^d)$   we have
\begin{equation}\label{eq:independent1}
\Limsup\bigl\vert\Av\, \big(w(\bn)\cdot
\int f_0\cdot T_1^{L_1(\bn)}f_1\cdot\ldots \cdot T_\ell^{L_\ell(\bn)} f_\ell\, d\mu\big)\bigr\vert
\leq C\norm w_{U^2(\N^d)} +\varepsilon\norm w_\infty.
\end{equation}
\end{proposition}
\begin{remark}
 This proves that the correlation sequence defined by the integral is $2$-anti-uniform with constants  $C$ that do not depend on the functions
$f_0,\ldots, f_\ell$ as long as they are bounded by $1$. This condition is essential in  the proof of Theorem~\ref{th:correl-lin} and Theorem~\ref{th:basis-eigen} is key in establishing the condition.
\end{remark}
 \begin{proof}
 Let $\varepsilon\in (0,1)$, $w\in \ell^\infty(\N^d)$,  and $\bI=(I_j)_{j\in \N}$ be a F\o lner sequence in $\N^d$. By passing to a  F\o lner  subsequence  we can assume that $w$ admits correlations along $\bI$.

By the defining property of the factor $\CZ_1$ (see Section~\ref{subse:Zk}) we have that
$\nnorm {f_i-\E_\mu(f\mid\CZ_1(X,\mu,T_i))}_{\mu,T_i,2}=0$
for $i=1,\dots,\ell$.
 Hence, by Lemma~\ref{lem:estimate}, the $\limsup$ in~\eqref{eq:independent1} remains unchanged if we replace each function $f_i$ with $\E_\mu(f_i\mid\CZ_1(X,\mu,T_i))$. Therefore,
 we can and will assume that  for $i=1,\ldots, \ell$, the function $f_i$ is measurable with respect to $\CZ_1(X,\mu,T_i)$.
For $i=1, \ldots, \ell$ let $\mu=\int \mu_{i,x}\, d\mu(x)$ be the ergodic decomposition of the system
$(X,\mu,T_i)$.

By Theorem~\ref{th:basis-eigen}, for $i=1,\ldots, \ell$ the space $L^2(\CZ_1(X,\mu,T_i),\mu)$ admits a relative orthonormal basis $(\phi_{i,j})_{j\in \N}$ such that $\phi_{i,j}$ is an eigenfunction of $(X,\mu,T_i)$ with eigenvalue $\lambda_{i,j}$ for  $j\in \N$.  We write $(f_{i,j})_{j\in \N}$ for the coordinates of $f_i$ in this base. We recall that $f_{i,j}$ is invariant under $T_i$ and that \
 $$
 f_{i,j}=\E_\mu(f_i\, \overline{\phi_{i,j}}\mid\CI(T_i)). 
 $$
 Then $\norm{f_{i,j}}_{L^\infty(\mu)}\leq 1$ for all $i,j\in \N$. Moreover,
  by Part~\eqref{it:relative2} of Proposition~\ref{prop:relative} we have
\begin{gather}
\label{eq:fi-fij1}
\E_\mu(|f_i|^2\mid \CI(T_i))=\sum_{j\in \N}|f_{i,j}|^2, \quad \mu\text{-a.e.};\\
\label{eq:fi-fij3}
f_i=\sum_{j\in \N} f_{i,j}\, \phi_{i,j},
\end{gather}
where convergence in \eqref{eq:fi-fij1} is pointwise and in \eqref{eq:fi-fij3} is in $L^2(\mu)$ and in $L^2(\mu_{i,x})$ $\mu$-a.e..

For $i=1,\dots,\ell$ we separate the series~\eqref{eq:fi-fij3} in two parts. For $x\in X$ we let
\begin{gather*}
E_i(x):=\{ j\in \N \colon |f_{i,j}(x)|^2\geq \ve^{10^i}\};\\
g_i(x):=\sum_{j\in E_i(x)} f_{i,j}(x)\, \phi_{i,j}(x)\ \text{ and }\ h_i:=f_i-g_i.
\end{gather*}
By~\eqref{eq:fi-fij1}  and since all functions are bounded by $1$, we have
  \begin{equation}\label{E:two}
  |E_i(x)|\leq \ve^{-10^i} \quad  \mu\text{-a.e.},  \quad
  \norm{g_i}_{L^\infty(\mu)}\leq \ve^{-10^i}, \quad i=1,\ldots, \ell.
  \end{equation}
Furthermore, since $f_{i,j}$ are $T_i$-invariant we have
$E_i(T_ix)=E_i(x) $ $\mu$-a.e.  and the set $A_{i,j}=\{x\in X\colon j\notin E_i(x)\}$ is invariant under $T_i$. We have
$$h_i=\sum_{j\in\N}\one_{A_{i,j}}\, f_{i,j}\, \phi_{i,j}$$
and thus the coordinates of the function $h_i$ in the base
$(\phi_{i,j})_{j\in \N}$ are the functions $\one_{A_{i,j}}f_{i,j}$.
By
Proposition~\ref{prop:basis-eigen2}, we obtain
\begin{equation}\label{E:semihi}
 \nnorm{h_i}_{T_i,\mu_{i,x},2}^4=\sum_{j\in\N}|\one_{A_{i,j}}(x)f_{i,j}(x)|^4
 = \sum_{j\notin E_i(x)} |f_{i,j}(x)|^4\quad \mu\text{-a.e.}, \quad i=1,\ldots, \ell.
 \end{equation}

 We therefore  have
$$
\nnorm{h_i}_{T_i,\mu,2}^4= \int \nnorm{h_i}_{T_i,\mu_{i,x},2}^4\, d\mu= \int \sum_{j\notin E_i(x)}|f_{i,j}(x)|^4\, d\mu\leq \ve^{ 10^i}\int \sum_{j\in \N}|f_{i,j}(x)|^2\, d\mu
\leq  \ve^{ 10^i},
$$
where we used \eqref{E:seminonerg} in the appendix to get the first identity,  \eqref{E:semihi} to get the second identity,  the definition of the sets $E_i(x)$ to get the first estimate, and \eqref{eq:fi-fij1} combined with the fact that the functions $f_i$ are bounded by $1$ to get the last estimate.

 Let $C$ be the constant defined in Lemma~\ref{lem:estimate}. Combining this lemma with the preceding estimates we deduce  for  $m=1,\dots,\ell$ that
\begin{multline*}
\Limsup\;\Av_\bI\Bigl|
\int f_0\cdot\Bigl(\prod_{i=1}^{m-1}T_i^{L_i(\bn)}g_i\Bigr)\cdot T_m^{L_m(\bn)}h_m\cdot\Bigl(\prod_{i=m+1}^\ell T_i^{L_i(\bn)}f_i\Bigr)\,d\mu\Bigr|^2\\
\leq C\,\prod_{i=1}^{m-1}\norm{g_i}_{L^\infty(\mu)}^2\cdot \nnorm{h_m}_{T_m,\mu,2}^2\leq C \ve^{-2\sum_{i=1}^{m-1}10^i} \, \ve^{5\cdot 10^{m-1}}
\leq C\ve^2.
\end{multline*}
Using the Cauchy-Schwarz inequality and telescoping, we obtain
$$
\Limsup\Bigl\vert\Av_\bI\,
\bigl(w(\bn)\int f_0\cdot \prod_{i=1}^\ell T_i^{L_i(\bn)}f_i\,d\mu-
w(\bn)\int  f_0\cdot \prod_{i=1}^\ell T_i^{L_i(\bn)}g_i\,d\mu\bigr)\Bigr|
\leq \ell C^{1/2}\ve\norm w_\infty.
$$

On the other hand, using the definition of the functions $g_1,\ldots, g_\ell$   and recalling that  for $i=1,\ldots, \ell$ and $j\in \N$ the function $\phi_{i,j}$ is a
$T_i$-eigenfunction with eigenvalue $\lambda_{i,j}$, we get
\begin{multline*}
\int  f_0\cdot \prod_{i=1}^\ell T_i^{L_i(\bn)}g_i\,d\mu=
\int \sum_{j_1\in E_1,\ldots,j_\ell\in E_\ell} f_0\cdot \prod_{i=1}^\ell f_{i,j_i}\cdot T_i^{L_i(\bn)} \phi_{i,j_i}\,d\mu
\\
=\int \sum_{j_1\in E_1,\ldots,j_\ell\in E_\ell} g_{j_1,\ldots,j_\ell}
\cdot \prod_{i=1}^\ell \lambda_{i,j_i}^{L_i(\bn)}\, d\mu
\end{multline*}
where
$$
g_{j_1,\dots,j_\ell}:=f_0\cdot\prod_{i=1}^\ell  f_{i,j_i}\cdot\phi_{i,j_i}.
$$
Since for $\mu$-a.e. $x\in X$ we have by \eqref{E:two} that  $|E_i(x)|\leq \ve^{-10^i}$  for $i=1,\ldots, \ell$, we deduce that
for $\mu$-a.e. $x\in X$ the sum contains at most $\varepsilon^{-10^{\ell+1}}$ terms. Moreover, since $L_1,\ldots, L_\ell$ are linear, using the van der Corput Lemma on $\C$ and the fact that  the functions
$g_{j_1,\dots,j_\ell}$ are bounded by $1$, we have the pointwise estimate
\begin{multline*}
\Limsup\Bigl\vert \Av_\bI\, w(\bn)
\cdot \prod_{i=1}^\ell \lambda_{i,j_i}^{L_i(\bn)}\Bigr|^2
\leq 4\limsup_{H\to+\infty}\frac 1{H^d}\sum_{\bh\in [H|^d}
\bigl\vert\Limav{\bn,\bI} \,
w(\bn+\bh)\overline{w(\bn)}\bigr\vert
\\
 \leq 4\lim_{H\to+\infty}\frac 1{H^d}\sum_{\bh\in [H|^d}
\norm{\sigma_\bh \, w\cdot\overline w}_{\bI,1}\leq 4\norm w_{\bI,2}^2\leq 4\norm w_{U^2(\N^d)}^2
\end{multline*}
where we used~\eqref{eq:unif-seminorm1} and~\eqref{eq:unif-seminorms2}.
Combining the above estimates, we deduce that  the left hand side of~\eqref{eq:independent1} is bounded by
$2\ve^{-10^{\ell+1}}\norm w_{U^2(\N^d)}+\ell C^{1/2}\ve\norm w_\infty$. This completes the proof.
\end{proof}


\begin{proof}[End of proof of Theorem~\ref{th:correl-lin}] We can  assume that the functions $f_0,\ldots, f_\ell$ are bounded by $1$.
Furthermore, we can extract a linearly independent subset of $r$ elements of $\{L_1,\dots,L_\ell\}$;  hence, after  reordering the linear forms we can assume that the first $r$ ones  are linearly independent.

Let $k:=\ell-r+1$. In order to prove that the sequence
$a\colon \N^d\to \C$  defined  by
$$
a(\bn):=\int f_0\cdot T_1^{L_1(\bn)}f_1\cdot\ldots\cdot T_\ell^{L_\ell(\bn)}f_\ell\,d\mu, \quad \bn \in \N^d,
$$
admits a decomposition of the announced type, it suffices
by Theorem~\ref{th:regular+antiunif} to show that it is $k$-regular and $(k+1)$-anti-uniform.

The regularity follows from  Proposition~\ref{prop:weight-nil}.

 Next we verify $(k+1)$-anti-uniformity.  Let $C=C(d, L_1,\ldots, L_r,\varepsilon)$ be the constant defined by Proposition~\ref{prop:independent}. We can assume that $C>1$.  For fixed $d\in \N$ we will prove by induction  on $k\in \N$
that if the functions $f_0, \ldots, f_\ell$ are bounded by $1$, then
for every $\varepsilon>0$ and every  $w\in \ell^\infty(\N^d)$ with $\norm{w}_\infty=1$
we have
$$
\Limsup\bigl\vert\Av\, a(\bn)\, w(\bn)\bigr\vert\leq 4 C\norm b_{U^{k+1}(\N^d)}+4\ve^{\frac{1}{2^{k-1}}}.
$$
 It will then follow that the sequence $a$ is $(k+1)$-anti-uniform with anti-uniformity constant $C':=4C(d, L_1,\ldots, L_r,\frac{1}{4}\varepsilon^{2^{k-1}})$.

 For $k=1$ the statement follows from Proposition~\ref{prop:independent}.

 Suppose that $k\geq 2$ and that the statement holds for $k-1$. Let $\varepsilon>0$ and $w\in \ell^\infty(\N^d)$.
Let also $\bI=(I_j)_{j\in \N}$ be a F\o lner sequence in $\N^d$. By passing to a  F\o lner  subsequence  we can assume that $w$ admits correlations along $\bI$.

Composing with $T_\ell^{-L_\ell(\bn)}$ we rewrite the sequence $a\colon \N^d\to \C$ as\footnote{This maneuver is necessary;  a direct application of the van der Corput Lemma  produces a weaker   estimate involving a seminorm of order $k+2$.}
$$
a(\bn)=\int T_\ell^{-L_\ell(\bn)}f_0\cdot
T_\ell^{-L_\ell(\bn)}T_1^{L_1(\bn)} f_1\cdot\ldots
\cdot T_\ell^{-L_\ell(\bn)}T_{\ell-1}^{L_{\ell-1}(\bn)} f_{\ell-1}\cdot f_\ell\ d\mu, \quad \bn\in\N^d.
$$
Using the Cauchy-Schwarz inequality and that $\norm{f_\ell}_{L^\infty(\mu)}\leq 1$, we get that
\begin{equation}\label{E:cs}
\Limsup\bigl\vert \Av_\bI \,a(\bn)\, w(\bn)\big|^2\leq \Limsup\bigl\Vert\Av_\bI\,\xi_\bn\bigr\Vert_{L^2(\mu)}^2
\end{equation}
where $\xi\colon \N^d\to L^2(\mu)$ is defined by
$$
\xi_\bn:=w(\bn)\cdot  T_\ell^{-L_\ell(\bn)}f_0\cdot
T_\ell^{-L_\ell(\bn)}T_1^{L_1(\bn)} f_1\cdot\ldots
\cdot T_\ell^{-L_\ell(\bn)}T_{\ell-1}^{-L_{\ell-1}(\bn)} f_{\ell-1}, \quad \bn\in \N^d.
$$
Using  van der Corput's Lemma in $L^2(\mu)$ for the sequence $\xi$,  we obtain that the right hand side of \eqref{E:cs}
  is bounded by
\begin{equation}
\label{eq:correl-lin3}
4\limsup_{H\to+\infty}\frac 1{H^d}
\sum_{\bh\in[H]^d}
\Limsup\bigl\vert \Av_{\bn,\bI}\,
\langle\xi_{\bn+\bh}, \xi_\bn\rangle\bigr|.
\end{equation}
Recall that $\sigma_\bh w(\bn)= w(\bn+\bh)$ for all $\bh, \bn\in\N^d$. A simple computation gives that for every $\bh\in\N^d$ we have
\begin{equation}\label{eq:correl-lin4}
\frac 1{|I_j|}
\sum_{\bn\in I_j}\langle\xi_{\bn+\bh}, \xi_\bn\rangle\
=\frac 1{|I_j|}\sum_{\bn\in I_j} \sigma_\bh w(\bn)\cdot\overline{w(\bn)}\int
\wt f_{0,\bh}\cdot T_1^{L_1(\bn)}\wt f_{1,\bh}\cdot\ldots\cdot T_{\ell-1}^{L_{\ell-1}(\bn)}\wt f_{\ell-1}\, d\mu
\end{equation}
where  $\wt f_{j,\bh}:=T_\ell^{-L_\ell(\bh)}T_j^{L_j(\bh)} f_j\cdot\overline{f_j}$ for $j=0,\dots,\ell-1$ and $L_0:=0$.

Note for every $\bh\in\N^d$ the sequence $(w(\bn+\bh)\, \overline{w(\bn)})_{\bn\in\N^d}$ admits correlations along $\bI$ and that $ \norm{\wt f_{j,\bh}}_{L^\infty(\mu)}\leq 1$ for $j=0,\dots,\ell-1$.
   The expression on the right hand side of~\eqref{eq:correl-lin4} is thus of the type studied, with $(\ell-1)$ in place of $\ell$. Hence, the  induction hypothesis applies and gives  that  for every $\bh\in\N^d$ we have
$$
\Limsup\bigl\vert\Av_{\bn,\bI}\,
\langle\xi_{\bn+\bh}, \xi_\bn\rangle\bigr\vert
\leq  4\,  C\, \norm{\sigma_\bh w\cdot w}_{\bI,k}+4\, \ve^{\frac{1}{2^{k-2}}}.
$$
It is important for the last part of the  argument that the constant $C$  is independent of  the parameter $\bh$.
Combining the above, we deduce that  the expression~\eqref{eq:correl-lin3} is bounded by
\begin{multline*}
16\, C\, \limsup_{H\to+\infty} \frac 1{H^d}\sum_{\bh\in [H]^d}
\norm{\sigma_\bh w\cdot\overline w}_{\bI,k}+16\, \ve^{\frac{1}{2^{k-2}}}\leq
\\
16\, C\, \limsup_{H\to+\infty}\Bigl( \frac 1{H^d}\sum_{\bh\in [H]^d}\norm{\sigma_\bh w\cdot\overline w}_{\bI,k}^{2^{k}}\Bigr)^{\frac{1}{2^{k}}}+ 16\, \ve^{\frac{1}{2^{k-2}}}
=16\,  C\, \norm w_{\bI,k+1}^2 + 16\, \ve^{\frac{1}{2^{k-2}}}
\end{multline*}
where the last identity follows from the inductive property~\eqref{eq:unif-seminorms2} of the uniformity seminorms. Putting together  the previous estimates and taking  square roots we get the announced bound.  This completes the induction and the  proof.
\end{proof}


\subsection{Proof of Theorem~\ref{th:cor-infty}}
We use the   variant of Furstenberg's correspondence principle already used in Section~\ref{subsec:def-seminorms} for a single sequence. In the case of several sequences it is proved in a
similar fashion and gives the following:
\begin{proposition}\label{P:correspondence}
 Let
$\ell, s\in \N$ and $a_1,\ldots, a_s\colon \Z^\ell\to \C$ be bounded sequences such that the family $\CF=\{a_1,\dots,a_s\}$ admits
correlations along the F\o lner sequence $\bf I$ in $\N^\ell$.
Then there exists a  topological dynamical system  $(X,T_1,\ldots, T_\ell)$, where $T_1,\ldots, T_\ell$ are commuting homeomorphisms,
functions $f_1,\ldots,f_s\in C(X)$, and a Borel probability measure $\mu$ on $X$ that is $T_i$-invariant for $i=1,\ldots, \ell$, such that
 $$
 \int    T_{\vec{n_1}}
f_1\cdot \ldots \cdot T_{\vec{n_s}}
f_s\, d\mu=
\Limav{\bk,\bI}\bigl(\prod_{i=1}^s a_i(\bk+\vec{n_i})\bigr)
$$
for every $\vec{n_1}, \ldots, \vec{n_s}\in \Z^\ell$.
\end{proposition}
Combining  this with Theorem~\ref{th:correl-poly} we immediately deduce
Theorem~\ref{th:cor-infty}. \qed

\subsection{Proof of Theorem~\ref{th:cor-finite}}
We start by recalling the definition of the Gowers norms in $\Z_N^d$.
\begin{definition}
Let $d,N\in \N$  and $f\colon \Z_N^d\to \C$ be a function. For every $\bh\in\Z_N^d$ we write $f_\bh(\bn):=f(\bn+\bh),$ $\bn\in \Z_N^d$. For $s\in \N$
we denote by $\norm f_{U^s(\Z_N^d)}$ the \emph{Gowers $U^s(\Z_N^d)$-norm} of $f$ that is defined inductively as follows:    We let
$$
\norm f_{U^1(\Z_N^d)}:=|\E_{\bn\in\Z_N^d}f(\bn)|
=\Bigl(\E_{\bn,\bh\in\Z_N^d} f(\bn)\, \overline f(\bn+\bh)\Bigr)^{1/2},
$$
and for every $s\geq 1$ we let
$$
\norm f_{U^{s+1}(\Z_N^d)}:=\Bigl(\E_{\bh\in\Z_N^d}\norm{f\cdot \overline f_\bh}_{U^s(\Z_N^d)}^{2^s}\Bigr)^{1/2^{s+1}}.
$$
\end{definition}
For $d=1$ the next result is deduced in \cite[Section~3.1]{HK12} from the inverse theorem for the Gowers norms in $\Z_N$ \cite{GTZ12}.
A multidimensional extension of this inverse theorem was established in
 \cite{ CaSz12, Sz12} (for an alternate proof see \cite{T15}) and the argument in \cite[Section~3.1]{HK12}
 allows us to deduce the following result:
 \begin{proposition}[\mbox{\cite[Section~3.1]{HK12}} for $d=1$]
 Let $d,k\in \N$, $C>0$, and $\ve>0$. Then there exists a $k$-step nilmanifold $X$ such that the following holds: If  $N\in \N$ and  $a\colon \Z_N^d\to \C$ satisfies
 \begin{equation}\label{E:estfin}
\Big| \E_{\bn\in \Z_N^d}a(\bn)\, b(\bn)\Big|\leq C \norm{b}_{U^{k+1}(\Z_N^d)} \quad \text{for every } \ b\colon \Z_N\to \C,
 \end{equation}
then we have the decomposition $a=a_\st+a_\er$ where
\begin{enumerate}
\item
$a_\st$ is a convex combination of $k$-step nilsequences   defined by functions on $X$ with  Lipschitz norm at most $1$;
\item
$
\E_{\bn\in \Z_N^d}|a_\er(\bn)|\leq\ve$.
\end{enumerate}
 \end{proposition}
It follows from  the previous result  that in order to complete the proof of
 Theorem~\ref{th:cor-finite} it suffices to show that the sequence $b\colon \Z_N^d\to \C$
 defined by
\eqref{E:corfin} satisfies the property \eqref{E:estfin} for some $k\in \N$ and $C>0$ that depend only on the integers $d,\ell, s,t$.
This can be verified directly with  $C=1$ by using a PET-induction argument  (as in \cite[Lemma~3.5]{FHK13})) and the estimate
$
|\E_{\bn\in \Z_N^d} a(\bn)|^2\leq \E_{\bh\in \Z_N^d}|\E_{\bn\in \Z_N^d} a(\bn+h)\, \overline {a(\bn)}|,
$
but it turns out to be simpler  to deduce  \eqref{E:estfin} directly from   Proposition~\ref{prop:unif-in-norm} as follows:

For fixed $N\in \N$ we interpret $b\colon \Z_N^d\to \C$ as a periodic sequence in  $\N^d$, and note that
$b$ can be represented as
$$
b(\bn)=\int  T_{\vec{p_1}(\bn)}f_1\cdot\ldots\cdot T_{\vec {p_s}(\bn)}f_s\, d\mu, \quad \bn\in \N^d,
$$
where $X:=\Z_N^\ell$, $\mu$ is the Haar measure on $X$, and for $i=1,\ldots, s$,   $f_i:=a_i$ and $T_i$ is the measure preserving transformation on $\Z_N^\ell$ (with the Haar measure) that shifts the $i$-th coordinate of an element of $\Z_N^\ell$ by $1$ and leaves the other coordinates unchanged.
Notice also that if the sequence $c\colon \N^d\to \N$ is  $N$-periodic on every coordinate direction,
then
$$
\limav c(\bn)=\E_{\bn\in\Z_N^d}c(\bn) \ \text{ and }\  \norm{c}_{U^{k}(\Z^d)}=\norm{c}_{U^{k}(\Z_N^d)}
$$ for every $k\in \N$.
Keeping all these in mind and using    Proposition~\ref{prop:unif-in-norm}, we deduce    that  the estimate  \eqref{E:estfin} holds with $C=4$ for some $k=k(d,\ell,s,t)$ and  we can take $k=s$ if  the polynomial mappings are linear. This completes the proof of Theorem~\ref{th:cor-finite}. \qed


 \subsection{Proof of Theorem~\ref{T:4}}
 \label{subsec:proofT4}
First we check that  $\CN_d$ is a linear subspace of $\ell^\infty(\N^d)$. For  $i=1,2$, let $(\psi_i(\bn))_{\bn\in\N^d}$ be a $k_i$-step nilsequence given by $\psi_i(\bn)=\Psi_i(\prod_{j=1}^d \tau_{i,j}^{n_i}\cdot e_{X_i})$  where $X_i=G_i/\Gamma_i$ is a $k_i$-step nilmanifold and $\tau_i\in G_i$.
 Then for $k:=\max(k_1,k_2)$ their sum is  the $k$-step nilsequence $(\Psi(\prod_{i=1}^d \tau_i^{n_i}\cdot e_X))_{n_1,\ldots, n_d\in \N}$ where  $X:=X_1\times X_2=(G_1\times G_2)/(\Gamma_1\times\Gamma_2)$,
 $\Psi(x_1,x_2):=\Psi_1(x_1)+\Psi_2(x_2)$, $e_X:=(e_{X_1},e_{X_2})$, and
$\tau_j:=(\tau_{1,j},\tau_{2,j})$ for $j=1,\ldots, d$.

 Next we  show that the space
 $\mathcal{MC}_\text{d,pol}$ is linear, a similar argument works for the space $\mathcal{MC}_\text{d,lin}$.
  Let  $a\colon \N^d\to \C$ be given by
   $$
a(\bn)=   \int  f_0\cdot  T_{\vec{p_1}(\bn)}
f_1\cdot \ldots \cdot T_{\vec{p_s}(\bn)}
f_s\, d\mu, \quad \bn\in \N^d,
$$
  for some system $(X, \mu, T_1,\ldots, T_\ell)$, functions $f_0,\ldots, f_s\in L^\infty(\mu)$, and polynomial mappings $\vec{ p_1},\ldots, \vec{ p_s}\colon \N^d\to \Z^\ell$.
  Let also $a'\colon \N^d\to \C$ be defined by a similar formula, with $\ell'$ in place of $\ell$, $(X', \mu', T'_1,\ldots, T'_{\ell'})$ in place of $(X, \mu, T_1,\ldots, T_\ell)$, $s'$ in place of $s$, $\vec {p'_1},\dots,\vec{p'_{s'}}\colon\Z^d\to\Z^{\ell'}$ in place of $\vec{p_1},\dots,\vec{p_s}$, and
$f'_0, \ldots, f'_{s'}\in L^\infty(\mu')$ in place of $f_0,\dots,f_s$. We define a system
$(Y,\nu,S_1,\dots S_{\ell+\ell'})$ by letting $Y$ to be the disjoint union $X\uplus X'$ of $X$ and $X'$, $\nu=\frac 12(\mu+\mu')$, and defining the transformations $S_1,\dots,S_{\ell+\ell'}$ of $Y$  by
\begin{align*}
S_j|_X:=T_j &\ \text{ and }\ S_j|_{X'}:=\id, & \text{for } &1\leq j\leq \ell ;\\
S_j|_X:=\id &\ \text{ and }\ S_j|_{X'}:=T'_j, &  \text{for }& \ell<j\leq \ell+\ell'.
\end{align*}
We also define the polynomial mappings  $\vec{q_1},\ldots, \vec{q_s}\colon \N^d\to \Z^{\ell+\ell'}$ and
$\vec{ q'_1},\ldots, \vec{ q'_{s'}}\colon \N^d\to \Z^{\ell+\ell'}$ as follows:  
\begin{align*}
\vec{q_i}(\bn)&:=(p_{i,1}(\bn),\dots,p_{i,\ell}(\bn),0,\dots,0), & \text{for } &1\leq i\leq s;\\
\vec{q'_i}(\bn)&:=(0,\dots,0,p'_{i,1}(\bn),\dots,p'_{i,\ell'}(\bn)), &   \text{for } &1\leq i\leq s'.
\end{align*}
We also let $\vec{q_0}=\vec{q'_0}:=0$. Finally, we define  the functions $g_i\in L^\infty(\nu)$, $0\leq i\leq s,$ and $g'_i\in L^\infty(\nu)$, $0\leq i\leq s'$, by
$$
g_i:=\one_X f_i+\one_{X'} \ \text{ and }\ g'_i:=\one_X+\one_{X'} f'_i.
$$
Then  we have
$$
\int_Y\;\prod_{i=0}^s S_{\vec{q_i}(\bn)}g_i\cdot\prod_{i=0}^{s'}
S_{\vec{q'_i}(\bn)}g'_i\;d\nu=\frac 12(a(\bn)+a'(\bn)), \quad \bn\in\N^d.
$$
This completes the proof of the linearity of the space $\mathcal{MC}_\text{d,pol}$.

The inclusion $\overline{\mathcal{N}_{d}}\subset \overline{\mathcal{MC}}_{\text{d,lin}}$
follows from Proposition~\ref{prop:represent}.
The inclusion $\overline{\mathcal{MC}}_{\text{d,lin}}\subset \overline{\mathcal{MC}}_{\text{d,pol}}$
is obvious.
The inclusion $\overline{\mathcal{MC}}_{\text{d,pol}}\subset \overline{\mathcal{N}_{d}}$
follows from Theorem~\ref{th:correl-poly}. This completes the proof of Theorem~\ref{T:4}. \qed

\section{Convergence criteria for weighted averages}
In this short section we use the machinery developed in the previous  sections  in order to prove the convergence criteria stated in Sections~\ref{SS:2.2} and \ref{SS:2.3}.
\subsection{Proof of Theorem~\ref{th:weighted-L2}}
\label{subsec:Proof-Th2}
Let  $d,\ell,s, t \in\N$,  $(X,\mu,T_1,\dots,T_\ell)$ be a system,    $f_1,\dots,f_s$ be functions in $ L^\infty(\mu)$, and   $\vec{p_i}\colon \N^d\to\Z^\ell$, $i=1,\dots,s$, be polynomial mappings of degree at most $t$. We can assume that $\norm{f_i}_{L^\infty(\mu)}\leq 1$ for $i=1,\ldots, s$.

 Let $\delta>0$ and $k=k(d,\ell,s,t)$
be given by Proposition~\ref{prop:unif-in-norm} and suppose  that the bounded sequence $w\colon\N^d\to\C$ is  $k$-regular. As remarked there, if all the polynomials are linear, then we can take $k=s$.
By Theorem~\ref{th:regular}, the sequence $w$ can be written as $w=\psi+u$ where $\psi$ is a $k$-step nilsequence in $d$ variables and $u\in \ell^\infty(\N^d)$ satisfies  $\norm{u}_{U^{k+1}(\N^d)}<\delta$.
By Proposition~\ref{prop:weight-nil}, the limit
$$
\limav \psi(\bn)\cdot T_{\vec{p_1}(\bn)}f_1\cdot\ldots\cdot T_{\vec {p_s}(\bn)}f_s
$$
exists in $L^2(\mu)$ and by Proposition~\ref{prop:unif-in-norm} for every F\o lner sequence $\bI=(I_j)_{j\in \N}$ we have
$$
\Lim\Bigl\Vert\Av_{\bI}\, u(\bn)\cdot T_{\vec{p_1}(\bn)}f_1\cdot\ldots\cdot
T_{\vec{p_s}(\bn)}f_s\Bigr\Vert_{L^2(\mu)}
\leq  4\delta.
$$
It follows that for all sufficiently large $j,j'\in \N$, the difference between
\begin{equation}
\label{eq:av:L2}
\frac 1{|I_j|}\sum_{\bn\in I_j}
u(\bn)\,  T_{\vec{p_1}(\bn)}f_1\cdot\ldots\cdot
T_{\vec{p_s}(\bn)}f_s
\end{equation}
and the similar average on $I_{j'}$ has a norm in $L^2(\mu)$ bounded by $8\delta$.
Therefore, the averages~\eqref{eq:av:L2} form a Cauchy sequence and thus converge in $L^2(\mu)$.

Furthermore, if  we assume  that the averages of $w(\bn)\psi(\bn)$ converge to $0$ for every $k$-step nilsequence $\psi$ in $d$ variables, then by Corollary~\ref{cor:unif-approx2} we have $\norm w_{U^k(\N^d)}=0$ and the averages~\eqref{eq:av:L2} converge to $0$ in $L^2(\mu)$ by  Proposition~\ref{prop:unif-in-norm}. This completes the proof.
\qed

\subsection{Proof of Theorem~\ref{th:Cesaro}}
\label{subsec:Proof-Th4}
Let  $k=k(d,\ell,s,t)$ be as in  Theorem~\ref{th:correl-poly} and the sequence $a\colon \N^d\to \C$ be defined by
$$
a(\bn):=\int f_0\cdot  T_{\vec{p_1}(\bn)}f_1\cdot\ldots\cdot T_{\vec{p_s}(\bn)}f_s\,d\mu, \quad \bn\in \N^d.
$$
 By Theorem~\ref{th:correl-poly} $a$ is an approximate  $k$-step nilsequence in $d$ variables and  by hypothesis~\eqref{eq:Cesaro-av1} the averages~\eqref{eq:Cesaro-av2} converge. This proves Part~\eqref{it:Cesaro1} of Theorem~\ref{th:Cesaro}. Furthermore, we deduce that if the averages~\eqref{eq:Cesaro-av1} converge to $0$ for every nilsequence $\psi$ in $d$ variables, then the averages~\eqref{eq:Cesaro-av2} converge  to $0$.

Next we prove Part~\eqref{it:Cesaro2}.
Let  $k=k(2d,\ell,2s-1,t)$ be as in  Theorem~\ref{th:correl-poly}
and let the sequence $b\colon\N^d\times\N^d\to \C$ be defined by
$$
b(\bn,\bn'):=\int  T_{\vec{p_1}(\bn)}f_1\cdot\ldots\cdot T_{\vec{p_s}(\bn)}f_s\cdot  T_{\vec{p_1}(\bn')}\overline{f}_1\cdot\ldots\cdot T_{\vec{p_s}(\bn')}\overline{f}_s\,d\mu, \quad \bn,\bn'\in\N^d.
$$
By Theorem~\ref{th:correl-poly} the sequence $b$ is an approximate $k$-step nilsequence in $2d$ variables. Hence, by hypothesis~\eqref{eq:Cesaro-av3}, the averages
$$
\frac 1{N^dN'^d}\sum_{\bn\in[N]^d}\sum_{\bn'\in[N']^d}
w(\bn)\, w(\bn') \, b(\bn,\bn')
$$
converge to some limit $L$  when $N$ and $N'$ tend to $+\infty$. Let $N_0$ be such that the difference between this average and $L$ is bounded in absolute value by $\ve$ for all $N,N'>N_0$.
Expanding the square
$$
\Bigl\Vert \frac{1}{N^d} \sum_{\bn\in[N]^d} w(\bn)  \cdot  T_{\vec{p_1}(\bn)}f_1\cdot\ldots\cdot T_{\vec{p_s}(\bn)}f_s
-\frac{1}{N'^d}\sum_{\bn'\in[N']^d} w(\bn') \cdot   T_{\vec{p_1}(\bn')}f_1\cdot\ldots\cdot T_{\vec{p_s}(\bn')}f_s\Bigr\Vert_{L^2(\mu)}^2
$$
we obtain $4$ terms of this form with alternate signs, and thus this square is bounded by $4\ve$. Hence,  the averages in~\eqref{eq:Cesaro-av4} form a Cauchy sequence in $L^2(\mu)$ and thus converge. This proves Part~\eqref{it:Cesaro2} of Theorem~\ref{th:Cesaro}. As above, if the limit~\eqref{eq:Cesaro-av3} is zero for all  $k$-step nilsequences $\psi$ in $2d$ variables, then  the limit~\eqref{eq:Cesaro-av4} is equal to $0$.

 In order to get the last part of Theorem~\ref{th:Cesaro} we argue as before and use the last parts of Theorems~\ref{th:correl-poly} and \ref{th:correl-lin}.\qed

\section{Applications to  arithmetic weights}
\label{sec:arithmrtic}
In this section we prove the main results stated in  Section~\ref{SS:arithmetic}.
\subsection{Proof of Theorem~\ref{th:cesaro-mult}}
\label{subsec:mult}
Theorem~\ref{th:cesaro-mult} follows immediately
by the next result which shows that  hypothesis~\eqref{it:Cesaro2} of Theorem~\ref{th:Cesaro} is satisfied.
\begin{theorem}
\label{th:double-av1}
Let $d\in \N$ and $\phi\colon \N^d\to \C$ be a good  multiplicative function.
Then for every nilsequence $\psi\colon \N^d\to \C$ the limit
\begin{equation}
\label{eq:double-av1}
\lim_{N_1,\ldots, N_d\to+\infty}\frac 1{N_1\cdots N_d}
\sum_{\bn \in [N_1]\times \cdots \times [N_d]}
\phi(\bn)\, \psi(\bn)
\end{equation}
exists.
Moreover, if  the multiplicative function $\phi$  is aperiodic, then the limit
is equal to $0$ for all nilsequences $\psi$ in $d$ variables.
\end{theorem}
\begin{remarks}
$\bullet$ One can also deduce variants of the previous result that deal with polynomial sequences on nilmanifolds;  indeed, any such sequence can be represented as a   linear sequence on a different nilmanifold \cite[Section~2.11]{L05b}.

$\bullet$ The case where $d=1$ can be deduced from
Theorem~6.1 in \cite{FH15a}, but due to the finitary  nature  of the statement in \cite{FH15a}
the argument used there is vastly more
complicated.
\end{remarks}
We prove Theorem~\ref{th:double-av1}  in  Section~\ref{sec:multiplicative}.

\subsection{Proof of Theorems~\ref{T:convarith} and~\ref{T:recarith}}\label{SS:pf}
For $b\in \N$ we let $\zeta$ be a root of unity of order $b$ and let $f_b$  be the multiplicative function
 defined by $f_b(p^k)=\zeta$ for all primes $p$ and all  $k\in \N$.\footnote{In order to get results when
 the set $S_{a,b}$ is defined by counting prime factors with multiplicity, we define the completely multiplicative
 function $f_b$ by $f_b(p^k)=\zeta^k$ for all primes $p$ and all $k\in \N$.}
 Note that
\begin{equation}\label{E:key}
 \one_{S_{a,b}}(n)= \frac{1}{b}\sum_{j=0}^{b-1} \zeta^{-aj}(f_{b}(n))^j, \quad n\in \N.
 \end{equation}
 It can be seen using Theorem~\ref{T:Halasz} (see~\cite[Proposition~2.10]{FH15b})
  that for $j=1,\ldots, b-1$ the multiplicative function  $f^j$ is aperiodic.

 Let
 $$
 V_\bn:= T_{\vec{p_1}(\bn)}f_1\cdot\ldots\cdot T_{\vec {p_s}(\bn)}f_s, \quad \bn \in \N^d.
 $$
Recall that $S:=(S_{a_1,b_1}+c_1)\times \cdots \times (S_{a_d,b_d}+c_d)$. Using \eqref{E:key} we get  that the averages
 $$
  \frac{1}{N^d}\sum_{\bn\in S\cap [N]^d}    V_\bn
 $$
  are asymptotically equal (meaning, the relevant difference converges to $0$ in $L^2(\mu)$)  to the averages
  \begin{equation}\label{E:zi}
  \frac{1}{N^d}\sum_{\bn\in [N]^d} \prod_{i=1}^d  \frac{1}{b_i}\sum_{j=0}^{b_i-1} \zeta^{-a_ij}(f_{b_i}(n_i))^j\cdot V_{\bn+\bc}
  \end{equation}
  where $\bc:=(c_1,\ldots, c_d)$ and as usual $n_1,\ldots, n_d$ denote the coordinates of the vector $\bn\in \N^d$.
  Since for $i=1,\ldots, d$ and $j=1,\ldots, b_i-1$ the multiplicative function $f_{b_i}^j$ is aperiodic,
  we get that for
 $j_i\in\{0,\ldots, b_i-1\}$, $i=1,\ldots, d$,  the multiplicative function
 $(n_1,\ldots,n_d)\mapsto  \prod_{i=1}^d (f_{b_i}(n_i))^{j_i}$ is aperiodic unless $j_1=\cdots=j_d=0$.
  Keeping this in mind,  expanding the product in \eqref{E:zi},  and using the second part of Theorem~\ref{th:cesaro-mult},
  we deduce  that the  averages \eqref{E:zi} are asymptotically equal in $L^2(\mu)$ to the averages
  $$
  \prod_{i=1}^d \frac{1}{b_i} \cdot \frac{1}{N^d}\sum_{\bn\in [N]^d}  V_{\bn}.
$$
  Furthermore,  the previous argument, applied for $V_\bn:=1,$  $\bn\in [N]^d$, gives that
  $$
  \lim_{N\to+\infty}\frac{|S\cap [N]^d|}{N^d}=\prod_{i=1}^d \frac{1}{b_i}.
  $$
 Combining the above, we deduce that the difference
 $$
 \frac{1}{|S\cap [N]^d|}\sum_{\bn\in S\cap [N]^d}    V_\bn- \frac{1}{N^d}\sum_{\bn\in [N]^d}  V_{\bn}
 $$
 converges to $0$ in $L^2(\mu)$. Using  this, in conjuction with  Theorem~\ref{T:Walsh},
 completes the proof of Theorem~\ref{T:convarith}.

 Using the previously established identity for $f_i:={\bf 1}_A$, $i=1,\ldots,d$,  where $A\in \CX$ is a set of positive measure,
 we deduce that
 \begin{align*}
\lim_{N\to + \infty} \frac{1}{|S\cap [N]^d|}\sum_{\bn \in S\cap [N]^d} \mu(A\cap T_{-\vec{p_1}(\bn)}A\cap \cdots\cap
T_{-\vec {p_s}(\bn)}A)=\\
\lim_{N\to + \infty} \frac{1}{N^d}\sum_{\bn \in  [N]^d} \mu(A\cap T_{-\vec{p_1}(\bn)}A\cap \cdots\cap T_{-\vec {p_s}(\bn)}A)>0,
\end{align*}
 where the positiveness of the last limit follows from the multiparameter polynomial Szemer\'edi theorem~\cite[Theorem~0.9]{BMcC}.
This proves  Theorem~\ref{T:recarith}. \qed

\subsection{Correlations of multiplicative functions with nilsequences}
\label{sec:multiplicative}

In this subsection we prove Theorem~\ref{th:double-av1} for $d=2$; the general
 case is identical to this one modulo changes in notation.
We begin with some preliminaries.

\subsubsection{Some classical facts about commutators}

\begin{lemma}
\label{lem:commut1}
Let $G$ be a group and $i,j\in\N$. Then the commutator map $(g,h)\to[g,h]$ maps $G_i\times G_j$ to $G_{i+j}$. Moreover, it induces a bi-homomorphism from $(G_i/G_{i+1})\times(G_j/G_{j+1})$ to $G_{i+j}/G_{i+j+1}$.
\end{lemma}
\begin{proof}[Sketch of the proof]
The first statement follows by induction from the classical three subgroups lemma:

{\em Let $H,K, L\subset G$ and $N$ be a normal subgroup of $G$.  If $[[H,K],L]\subset N$
 and $[[L,H],K]\subset N$,  then $[[K,L],H]\subset N$.}

The second statement follows from the identity:
\begin{equation}
\label{eq:commut}
[xy,z]=x[y,z]x\inv \cdot[x,z] = [x,[y,z]]\cdot [y,z]\cdot [x,z]
\quad \text{for all } \ x,y,z\in G.
\end{equation}
\end{proof}

\begin{lemma}
\label{lem:commut2}
Let $G$ be a group, $H^{(1)},\dots,H^{(k)}$ and $Q^{(1)},\dots,Q^{(\ell)}$ be  normal subgroups of $G$,  $H=H^{(1)}\cdot\ldots\cdot H^{(k)}$ and $Q=Q^{(1)}\cdot\ldots\cdot Q^{(\ell)}$. Then the commutator group $[H,Q]$ is the product of  the  groups $[H^{(i)},Q^{(j)}]$ for $i=1,\dots, k$ and $j=1,\dots,\ell$.
\end{lemma}
\begin{proof}
All the groups  $[H^{(i)},Q^{(j)}]$ are normal and included in $[H,Q]$, thus   it suffices to prove that $[H,Q]$ is contained in the product of these groups.
If $h_i\in H^{(i)}$ for $i=1,2$ and $q\in Q^{(1)}$, it follows from~\eqref{eq:commut} and from the normality of  $[H^{(1)},Q^{(1)}]$  that
$[h_1h_2,q]\in[H^{(1)},Q^{(1)}]\cdot[H^{(2)},Q^{(1)}]$.
 This proves the result in the case $k=2$, $\ell=1$. The result for $k$ arbitrary and $\ell=1$ follows by induction on $k$. The result for $k,\ell$ arbitrary follows by exchanging the roles of $H$ and $Q$.
\end{proof}

\subsubsection{Some reductions}
We prove  Theorem~\ref{th:double-av1} for $d=2$.
Recall that  $\psi$ is an $s$-step nilsequence in  $2$ variables  defined by
$$
\psi(n_1,n_2):=\Psi(\tau_1^{n_1}\tau_2^{n_2}\cdot e_X), \quad n_1,n_2\in \N,
$$ where $X=G/\Gamma$ is a   nilmanifold, $\Psi\in C(X)$,  and $\tau_1,\tau_2$ are two commuting elements of $G$.
It is known \cite{L05b,Le91}  that the closure $X'$ of $\{\tau_1^{n_1}\tau_2^{n_2}\cdot e_X\colon n_1,n_2\in\N\} $ in $X$ is a sub-nilmanifold of $X$.
Substituting  $X'$ for $X$ we can assume that
$\{\tau_1^{n_1}\tau_2^{n_2}\cdot e_X\colon n_1,n_2\in\N\} $
 is dense in $X$; it then  follows by \cite[Theorem~1.4]{L05b} that  it is  equidistributed in $X$.


We can assume  without loss of generality that $G$ is spanned by the connected component  of the unit element and the elements $\tau_1$ and $\tau_2$. This implies   that all commutator subgroups $G_i$, $i\geq 2$, are connected (see for example~\cite[Theorem~4.1]{BHK05}).

 Suppose that $s\geq 2$ and that $G$ is  $s$-step  but not $(s-1)$-step nilpotent. Then the group $K_s:=G_s/(G_s\cap\Gamma)$ is a finite dimensional torus, sometimes called the \emph{vertical torus} \cite{GT12c},
and acts freely on $X$. We denote this action by $(u,x)\mapsto u\cdot x$ for $u\in K_s$ and $x\in X$. Let $\wh{K_s}$ be the dual group of $K_s$, that is, the group of continuous group homomorphisms from $K_s$ to the circle group $\CS^1$. If for some $\chi\in\wh{K_s}$ the
 function $\Psi\in C(X)$    satisfies
$\Psi(u\cdot x)=\chi(u)\Psi(x)$ for every $u\in K_s$ and $x\in X$, then  it is  called a \emph{nilcharacter} of $X$  with \emph{vertical frequency} $\chi$.
The linear span of nilcharacters is dense in $C(X)$ for the uniform norm. Therefore, it suffices to prove \eqref{eq:double-av1} when the function $\Psi$ defining the nilsequence $\psi$
  is a nilcharacter.

If the vertical frequency $\chi$ of $\Psi$ is the trivial character of $K_s$, then the function $\Psi$ factorizes through the quotient of $X$ under the action of this group. This quotient is the $(s-1)$-step nilmanifold $X'=G/(G_s\Gamma)=(G/G_s)/((\Gamma\cap G_s)/G_s)$.  Writing $\tau'_1,\tau'_2$ for the images of $\tau_1,\tau_2$ in $G/G_s$, we have that $\psi(n_1,n_2)=\Psi'(\tau_1'^{n_1}\tau_2'^{n_2}\cdot e_{X'})$ for some $\Psi'\in C(X')$.

Iterating this procedure, we reduce matters to considering the following
two cases: (i) $X$ is a $1$-step nilmanifold, and (ii) $X$ is an $s$-step nilmanifold for some $s\geq 2$ and $\Psi$ is a nilcharacter with vertical frequency $\chi\neq 1$.

If $(i)$ holds, then $X$ is a compact Abelian group and
we can further  reduce matters  to the case where  $\Psi$ is a character of $X$.  Then the
average in~\eqref{eq:double-av1} can be rewritten as
\begin{equation}
\label{eq:double-av3}
\Bigl(\frac 1{N_1}\sum_{n_1=1}^{N_1}\phi_1(n_1)\, \e(n_1t_1)\Bigr)\cdot
\Bigl(\frac 1{N_2}\sum_{n_2=1}^{N_2}\phi_2(n_2)\, \e(n_2t_2)\Bigr)
\end{equation}
for some $t_1,t_2\in\R$ and good multiplicative functions $\phi_1,\phi_2\colon \N\to \C$. If $t_1$ or $t_2$ is irrational, then the limit is equal to $0$ by a Theorem of Daboussi~\cite{D74} (see also~\cite{DD74,DD82}).
If both $t_1$ and $t_2$ are rational, then  the limit exists since by hypothesis the multiplicative functions $\phi_1$ and $\phi_2$ are  good. Furthermore, this limit is equal to $0$ if either $\phi_1$ or $\phi_2$ is aperiodic, that is, if $\phi$ is aperiodic.

Hence, it  suffices to consider  case $(ii)$, that is, we can  assume  that $X$ is an $s$-step nilmanifold for some $s\geq 2$ and $\Psi$ is a nilcharacter with vertical frequency $\chi\neq 1$.
Replacing $X$ with its quotient by the kernel of $\chi$, we are reduced to the case where $G_s=\CS^1$ and $\Psi$ has vertical frequency $1$, that is,
\begin{equation}\label{E:assumption}
\Psi(u\cdot x)=u\, \Psi(x) \quad \text{for every }x\in X\ \text{ and }\ u\in G_s=\CS^1.
\end{equation}
It therefore suffices to prove:
\begin{proposition}[Theorem~\ref{th:double-av1}-Nilcharacter form]
\label{pr:double-av2}
Let $X=G/\Gamma$ be an $s$-step nilmanifold for some $s\geq 2$. Suppose   that $G_s=\CS^1$, $\Psi\in C(X)$ satisfies \eqref{E:assumption},  and $\tau_1,\tau_2$ are two commuting elements of $G$ such that the sequence $(\tau_1^{n_1}\tau_2^{n_2}\cdot e_X)_{n_1,n_2\in\N} $
is equidistributed in $X$. Then
 for all  multiplicative functions $\phi_1,\phi_2\colon \N\to \C$ that are bounded by $1$  we have
\begin{equation}
\label{eq:double-av2}
\lim_{N_1,N_2\to+\infty}\frac 1{N_1N_2}
\sum_{n_1\in[N_1],\ n_2\in[N_2]}
\phi_1(n_1)\, \phi_2(n_2)\,\Psi(\tau_1^{n_1}\tau_2^{n_2}\cdot e_X)=0.
\end{equation}
\end{proposition}
We proceed now to establish Proposition~\ref{pr:double-av2}.
\subsubsection{A variant of Kat\'ai's Lemma}
We use  a two dimensional variant  of  a result of
Kat\'ai (\cite{K86}, see also~\cite{D74}); we omit its proof since it is identical with its one dimensional version modulo changes in notation.
\begin{proposition}
\label{prop:Katai}
Let $P_0\in\N$, $\phi_1,\phi_2\colon \N \to \C$  be  multiplicative functions, bounded by $1$,   and $a\in \ell^\infty(\N^2)$. Suppose that
$$
\lim_{N_1,N_2\to+\infty}
\frac 1{N_1N_2}\sum_{n_1\in[N_1],\; n_2\in [N_2]}
 a(p_1n_1,p_2n_2)\, \overline{a(p'_1n_1,p'_2n_2)}=0
$$
 for all  distinct primes $p_1,p_2,p'_1,p'_2\geq P_0$. Then
 $$
\lim_{N_1,N_2\to+\infty}
\sup_{\phi_1,\phi_2}\Big|\frac 1{N_1N_2}\sum_{n_1\in [N_1],\; n_2\in  [N_2]}
\phi_1(n_1)\, \phi_2(n_2)\,a(n_1,n_2)\Big|=0
$$
where the sup is taken over all $\phi_1,\phi_2\colon \N\to \C$ that are multiplicative and bounded by $1$.
\end{proposition}

 \subsubsection{The nilmanifold $Y$}
 In the sequel, $p_1,p_2,p'_1,p'_2$ are fixed distinct primes.  By Proposition~\ref{prop:Katai}, in order to prove Proposition~\ref{pr:double-av2},
 it suffices to show that
 \begin{equation}
\label{eq:Katai}
\lim_{N_1,N_2\to+\infty}
\frac 1{N_1N_2}\sum_{n_1\in[N_1], n_2\in[N_2]}
 \Psi(\tau_1^{p_1n_1}\tau_2^{p_2n_2}\cdot e_X)\, \overline{\Psi(\tau_1^{p'_1n_1}\tau_2^{p'_2n_2}\cdot e_X)}=0
 \end{equation}
where $\Psi$ satisfies $\eqref{E:assumption}$.
We let
$$
Y:=\overline{\bigl\{
\bigl( \tau_1^{p_1n_1}\tau_2^{p_2n_2}\cdot e_X\;,\; \tau_1^{p'_1n_1}\tau_2^{p'_2n_2}\cdot e_X\bigr)
\colon n_1,n_2\in\N\bigr\}}\subset X\times X.
$$
Then $Y$ is a sub-nilmanifold of $X\times X$ and  by \cite[Theorem~1.4]{L05b} the $2$-variable sequence above is equidistributed in $Y$.
Writing $m_Y$ for the Haar measure of $Y$, the limit in~\eqref{eq:Katai} is equal to
\begin{equation}
\label{eq:int-Y}
\int_Y\Psi(x)\, \overline{\Psi(x')}\,dm_Y(x,x').
\end{equation}
Hence, it remains to  show that this integral is zero.

Let $H$ be the smallest closed subgroup of $G$ containing $\Gamma\times\Gamma$ and the shift elements $(\tau_1^{p_1},\tau_1^{p'_1})$ and $(\tau_2^{p_2},\tau_2^{p'_2})$. We claim that
$$
Y=H/(\Gamma\times\Gamma).
$$
Indeed, by the definition of $H$ and $Y$, we have   $H\cdot(e_X,e_X)\supset Y$.
Furthermore, by the remark  following~\cite[Theorem~2.21]{L05}, we have
 $Y=H_1\cdot(e_X,e_X)$ for some closed subgroup $H_1$ of $G\times G$ containing the shift elements  $(\tau_1^{p_1},\tau_1^{p'_1})$ and $(\tau_2^{p_2},\tau_2^{p'_2})$.  Since $Y$ is compact,
 $H_1\cap(\Gamma\times\Gamma)$
 is cocompact in $H_1$ and thus $H_2:=H_1\cdot(\Gamma\times\Gamma)$ is a closed subgroup of $G$. Since $H_2$ is a closed subgroup that contains  the shift elements and $\Gamma\times\Gamma$, we have  $H_2\supset H$, hence
 $H\cdot(e_X,e_X)\subset H_2\cdot(e_X,e_X)=H_1\cdot(e_X,e_X)=Y$. Therefore, $H\cdot(e_X,e_X)= Y$,
 which implies that  $Y=H/(\Gamma\times\Gamma)$. This proves the claim.
\subsubsection{Projection on the Kronecker factor}
We denote by  $Z$ the  compact Abelian group $G/(G_2\Gamma)$
and let  $\pi\colon X\to Z$ and $p\colon G\to Z$ be the natural projections. For $i=1,2$ let $\alpha_i:=p(\tau_i)$ and $Z_i$ be the closed subgroup of $Z$ spanned by $\alpha_i$. Since $\{\tau_1^{n_1}\tau_2^{n_2}\cdot e_X\colon n_1,n_2\in\N\} $
 is dense in  $X$, we have that  $\{\alpha_1^{n_1}\alpha_2^{n_2}\cdot e_X\colon n_1,n_2\in\N\} $
 is dense in  $Z$  and $Z=Z_1Z_2$.
 For $i=1,2$ we let
 $$
G^{(i)}:=p\inv(Z_i).
$$
Then $G^{(i)}$ is a closed subgroup of $G$ containing $\Gamma$ and $G_2$, hence normal in $G$, and
$$
G=G^{(1)}G^{(2)}.
$$

Let
$$
W:=(p\times p)(H)=(\pi\times\pi)(Y).
$$
By the definition of $Y$, we have that $W$ is the closure in $Z\times Z$ of
$$
\bigl\{
\bigl( \alpha_1^{p_1n_1}\alpha_2^{p_2n_2}\cdot e_X\;,\; \alpha_1^{p'_1n_1}\alpha_2^{p'_2n_2}\cdot e_X\bigr)
\colon n_1,n_2\in\N\bigr\}
$$
and thus
\begin{equation}
\label{eq:W}
W=\bigl\{ \bigl(z_1^{p_1}z_2^{p_2}\;,\; z_1^{p'_1}z_2^{p'_2}\bigr)\colon z_1\in Z_1,\ z_2\in Z_2\bigr\}.
\end{equation}

\subsubsection{Starting the induction}
For $i=1,2$ let $g_i\in G^{(i)}$. Then by~\eqref{eq:W} we  have
$$
(p\times p)\bigl(g_1^{p_1}, g_1^{p'_1}
\bigr)\in W\ \text{ and }\
(p\times p)\bigl(g_2^{p_2}, g_2^{p'_2}
\bigr)\in W.
$$
We have $(p\times p)(H)=W$,   the kernel of $p\times p$ is $(G_2\times G_2)(\Gamma\times\Gamma)$, and $\Gamma\times\Gamma\subset H$, thus
$$
(g_1^{p_1}, g_1^{p'_1})\in H(G_2\times G_2)\
\text{ and }\
(g_2^{p_2}, g_2^{p'_2})\in H(G_2\times G_2).
$$
For $i=1,2$, let  $h_i\in G^{(i)}$.
By~\eqref{eq:commut} we have
$$
\bigl( H(G_2\times G_2)\bigr)_2\subset H_2(G_3\times G_3)
$$
 and thus
$$
\bigl[  (g_1^{p_1}, g_1^{p'_1})\;,\;(h_1^{p_1}, h_1^{p'_1})\bigr], \
\bigl[  (g_1^{p_1}, g_1^{p'_1})\;,\; (h_2^{p_2}, h_2^{p'_2})\bigr],  \text{ and }
\bigl[  (g_2^{p_2}, g_2^{p'_2})\;,\;(h_2^{p_2}, h_2^{p'_2})\bigr]
$$
belong to $H_2(G_3\times G_3)$. By Lemma~\ref{lem:commut1}, these elements are equal modulo $G_3\times G_3$ to
$$
\bigl([g_1,h_1]^{p_1^2},[g_1,h_1]^{p_1'^2}\bigr),\
\bigl([g_1,h_2]^{p_1p_2}, [g_1,h_2]^{p'_1p'_2}\bigr),
\text{ and }
\bigl([g_2,h_2]^{p_2^2}, [g_2,h_2]^{p_2'^2}\bigr)
$$
respectively.

Henceforth, for $i=1,2$ we denote $(G^{(i)})_2$ by $G^{(i)}_2$. For $u,v\in G_2$ we have $(uv)^{p_1^2}=u^{p_1^2}v^{p_1^2}$ mod $G_3$ and $(uv)^{p_1'^2}=u^{p_1'^2}v^{p_1'^2}$ mod $G_3$ and thus the set
$$
L:=\{u\in G^{(1)}_2\colon (u^{p_1^2}, u^{p_1'^2})\in H_2(G_3\times G_3)\}
$$
is a subgroup of $G^{(1)}_2$. By the previous discussion, for $g_1,h_1\in G^{(1)}$ the set $L$
 contains $[g_1,h_1]$ and thus it is equal to $G^{(1)}_2$. Hence,
$$
(u^{p_1^2}, u^{p_1'^2})\in H_2(G_3\times G_3)\  \text{ for every }\ u\in G^{(1)}_2.
$$
In the same way, we have that
\begin{gather*}
(u^{p_2^2}, u^{p_2'^2})\in H_2(G_3\times G_3)\  \text{ for every }\ u\in G^{(2)}_2;\\
(u^{p_1p_2}, u^{p'_1p'_2})\in H_2(G_3\times G_3) \ \text{ for every }\ u\in [G^{(1)},G^{(2)}].
\end{gather*}

\subsubsection{The induction}
By induction on $r$ we show:
\begin{lemma}
\label{lem:Gs}
For $2\leq r\leq s$ and $i_1,\dots,i_r\in\{1,2\}$ let
\begin{gather*}
s_i:=|\{j\colon 1\leq j\leq r,\ i_j=i\}|\ \text{ for }i=1,2;\\
G^{(i_1,\dots,i_r)}:=
\bigl[G^{(i_1)}, \bigl[ G^{(i_2)},\dots,\bigl[G^{(i_{r-1})},G^{(i_r)}\bigr]\dots\bigr]\bigr]
\subset G_r.
\end{gather*}
Then, for every
$u\in G^{(i_1,\dots,i_r)}$ we have
$\bigl( u^{p_1^{s_1}p_2^{s_2}},
u^{p_1'^{s_1}p_2'^{s_2}}\bigr)\in H_r(G_{r+1}\times G_{r+1})$.
\end{lemma}

\begin{proof}
For $r=2$ this was proved in the preceding subsection and the inductive step is proved by the same method.
\end{proof}
In the sequel we only use this lemma for $r=s$.

\subsubsection{Conclusion of the proof}
We argue by contradiction. Suppose  that
$$
I:=\int \Psi(x)\, \overline{\Psi}(x')\, dm_Y(x,x')\neq 0.
$$
Recall that  $G_s=\CS^1$ and that $\Psi$ satisfies \eqref{E:assumption}.

Let  $i_1,\dots,i_s\in\{1,2\}$ and $s_i=|\{j\colon 1\leq j\leq s,\ i_j=i\}|$ for $i=1,2$.
Since $G$ is $s$-step nilpotent, the subgroup  $G_{s+1}$ is trivial. Hence,  by Lemma~\ref{lem:Gs},
for every $u\in G^{(i_1,\dots,i_s)}$,
 we have $\bigl( u^{p_1^{s_1}p_2^{s_2}},
u^{p_1'^{s_1}p_2'^{s_2}}\bigr)\in H_s$.
Therefore, the measure $m_Y$ is invariant under translation by $\bigl( u^{p_1^{s_1}p_2^{s_2}},
u^{p_1'^{s_1}p_2'^{s_2}}\bigr)$. Hence,
$$
I= \int \Psi(u^{p_1^{s_1}p_2^{s_2}}\cdot x)\, \overline{\Psi}(u^{p_1'^{s_1}p_2'^{s_2}}\cdot  x')\,dm_Y(x,x')=
u^{ p_1^{s_1}p_2^{s_2} - p_1'^{s_1}p_2'^{s_2}}\cdot I.
$$
Thus $u^{ p_1^{s_1}p_2^{s_2} - p_1'^{s_1}p_2'^{s_2}}=1$ for every $g\in G^{(i_1,\dots,i_s)}$ .
Since $G^{(i_1,\dots,i_s)}$ is a subgroup of the torus $G_s$ and $p_1,p'_1,p_2,p'_2$ are distinct primes, it  follows that
 the group
$G^{(i_1,\dots,i_s)}$ is finite.

Furthermore, since $G=G^{(1)}G^{(2)}$,   using Lemma~\ref{lem:commut2} and induction, we have that the group
 $G_s$ is the product of  the groups $G^{(i_1,\dots,i_s)}$ for $i_1,\dots,i_s\in\{1,2\}$  and thus is finite, a contradiction since $G_s=\mathcal{S}^1$. This completes the proof of Proposition~\ref{pr:double-av2} and hence of  Theorem~\ref{th:double-av1}.\qed

\section{Applications to Hardy field weights}


\subsection{Hardy field functions} \label{SS:Hardy}
 Let $B$ be the collection of equivalence classes of real
valued functions defined on some half line $[c,+\infty)$, where we identify two functions if they agree
on some half line.  A \emph{Hardy field} $H$ is a subfield of the ring $(B,+,\cdot)$
 that is closed under differentiation.
A \emph{Hardy field function} is a function that
belongs to some Hardy field.
An example of a  Hardy field consists of all \emph{logarithmic-exponential functions}, meaning, all functions defined on
some half line $[c,+\infty)$ by a finite combination of the symbols $+,-,\times, :,\log,\exp$ operating on
the real variable t and on real constants. Examples include functions of the form $t^a(\log t)^b$ for every $a,b\in \R$. An important property of
Hardy field functions is
that we can relate their growth rates with the growth rates of
their derivatives.  The reader can find further discussion about Hardy fields  in \cite{Bo81, Bo94, BKQW94, F09}
and the references therein.

Let $f$ be a  Hardy field function. We say that it
\begin{enumerate}
\item
has \emph{at most polynomial growth} if $f(t)/t^m\to 0$ for some $m\in \N$;

\item \emph{stays away from polynomials} if
$|f(t)-p(t)|/\log t\to + \infty$ for every $p\in \R[t]$;

\item is \emph{asymptotically polynomial} if $f(t)-p(t)\to 0$ for some $p\in \R[t]$.
\end{enumerate}

For our purposes, the key property of  Hardy field functions that stay away from polynomials is that they satisfy Lemma~\ref{L:basichardy} below.
Examples of  such functions  include:
\begin{itemize}
\item  $t^a$ where $a$ is a positive non-integer;

\item $t^a(\log t)^b$, where
$a> 0$ and  $b\in \R\setminus \{0\}$;

 \item $t^a+(\log t)^b$, where $a\in \R$ and $b>1$.
 \end{itemize}

\subsection{Convergence and recurrence results}
 The main result of this section is the following:
 \begin{theorem}
\label{T:hardymain}
Let $d\in \N$ and $f_1,\ldots, f_d$  be
   Hardy field functions with at most polynomial growth that stay away from polynomials. We define the
   sequence
$w\colon \N^d\to \C$  by
$$
w(\bn):=
e\big(\sum_{i=1}^df_i(n_i)\big), \quad  \bn=(n_1,\ldots, n_d)\in \N^d.
$$
Then for every system $(X,\mu,T_1,\dots,T_\ell)$, functions  $F_1,\ldots, F_s\in L^\infty(\mu)$, and polynomial mappings
 $\vec{p_i}\colon\N^d\to\Z^\ell$, $i=1,\dots,s$,  we have
$$
\lim_{N\to+\infty}\frac 1{N^d}\sum_{\bn \in [N]^d}w(\bn)  \cdot T_{\vec{p_1}(\bn)}F_1\cdot\ldots\cdot T_{\vec{p_s}(\bn)}F_s=0
$$
where the limit is taken in $L^2(\mu)$.
\end{theorem}
\begin{remark}
Related
 work for pointwise convergence when $d=\ell=s=1$  appears in \cite{EK15}.
 \end{remark}
Using the $d=1$ case of the previous result  we deduce in Section~\ref{SS:93} the following:
 \begin{corollary}
\label{C:hardymain}
Let $f$  be a   Hardy field function with at most polynomial growth. Then the
sequence
$w\colon \N\to \C$ defined by $w(n):=e(f(n)), n\in \N,$
is a good universal weight for mean convergence of the
averages \eqref{eq:Cesaro-av4} if and only if  either $f$ is asymptotically polynomial or $f$
stays away from polynomials.
\end{corollary}

 Theorem~\ref{T:hardymain}  follows immediately from Part~$(ii)$ of Theorem~\ref{th:Cesaro} and the next result.

\begin{proposition}
\label{prop:hardyequi}
Let $d\in \N$ and $f_1,\ldots, f_d$  be  Hardy field functions
with at most polynomial growth that stay away from polynomials.
Then for every nilsequence $\psi\colon \N^d\to \C$  we have
$$
\lim_{N_1,\ldots, N_d\to+\infty}\frac 1{N_1\cdots N_d}
\sum_{\bn \in [N_1]\times \cdots \times [N_d]}
e\big(\sum_{i=1}^df_i(n_i)\big)\, \psi(\bn)=0.
$$
\end{proposition}
  We prove Proposition~\ref{prop:hardyequi}
 in Section~\ref{SS:93}.

Next,    we give  some applications.
Note that for $0< a<b< 1/2$ and $t\in [0,1)$  we have  $$
\one_{[a,b]}(\norm{t})=\one_{[a,b]}(t)+\one_{[1-b,1-a]}(t).
$$
Since  $\one_{[c,d]}(t)$   is Riemann integrable,  for all $c,d\in \R$ with  $0\leq c<d<1$,
we have that for every $\varepsilon>0$ there exist ($1$-periodic) trigonometric polynomials $P_1,P_2$ with zero constant terms such that
\begin{equation}\label{E:approximate}
P_1(t)-\varepsilon\leq  \one_{[c,d]}(t)-(d-c)\leq P_2(t)+\varepsilon, \quad t\in [0,1].
\end{equation}
Using  Proposition~\ref{prop:hardyequi} with $2\pi k_if_i$ in place of $f_i$ for $k_i\in \Z$ not all of them zero, for $i=1,\ldots, d$, we deduce using the estimate   \eqref{E:approximate}   the following:
\begin{corollary}
\label{cor:hardyequi'}
Let $d\in \N$ and $f_1,\ldots, f_d$  be  Hardy field functions
with at most polynomial growth that stay away from polynomials.  Let  $a_i,b_i \in \R$ with $0\leq a_i<b_i\leq 1/2$, $i=1,\ldots, d$,
and
$$
S:=\{(n_1,\ldots, n_d)\in \N^d\colon \norm{f_1(n_1)}\in [a_1,b_1], \ldots, \norm{f_d(n_d)}\in [a_d,b_d]\}.
$$
 Then for  every nilsequence $\psi\colon \N^d\to \C$ we have
$$
\lim_{N\to +\infty}\frac 1{|S\cap [N]^d|}
\sum_{n\in S\cap [N]^d}
 \psi(\bn)=\lim_{N\to +\infty}\frac{1}{N^d}\sum_{\bn \in [N]^d} \psi(\bn).
$$
\end{corollary}
We also deduce the following mean  convergence and multiple recurrence result:

\ \\

\begin{theorem}
\label{C:hardyrec}
Let $S\subset \N^d$ be as in Corollary~\ref{cor:hardyequi'}.
Then  the density $d(S)$ of $S$ is $(\prod_{i=1}^d2(b_i-a_i))^{-1}$ and
\begin{enumerate}
\item The sequence $w:={\bf 1}_S$ is a good universal weight for mean convergence of the averages
\eqref{eq:Cesaro-av4} and the limit of these averages is equal to the limit obtained when $w:=d(S)$.

\item For every  $d,\ell,s\in\N$, polynomial mappings  $\vec{p_1},\ldots, \vec{p_s}\colon \N^d\to\Z^\ell$
with zero constant term, system $(X,\mu,T_1,\ldots,T_\ell)$,   and set
$A\in \CX$ with $\mu(A)>0$, we have
\begin{equation}\label{456}
\lim_{N\to +\infty} \frac{1}{N^d}\sum_{\bn \in [N]^d} {\bf 1}_{S}(\bn)\, \mu(A\cap T_{-\vec{p_1}(\bn)}A\cap \cdots\cap T_{-\vec {p_s}(\bn)}A)>0.
\end{equation}
\end{enumerate}
\end{theorem}
\begin{proof}
 If $f$ is a Hardy field function of at most polynomial growth that stays away from polynomials, then the sequence
$(f(n))_{n\in \N}$ is uniformly distributed  $\mod{1}$  (see \cite{Bo94}). The statement about the density of $S$ follows from this fact.

  Using  Theorem~\ref{T:hardymain} with $2\pi k_if_i$ in place of $f_i$ for $k_i\in \Z$ not all of them zero, for $i=1,\ldots, d,$ and the estimate   \eqref{E:approximate},
   we deduce Part $(i)$.

To prove Part $(ii)$ we use Part $(i)$ for $f_1=\cdots=f_s={\bf 1}_A,$ multiply by ${\bf 1}_A$,  and integrate with respect to $\mu$. We deduce that the limit in \eqref{456} is the same as the one obtained when the constant sequence
$d(S)$ takes the place of ${\bf 1}_S$.
The asserted   positiveness then follows from the multiparameter polynomial Szemer\'edi theorem~\cite[Theorem~0.9]{BMcC}.
\end{proof}

\subsection{Proof of Proposition~\ref{prop:hardyequi} and  Corollary~\ref{C:hardymain}}\label{SS:93}
We start with some preliminary facts.
Our assumptions on the functions  $f_1,\ldots, f_d$ are used via the following lemma:
\begin{lemma}\label{L:basichardy}
Let $f$ be a Hardy field function that stays away from polynomials and satisfies   $f(t)/t^m\to 0$ as $t\to +\infty$ for some $m\in \N$. Let    $k\geq m$ be an integer. Then
 there exist real numbers $\alpha_N$ with $\alpha_N\to 0$, polynomials $q_N\in \R[t]$ with $\deg(q_N)<k$, positive integers   $L_N$  with $L_N/N\to 0$  and $L_N^k|\alpha_N|\to +\infty$, such that
$$
f(N+n)=n^k\alpha_N+ q_N(n)+o_{N\to +\infty}(1), \quad n\in [L_N].
$$
\end{lemma}
\begin{proof}
This follows by  noticing that the proof of   \cite[Lemma~3.5]{F09}   applies to all $k\in \N$  with $k\geq m$   and then  following the argument in the proof of  \cite[Lemma~3.4]{F09}.
\end{proof}
In the proof of Proposition~\ref{prop:hardyequi} we use some quantitative equidistribution results
 from \cite{GT12c}. We record here some relevant notions:
 \begin{itemize}
\item If $G$ is a nilpotent group, then    $g\colon \N^d\to G$
 is  a \emph{polynomial sequence} if it has  the form $g(\bn)=\prod_{i=1}^s \tau_i^{p_i(\bn)}$, where for $i=1,\ldots, s$ we have
 $\tau_i\in G$ and $p_i\colon \N^d\to \Z$ are polynomials. The \emph{degree} of the polynomial sequence (with a given representation) is the maximum of the degrees of the polynomials $p_1,\ldots, p_s$.

\item For $N_1,\ldots, N_d \in \N$ we say that  the finite  sequence $(g(\bn)\cdot e_X)_{\bn \in [N_1]\times \cdots \times [N_d]}$ is \emph{ $\delta$-equidistributed}  in the nilmanifold $X$, if  for every Lipschitz function $\Psi\colon X\to \C$ with $\norm\Psi_{\text{Lip}(X)}\leq 1$ and $\int_X\Psi\,dm_X=0$, we have
$$
\Big|\frac{1}{N_1\cdots N_d}\sum_{\bn \in [N_1]\times \cdots \times [N_d]} 
\Psi(g(\bn)\cdot e_X)\Big|\leq\delta.
$$

\item An infinite sequence $(g(\bn)\cdot e_X)_{\bn\in \N^d}$ is \emph{equidistributed in $X$} if  for all   $\Psi\in C(X)$  with $\int_X\Psi\,dm_X=0$  we have (note that the averages below are uniform)
$$
\limav  \Psi(g(\bn)\cdot e_X)=0.
$$
It is
\emph{totally equidistributed in $X$} if the sequence $(\one_{P_1\times \cdots \times P_d}(\bn)\cdot g(\bn)\cdot e_X)_{\bn\in \N}$
is equidistributed in $X$ for
     all infinite  arithmetic progressions $P_1,\ldots, P_d\subset \N$.

\item The \emph{horizontal torus} of the nilmanifold $X=G/\Gamma$ is the compact Abelian group $Z:=G/(G_2\Gamma)$. If $G$ is connected, it is a finite dimensional torus.  A \emph{horizontal
    character} is a continuous group homomorphism $G\to \T$.  It factors through the horizontal torus and induces a character $\eta\colon Z\to \T$; when $G$ is connected,  it is     of the form
      $\bt\mapsto \bk\cdot \bt$, where $\bk\in \Z^s, \bt\in \T^s$, $s:=\text{dim}(Z)$.
   In this case, we define
$\norm{\eta}:=\norm{\bk}_1$, that is, the sum of the absolute values of the coordinates of $\bk$.

\item If $p\colon \N^d\to \T$ has  the form
$p(n_1,\ldots, n_d)=\sum_{j_1,\ldots, j_d} \alpha_{j_1,\ldots, j_d}n_1^{j_1}\cdots n_d^{j_d},$
we define
$$
\norm{p}_{C^\infty[N_1]\times \cdots \times [N_d]}:=\max_{(j_1,\ldots,j_d)\neq (0,\ldots, 0)}N_{1}^{j_1}\cdots N_{d}^{j_d}\norm{\alpha_{j_1,\ldots, j_d}}.
$$

      \item If $X=G/\Gamma$ is a nilmanifold, then $\gamma$ is a  \emph{rational element of}
      $G$ if $\gamma^k\in \Gamma$ for some $k\in \N$.
\end{itemize}

We will use the following quantitative  equidistribution result:
\begin{theorem}[\mbox{\cite[Theorem~8.6]{GT12c} and \cite{GT14}}]
\label{th:Leibman}
Let  $X:=G/\Gamma$ be a nilmanifold  with
 $G$  connected
and simply connected,
$d, t\in \N$, and  $ \ve>0$. There exists
$M:=M(X, d,t,  \ve)>0$ such that the following holds: For all
$N_1,\ldots, N_d\in\N$, greater than $M$,   if $g\colon \N^d\to G$ is a polynomial sequence of degree $t$ and  $(g(\bn)\cdot e_X)_{\bn\in[N_1]\times\cdots\times [N_d]}$ is not   $\ve$-equidistributed in $X$,
  then  there exists
a non-trivial horizontal character $\eta$ such that
$$
0<\norm\eta\leq M\quad \text{ and }\quad \norm{\eta\circ
g}_{C^\infty[N_1]\times\cdots\times [N_d]}\leq M.
$$
\end{theorem}
 \begin{remark} For   every horizontal character $\eta$, the sequence $\eta\circ g$ is a polynomial sequence  in $\T$ of degree at most  $t$.
\end{remark}
We will use the following  elementary result which is a two dimensional variant of \cite[Lemma~3.3]{F09}.
\begin{lemma}\label{L:averaging}
Let $a\in \ell^\infty(\N^2)$ be such that
\begin{equation}\label{E:assum}
\lim_{N,N' \to +\infty}\frac{1}{L_NL'_{N'}}\sum_{\bn\in (N+[L_N])\times (N'+[L'_{N'}])}a(\bn)=0
\end{equation}
for some  sequences  of positive integers $(L_N)_{N\in \N}$, $(L'_N)_{N\in \N}$   that satisfy  $L_N/N\to 0$ and $ L'_N/N\to 0$ as $N\to +\infty$.
Then
$$
\lim_{N,N' \to +\infty} \frac{1}{NN'}\sum_{\bn\in [N]\times [N']} a(\bn)=0.
$$
\end{lemma}
\begin{proof}
 Let the  sequence of positive integers $(k_i)_{i\in \N}$ be defined by
 $k_1:=1$, $k_{i+1}:=k_i+L_{k_i}$, $i\in \N$, and similarly let the sequence
  $(k'_i)_{i\in \N}$ be defined by
 $k'_1:=1$, $k'_{i+1}:=k'_i+L'_{k'_i}$, $i\in \N$.
  For $N\in \N$ let $i_N:=\max\{i\in \N\colon k_i\leq N\}$ and $i'_N:=\max\{i\in \N\colon k'_i\leq N\}$.
 Then the rectangles
 $(k_i,k_{i+1}]\times (k'_{i'},k'_{i'+1}]$, where  $i\in [i_N-1]$ and  $i'\in [i'_{N'}-1]$,
have the form  $(k, k+L_k]\times(k', k'+L'_{k'}]$, and together with a leftover set $E_{N,N'}$
form a
 partition  of  the rectangle $[N]\times[N']$.
 The  set $E_{N,N'}$ is  contained in  the union of the rectangles  $[N]\times (N'-L'_{k'_{i_{N'}}}, N']$ and
  $(N-L_{k_{i_N}},N]\times [N']$, and since $k_{i_N}\leq N$ and $k'_{i_N'}\leq N',$  we have
  $$
  |E_{N,N'}|\leq N\max_{k\leq N'}(L'_{k})+N'\max_{k\leq N}(L_k).
  $$
 Since  $L_N/N, L'_N/N\to 0$ as $N\to +\infty$, we get that
$|E_{N,N'}|/(NN')\to 0$ as $N,N'\to +\infty$.
 Using this,   the fact that $a\colon \N^2\to \C$ is bounded, and our assumption \eqref{E:assum}, we  deduce  that
$$
\lim_{N,N'\to +\infty} \frac{1}{NN'}\sum_{\bn\in [N]\times [N']}a(\bn)=
\lim_{k,k'\to +\infty} \frac{1}{L_kL'_{k'}}\sum_{\bn\in (k, k+L_k]\times(k', k'+L'_{k'}]} a(\bn)=0.
$$
 This completes the proof.
\end{proof}

\begin{proof}[Proof of Proposition~\ref{prop:hardyequi}]
 We give the proof for $d=2$, the proof in the general case is completely similar.

Suppose that the nilsequence $\psi$ has the form
$$
\psi(n,n')=\Psi(\tau^n \tau'^{n'}e_X), \quad  n,n'\in \N,
$$
  for some
nilmanifold $X=G/\Gamma$, commuting elements $\tau,\tau'\in G$, and  function $\Psi\in C(X)$.  By a remark made
in Section~\ref{SS:nil},
we can assume that the group $G$ is connected and simply connected.
 Moreover,  we can assume that $\Psi$  is a Lipschitz function  with $\norm\Psi_{\text{Lip}(X)}\leq 1$.

 By the infinitary factorization theorem  \cite[Corollary~1.12]{GT12c} (the same argument works for sequences in several variables) the sequence $g\colon \N^2\to G$ defined by $g(n,n'):=\tau^n\tau'^{n'}$, $n,n'\in \N,$ can be factorized as follows
$$
g(n,n')=g'(n,n')\, \gamma(n,n'), \quad n,n'\in \N,
$$
where
\begin{itemize}
\item $g'\colon \N^2\to G'$ is a polynomial sequence on a closed, connected and simply connected subgroup $G'$  of $G$
 such that $X':=G'/(G'\cap \Gamma)$ is a nilmanifold;

 \item the sequence $(g'(n,n')\cdot e_{X'})_{n,n'\in \N}$ is totally  equidistributed
on $X'$;

 \item the sequence $(\gamma(n,n')\cdot e_X)_{n,n'\in \N}$ is periodic and $\gamma(n,n')$ is a rational
element of $G$ for every $n,n'\in \N$.
\end{itemize}
Then for some $r\in \N$ and all
$(i_1,i_2)\in\{0,\ldots, r-1\}^2$ the sequence
$(\gamma(n,n')\cdot e_X)_{n,n'\in \N}$ is constant in the set $r\N^2+(i_1,i_2)$; say that it is equal to $\gamma_{i_1,i_2}\cdot e_X$ for some rational element $\gamma_{i_1,i_2}$ in $G$. After partitioning $\N^2$ as a union of such sets,
we are reduced to showing that
$$
\lim_{N, N'\to+\infty}\frac 1{N N'}
\sum_{n \in [N], n'\in  [N']}
e(f(rn+i_1)+f'(rn'+i_2))\, \Psi(g'(rn+i_1,rn'+i_2)\gamma_{i_1,i_2} \cdot e_X)=0
$$
for all $(i_1,i_2)\in\{0,\ldots, r-1\}^2$.
 Notice that if $h$ is a Hardy field function of at most polynomial growth that stays away from polynomials, then also $t\mapsto h(kt+l)$ has the same property for all  $k\in \N$ and   $l\in\Z$.
  Hence, it suffices to show that if $f, f'$ satisfy the assumptions of Proposition~\ref{prop:hardyequi}, then
\begin{equation}
\label{E:wanted}
\lim_{N, N'\to+\infty}\frac 1{N N'}
\sum_{n \in [N], n'\in  [N']}
e(f(n)+f'(n'))\, \Psi(g'(n,n')\gamma \cdot e_X)=0,
\end{equation}
where   $g'\colon \N^2\to G'$ is such that the infinite polynomial sequence $(g'(n,n')\cdot e_{X})_{n,n'\in \N}$ is equidistributed on $X'$ and $\gamma$ is a rational element of $G$.

Let $(L_N)_{N\in \N}$, $(L'_N)_{N\in \N}$ be  sequences of positive integers that will be specified later; for the moment we only assume that $L_N,L'_N\to +\infty$  and
$L_N/N, L'_N/N\to 0$ (we will also impose condition \eqref{E:qNa} later).
Using  Lemma~\ref{L:averaging} we see  that it suffices to show that
\begin{equation}
\label{eq:f2}
\lim_{N, N'\to+\infty}\frac 1{L_N L_{N'}}
\sum_{n \in [L_N], n'\in  [L_{N'}]}
e(f(N+n)+f(N'+n'))\, \Psi_\gamma(g_\gamma'(N+n,N'+n') \cdot e_X)=0,
\end{equation}
where  $g'_\gamma:=\gamma^{-1}g'\gamma$ is a polynomial sequence on $G'_\gamma:=\gamma^{-1}G\gamma$  and $\Psi_\gamma\in \text{Lip}(X)$ is defined by $\Psi_\gamma(g\cdot e_X):=\Psi(\gamma g\cdot e_X)$ for $g\in G$.

 In order to prove \eqref{eq:f2} we need to first gather some data. First, we claim that  for $N,N'\in \N$  the finite sequence
 $$
(g'(N+n,N'+n')\cdot e_X)_{n \in [L_N], n'\in  [L_{N'}]}
$$
  is $\delta_{N,N'}$-equidistributed on $X'$ for some $\delta_{N,N'}>0$ that satisfy $\delta_{N,N'}\to 0$ as $N,N'\to +\infty$. Indeed, if this is not the case, then there exists $\delta>0$,  $N_m,N'_m\to +\infty$, and $\Psi_m\in \text{Lip}(X')$ with $\norm{\Psi_m}_{\text{Lip}(X')}\leq 1$ and $\int_{X'}\Psi_m\,dm_{X'}=0$, such that
 \begin{equation}\label{E:Nmm}
\Big|\frac{1}{L_{N_m}L_{N_m'}}\sum_{(n,n')\in  (N_m+[L_{N_m}])\times (N_m'+[L_{{N_m'}}])} 
\Psi_m(g(n,n')\cdot e_{X'})\Big|\geq\delta, \quad \text{for all } m\in \N.
\end{equation}
By the Arzel\'a-Ascoli theorem, a subsequence of    $\Psi_m$ converges uniformly to some
$\Psi_0\in \text{Lip}(X')$ with $\norm{\Psi_0}_{\text{Lip}(X')}\leq 1$ and $\int_{X'}\Psi_0\,dm_{X'}=0$.
Then \eqref{E:Nmm} is satisfied
for infinitely many $m\in \N$  with $\Psi_0$ in place of $\Psi_m$ and $\delta/2$ in place of $\delta$.
Since $L_{N_m}, L_{N'_m}\to +\infty$,  $\big((N_m+[L_{N_m}])\times (N_m'+[L_{{N_m'}}])\big)_{m\in \N}$ is a F\o lner sequence in $\N^2$, and
we deduce that
$$
\limav  \Psi_0(g'(n,n')\cdot e_{X})\neq 0.
$$
 This contradicts our assumption that $(g'(n,n')\cdot e_{X})_{n,n'\in \N}$ is equidistributed in $X'$.

From the aforementioned equidistribution property, we deduce using  \cite[Corollary~5.5]{FH15a}
that the finite sequence
$$
(g_\gamma'(N+n,N'+n')\cdot  e_X)_{n \in [L_N], n'\in  [L_{N'}]}
$$ is $\delta'_{N,N'}$-equidistributed, where   $\delta'_{N,N'}\to 0$ as $N,N'\to +\infty$, on the nilmanifold $X'_\gamma:=G'_\gamma/\Gamma'_\gamma$,
where $\Gamma'_\gamma:=\Gamma\cap G'_\gamma$.

 We move now to the proof of \eqref{eq:f2}. We apply  Lemma~\ref{L:basichardy} for the functions $f,f'$ and we get sequences
 $(L_N)_{N\in \N}$, $(L'_N)_{N\in \N}$ satisfying the conditions in the lemma
 for some  $k,k'\in \N$ such that  $k,k'>\text{deg}(g_\gamma')$. After  ignoring negligible errors,
 we deduce that in \eqref{eq:f2} we can replace
  the finite sequences $(f(N+n))_{n\in [L_N]}$ and $(f'(N'+n'))_{n'\in [L'_{N'}]}$
  by finite
   polynomial sequences $(p_N(n))_{n\in [L_N]}$ and $({p'}_{N'}(n'))_{n'\in [L'_{N'}]}$, where
\begin{equation}\label{E:aN}
p_N(n)=n^{k}\alpha_N+q_N(n),\quad  p'_{N'}(n')=n'^{k'}\alpha'_{N'}+q'_{N'}(n'), \quad n,n',N,N'\in \N,
\end{equation}
for some real numbers $\alpha_N, \alpha'_N$ satisfying  $\alpha_N\to 0$ and $\alpha'_N\to 0,$ and    polynomials $q_N, q'_N\in \R[t]$, $N\in \N$,  that satisfy
\begin{equation}\label{E:qNa}
\text{deg}(q_N)<k,\  \text{deg}(q'_N)< k', \quad \text{and} \quad L_N^{k}|\alpha_N|, {L'}_N^{k'}|\alpha'_N|\to +\infty.
\end{equation}

We claim that the finite  polynomial sequence
\begin{equation}\label{E:2pol}
(p_N(n),p_{N'}(n'),  g_\gamma'(N+n,N'+n') \cdot e_X)_{n \in [L_N], n'\in  [L'_{N'}]}
\end{equation} is $\delta''_{N,N'}$-equidistributed in the nilmanifold $\T^2\times X'_\gamma$
where
$\delta''_{N,N'}\to 0$ as $N,N'\to +\infty$.
Arguing by contradiction, suppose that  this is not true. Then there exists $\delta>0$
such that
\begin{equation}\label{E:mm'}
\text{the sequence  } \eqref{E:2pol} \text{ is not } \delta\text{-equidistributed in }
\T^2\times X'_\gamma \text{ for some } N_m,N_m'\to +\infty.
\end{equation}
 The horizontal torus of the nilmanifold $X_\gamma'$ has the form $\T^s$ for some $s\in \N$. Then the horizontal torus of  the nilmanifold $\T^2\times X_\gamma'$ is $\T^2\times \T^s$. Let $\pi\colon G\to \T^s$ be the natural projection on the horizontal torus of $X_\gamma'$ and  $r_{N,N'}\colon \N^2\to \T^s$ be the polynomial sequence defined by $ r_{N,N'}(n,n'):=\pi(g_\gamma'(N+n,N'+n'))$.
Then
\begin{equation}\label{E:rNN'}
\text{deg}(r_{N,N'})<\min(k,k') \ \text{ for every } N,N'\in \N.
\end{equation}

By Theorem~\ref{th:Leibman}, we deduce that
there exists $M>0$ and $k_{m},k'_{m}\in \Z, k''_{m}\in \Z^s$ such that
for those $N_m, N'_m$ for which  \eqref{E:mm'} holds  and are greater than $M$ we have
\begin{equation}\label{E:km1}
0<|k_{m}| +|k'_{m}|+\norm{k''_{m}}_1\leq M
\end{equation}
and
\begin{equation}\label{E:Nm}
 \norm{k_{m}p_{N_m}(n)+k'_{m}p'_{N_m'}(n')+k''_{m}\cdot r_{N_m,N_m'}(n,n')}_{C^\infty([L_{N_m}]\times [L'_{N'_m}])}\leq M.
\end{equation}

 If $k_{m}=k'_{m}=0$ for infinitely many $m\in \N$, then $k''_{m}$ is non-zero for infinitely many $m\in \N$, and using \cite[Lemma~5.3]{FH15a} we get a contradiction from the fact that
the sequence  $(g_\gamma'(N+n,N'+n')\cdot  e_X)_{n \in [L_N], n'\in  [L'_{N'}]}
$ is $\delta'_{N,N'}$-equidistributed on the nilmanifold $X'_\gamma$ where   $\delta'_{N,N'}\to 0$ as $N,N'\to +\infty$.

 Suppose next that $k_{m}\neq 0$ for infinitely many $m\in \N$.  Using \eqref{E:Nm} (note that the polynomials $p_N$ and $p_{N'}$ depend on different variables) 
 in conjunction  with \eqref{E:aN}, \eqref{E:qNa}, \eqref{E:rNN'}, we obtain that
$$
L_{N_m}^k \norm{k_{m}\alpha_{N_m}}\leq M, \quad \text{for infinitely many } m\in \N.
$$
Since $\alpha_N\to 0$ and $1\leq |k_{m}|\leq M$,   we get  that
$\norm{k_{m}\alpha_{N_m}}=|k_{m}\alpha_{N_m}|\geq |\alpha_{N_m}|$ for infinitely many $m\in \N$. We deduce that
$$
L_{N_m}^k |\alpha_{N_m}|\leq M, \quad \text{for infinitely many } m\in \N.
$$
This contradicts \eqref{E:qNa}. The argument is similar if $k'_{m}\neq 0$ for infinitely many $m\in \N$.




Hence, the finite polynomial sequence \eqref{E:2pol}
is $\delta''_{N,N'}$-equidistributed in the nilmanifold $\T^2\times X'_\gamma$
where
$\delta''_{N,N'}\to 0$ as $N,N'\to +\infty$.
We deduce
that the limit in the left hand side of \eqref{eq:f2} is equal to
$$
\int e(t)\cdot e(t')\cdot F_\gamma(x)\, dm_{\T^2\times X'_\gamma}=0
$$
where the last identity follows since  $m_{\T^2\times X'_\gamma}=m_{\T^2}\times m_{X'_\gamma}$.  This completes the proof.
\end{proof}

\begin{proof}[Proof of Corollary~\ref{C:hardymain}]
Let $f$ be a Hardy field function of polynomial growth. We consider the following  three cases:

If $f$ stays away from polynomials, then  the conclusion follows from the $d=1$ case of Theorem~\ref{T:hardymain}
and the corresponding averages converge to $0$ in $L^2(\mu)$.

 Suppose next that $f$ is asymptotically polynomial, that is,  $f(t)-p(t)\to 0$ for some $p\in \R[t]$. In this case, the mean convergence of the  averages  \eqref{eq:Cesaro-av4}  follows from
Proposition~\ref{prop:weight-nil} and the well known fact that sequences of the form $n\mapsto e(p(n))$ are  nilsequences.

Lastly, suppose that  $f=p+g$ for some polynomial $p\in \R[t]$ and Hardy field function $g$ that satisfies $|g(t)|\to +\infty$ and $|g(t)|\leq C \log t $  for some $C>0$ and all sufficiently large  $t\in \R_+$. Let $p(t)=\sum_{i=0}^\ell\alpha_i t^i, t\in \R,$ for some $\ell\in\N$ and $\alpha_1,\ldots, \alpha_\ell\in \R$.  For $i=1,\ldots, \ell$, we consider the commuting transformations
$T_it:=t+\alpha_i, t\in \T,$ acting on $\T$ with the Haar measure $m_\T$, and the function $h\in L^\infty(m_\T)$ defined by $h(t):=e(-t), t\in \T$. Then
$$
\frac{1}{N}\sum_{n=1}^N e(f(n)) \, h\big(\prod_{i=1}^\ell T_i^{n^i}t\big)=e(-t+\alpha_0)\, \frac{1}{N}\sum_{n=1}^N e(g(n)), \quad \text{for every } N\in \N, \, t\in \T.
$$
By \cite[Proof of Theorem~3.1]{F09}  (see also \cite[Proof of Theorem~3.3]{BKQW94}) we get that the last averages do not converge as $N\to +\infty$. Hence, the sequence $n\mapsto e(f(n))$ is not a good universal weight for weak convergence of averages of the form  \eqref{eq:Cesaro-av4} even when $s=1$.

If $f$ is any Hardy field function and $p\in \R[t]$ is any polynomial, then  it is known that the limit $\lim_{t\to +\infty} (f(t)-p(t))/\log t$ either exists or else is $\pm \infty$.
   Hence, every Hardy field function with at most polynomial growth is covered in one of the previous three cases and  the proof is complete.
\end{proof}

\appendix

\section{Seminorms on $L^\infty(\mu)$ and related factors}
\label{sec:seminorms}
Let $(X,\mu,T_1,\dots,T_\ell)$ be a system.
We recall here the definition and some properties of the seminorms $\nnorm\cdot_k$ on $L^\infty(\mu)$ and of the factors $\CZ_k$ defined in~\cite{HK05} for the ergodic case  and in~\cite{CFH11} for the general case. These two papers deal only with the case of a single transformation, the generalization to the case  of several commuting transformations
 is completely similar and is given below.

\subsection{The seminorms $\nnorm\cdot_k$}
We write $\CI(\vec T)$ for the $\sigma$-algebra of sets invariant under all transformations $T_1,\dots,T_\ell$. For $f\in L^\infty(\mu)$, we define
\begin{equation}
\label{eq:seminorm1}
\nnorm f_1:=\bigl\Vert\E_\mu(f\mid\CI(\vec T))\bigr\Vert_{L^2(\mu)}
\end{equation}
and for $k\in \N$ we let
\begin{equation}
\label{eq:seminorm2}
\nnorm f_{k+1}:=\Bigl(\mathrm{Lim\, Av}_{\vec n}\ \nnorm{f\cdot T_{\vec n}\overline f\,}_k^{2^k}\Bigr)^{1/2^{k+1}}
\end{equation}
where, as usual, we use the notation $T_{\vec n}=\prod_{i=1}^\ell T_i^{n_i}$ for $\vec n=(n_1,\ldots, n_\ell)$. By induction, we have
$$\nnorm f_k\leq\norm f_{L^\infty(\mu)}\ \text{ for every }k\in \N.
$$
In case of ambiguity, we write $\nnorm f_{\mu,k}$ or $\nnorm f_{\vec T,\mu,k}$.
If $\mu=\int\mu_x\,d\mu(x)$ is the ergodic decomposition of $\mu$ under $\vec T$, then for every $f\in L^\infty(\mu)$ and every $k\in \N$ we have
\begin{equation}\label{E:seminonerg}
\nnorm f_{\mu,k}^{2^k}=\int\nnorm f_{\mu_x,k}^{2^k}\,d\mu(x).
\end{equation}

For $f\in L^\infty(\mu)$, by~\eqref{eq:seminorm1} we have
\begin{equation}
\label{eq:seminorm3}
\Bigl|\int f\,d\mu\Bigr|\leq\nnorm f_1.
\end{equation}
Writing $\vec T\times\vec T$ for the $\Z^\ell$-action on $X\times X$ induced  by
$T_1\times T_1$, \dots, $T_\ell\times T_\ell$, we have
\begin{align*}
\nnorm{f\otimes\overline f}_{\vec T\times \vec T,\mu\times\mu,1}^2
&=\bigl\Vert\E_{\mu\times\mu}\bigl(f\otimes\overline f\mid\CI(\vec T\times \vec T)\bigr)\bigr\Vert_{L^2(\mu\times\mu)}^2
\\
&=\Limav{\vec n}\,\Bigl|\int f\cdot T_{\vec n}\overline f\,d\mu\Bigr|^2\text{ by the ergodic theorem}\\
&\leq\Limsup\,\Av_{\vec n}\,\bigl\Vert\E_\mu\bigl(f\cdot T_{\vec n}\overline f\mid\CI(\vec T)\bigr)\bigr\Vert_{L^2(\mu)}^2\\
&=\Limsup\,\Av_{\vec n}\,\nnorm{f\cdot T_{\vec n}\overline f}_{\vec T,\mu,1}^2
=\nnorm f_{\vec T,\mu,2}^4 \ \ \text{by~\eqref{eq:seminorm1} and~\eqref{eq:seminorm2}.}\\
\end{align*}
By induction, using the relation~\eqref{eq:seminorm2} for the measures $\mu\times\mu$ and $\mu$, we deduce that for every $k\in \N$ we have
\begin{equation}
\label{eq:seminorm4}
\nnorm{f\otimes\overline f}_{\vec T\times \vec T,\mu\times\mu,k}\leq
\nnorm f_{\vec T,\mu,k+1}^2.
\end{equation}

\subsection{The factors $\CZ_k$}
\label{subse:Zk}
For $k\in \Z_+$, the factor $\CZ_k$ of $X$ is characterized by the following property
$$
\text{\em for } f\in L^\infty(\mu),\ \E_\mu(f|\CZ_k)=0\text{ \em if and
  only if }\  \nnorm f_{k+1} = 0.
$$
Equivalently, one has
$$
L^\infty(\CZ_{k},\mu)
=\Big\{f\in L^\infty(\mu)\colon \int f\cdot g \ d\mu=0
 \text{ for every } g\in L^\infty(\mu) \text{ with }
\nnorm g_{k+1} = 0\Big\}.
$$
In case of ambiguity, we write $\CZ_k(X,\mu,\vec T)$.

We say that $(X,\mu,\vec T)$ is a \emph{system of order $k$} if the $\sigma$-algebra $\CZ_k$ coincides with the $\sigma$-algebra $\CX$.
Equivalently, if $\nnorm\cdot_{k+1}$ is a norm on $L^\infty(\mu)$.

\begin{proposition}
\label{prop:ae-rotation}
Let $(X,\mu,\vec T)$ be a  system of order $1$ and
$\mu=\int\mu_x\,d\mu(x)$ be the ergodic decomposition of $\mu$ under $\vec T$. Then for $\mu$-a.e. $x\in X$ the system $(X,\mu_x,\vec T)$ is isomorphic to an ergodic rotation
on a compact Abelian group.
\end{proposition}
\begin{remark}
It is not hard to show that the converse also holds.
\end{remark}
\begin{proof}
Since $(X,\CX, \mu)$ is a Lebesgue space, there exists a countable sequence
$(f_n)_{n\in\N}$ of bounded Borel functions (defined everywhere) that is dense in $L^1(\mu)$ and in $L^1(\mu_x)$ for every $x\in X$. By~\cite[Corollary~3.3]{CFH11}, there exists a Borel set $X_1\subset X$  with $\mu(X_1)=1$, such that for every $x\in X_1$ and every $n\in \N$ the function $f_n$ belongs to $L^\infty(\CZ_1(X,\mu_x,\vec T))$. For $x\in X_1$ it follows by density that every $f\in L^1(\mu_x)$ belongs to $L^1(X,\CZ_1(X,\mu_x,\vec T))$. The $\sigma$-algebras $\CX$ and $\CZ_1(X,\mu_x,\vec T)$ coincide up to $\mu_x$-null sets and $(X,\mu_x,\vec T)$ is a system of order $1$ for $\mu$-a.e. $x\in X$. It is well known that  an ergodic system of order $1$ is isomorphic to an  ergodic rotation on  a compact Abelian group and the proof is complete.
\end{proof}




\end{document}